\documentclass{amsart}
\usepackage{amssymb, latexsym, amsmath}
\usepackage{bbm}
\usepackage{longtable}
\usepackage{lscape}

\newcommand{\C}{\mathbbm{C}}
\newcommand{\R}{\mathbbm{R}}
\newcommand{\N}{\mathbbm{N}}
\newcommand{\Z}{\mathbbm{Z}}
\newcommand{\Q}{\mathbbm{Q}}
\newcommand{\A}{\mathbb{A}}
\newcommand{\F}{\mathbbm{F}}
\newcommand{\GSp}{{\rm GSp}}
\newcommand{\PGSp}{{\rm PGSp}}
\newcommand{\GL}{{\rm GL}}
\newcommand{\Sp}{{\rm Sp}}
\newcommand{\Orth}{{\rm O}}
\newcommand{\SOrth}{{\rm SO}}
\newcommand{\St}{{\rm St}}
\newcommand{\triv}{{\bf 1}}

\theoremstyle{plain}
\newtheorem{thm}{Theorem}[section]
\newtheorem{prop}[thm]{Proposition}
\newtheorem{lemma}[thm]{Lemma}
\newtheorem{cor}[thm]{Corollary}

\theoremstyle{definition}

\theoremstyle{remark}
\newtheorem{remark}[thm]{Remark}

\usepackage{fancyhdr}

\usepackage{paralist}
\usepackage{hyperref}

\begin{document}

\title[Irreducible characters of $\GSp(4, q)$]{Irreducible characters of $\GSp(4, q)$ and dimensions of spaces of fixed vectors}

\thanks{The author would like to thank Ralf Schmidt and Alan Roche for their valuable notes and comments on this paper. The author would especially like to thank the referee for their detailed comments and careful reading of this paper.}

\author{Jeffery Breeding II}
\address{Department of Mathematics, Fordham University, Bronx, NY 10458}
\email{jbreeding@fordham.edu}
\urladdr{math.jeffbreeding.com}


\maketitle

\begin{abstract}
In this paper, we compute the conjugacy classes and the list of irreducible characters of $\GSp(4,q)$, where $q$ is odd. We also determine precisely which irreducible characters are non-cuspidal and which are generic. These characters are then used to compute 
dimensions of certain subspaces of fixed vectors of smooth admissible non-supercuspidal representations of $\GSp(4,F)$, where $F$ is a non-archimedean local field of characteristic zero with residue field of order $q$.\\
\end{abstract}

\tableofcontents


\section{Introduction}
\label{intro}
The passage of cuspidal Siegel modular forms of degree $n$ to cuspidal automorphic representations of $\GSp(2n,\mathbb{A})$ is a natural extension of the classical theory for GL(2) described in the seminal work \cite{JL} by Jacquet and Langlands. As in the GL(2)-case, these cuspidal automorphic representations $\pi$ can be written in terms of local components $\pi_v$, where $v$ is a place of $\Q$.

In the degree 2 case, let $\Gamma'$ be a discrete subgroup of $\Sp(4,\R)$. Dimensions of the spaces $S^2_k(\Gamma')$ can be computed using the Selberg trace formula. A general formula for these dimensions has been given by Hashimoto \cite{Hash}. In terms of the associated representations, these dimensions tell us essentially how many possibilities there are for the local factors of the representation by giving restrictions on the dimensions of its subspace of $\Gamma'$-fixed vectors.

Let $N\in\N$ be square-free. In the case $\Gamma'=\Gamma(N)$, the principal congruence subgroup of level $N$, the dimensions of spaces of $\Gamma(N)$-fixed vectors can be computed in terms of representations of the finite group $\GSp(4,p)$ for primes $p|N$. Let $\pi$ be an irreducible automorphic representation of $\GSp(4,\A)$, where $\A$ is the adele ring of $\Q$. Write $\pi$ in terms of its local components $$\pi=\otimes_{p\leq\infty}'\pi_p,$$ where $\pi_p$ is a representation of $\GSp(4,\Q_p)$. Suppose $\pi$ is associated to a Siegel cuspidal eigenform $f\in S_k^2(\Gamma(N))$. Then, for each $p|N$, the representation $\pi_p$ has a finite dimensional space of $\Gamma(p)$-fixed vectors. This dimension would make a contribution to the dimension of the space $S^2_k(\Gamma(N))$. The computation of the precise contribution can be achieved by determining the `$\GSp(4,p)$ analogue' of the representation $\pi_p$ of $\GSp(4,\Q_p)$ for each finite prime $p$. The dimension of the $\GSp(4,p)$ analogue is the dimension of the subspace of $\Gamma(p)$-fixed vectors.

The conjugacy classes and the complete list (Table \ref{GSpIrredCharacters}) of all of the irreducible characters of the finite group $\GSp(4,q)$ as well as their cuspidality and genericity is determined in this paper. Our method relies on Srinivasan's irreducible character tables for $\Sp(4,q)$, \cite{Sr}, and does not lead to the computation of the complete character table\footnote{The referee has informed us that Shinoda had independently obtained the conjugacy classes \cite{ShinodaC} and the complete character table \cite{Shinoda} of $\GSp(4,q)$, which he calls ${\rm CSp}(4,q),$ using the modern approach of Deligne-Lusztig theory. Shinoda remarks at the end of the introduction in \cite{Shinoda} that Reid had also obtained this character table in \cite{Reid}. We were not aware of these results when doing the computations in this paper. We give a correspondence between our notations and Shinoda's notations for the conjugacy classes and the irreducible characters of $\GSp(4,q)$ in section \ref{conjugacyclasses} and section \ref{Irredchar}, respectively.}. In particular, our method will not yield the character values of some of the irreducible characters on conjugacy classes that have non-square multipliers.

Cuspidality is determined by defining cuspidal representations on the (standard) Borel, Siegel parabolic, and Klingen parabolic subgroups and then inducing. The irreducible non-cuspidal representations are precisely the irreducible constituents of these induced representations. Criteria are determined for these induced characters to be irreducible. If an induced character is reducible, then its irreducible constituents are determined. We find that we did not need the complete character table to identify the cuspidal characters.

After computing what we need from the representation theory of $\GSp(4,q)$, we give formulas for dimensions of subspaces of $\Gamma(\mathfrak{p})$-fixed vectors of certain smooth admissible non-supercuspidal representations of $\GSp(4,F)$, where $F$ is a non-archimed-
ean local field of characteristic zero with ring of integers $\mathfrak{o}$ and maximal ideal $\mathfrak{p}$ such that $\mathfrak{o}/\mathfrak{p}$ is a finite field with $q$ elements.

\section{Definitions and notations}
\label{Basicdefs}
Let $G$ be a finite group of Lie type. A {\em Borel subgroup} of $G$ is a maximal closed connected solvable subgroup of $G$. For any such group $G$, any two Borel subgroups of $G$ are conjugate. A {\em parabolic subgroup} of $G$ is a subgroup of $G$ which contains a Borel subgroup of $G$.

Let $B$ be a Borel subgroup of $G$. Then $B$ has the Levi decomposition $B = T\ltimes U$, where $T$ is a maximal torus in $G$ such that the centralizer of $T$ in $G$ is contained in $B$ and $U$ is a connected unipotent group. Moreover, all maximal tori of $B$ are conjugate in $B$. Let $T\subset B$ be a maximal torus and let $W=N(T)/T$ be its Weyl group, where $N(T)$ is the normalizer of $T$ in $G$. Then $G$ has the following double coset decomposition, which is called its {\em Bruhat decomposition} with respect to $B$:
$$G=\bigsqcup_{w\in W} BwB.$$

Let $\F_q$ denote a finite field with $q=p^n$ elements, with $p$ an odd prime. Let $\kappa$ be a generator of $\F_{q^4}^\times$ and let $\zeta = \kappa^{q^2-1}$, $\theta = \kappa^{q^2+1}$, $\eta = \theta^{q-1}$, and $\gamma = \theta^{q+1}$. The element $\eta$ is the generator of the set of elements in $\F_{q^2}^\times$ whose norm over $\F_q$ is 1 and $\zeta$ is the generator of the set of elements in $\F_{q^4}^\times$ whose norm over $\F_{q^2}$ is 1. Fix a monomorphism from $\F_{q^4}^\times$ into $\C^\times$ and denote the images of $\zeta,\theta,\eta,$ and $\gamma$ under this monomorphism by $\tilde{\zeta},\tilde{\theta},\tilde{\eta},$ and $\tilde{\gamma}$, respectively. Define the sets
\begin{flushleft}
$R_1 = \{1,...,\frac{1}{4}(q^2-1)\}$, \\
$R_2$ is a set of $\frac{1}{4}(q-1)^2$ distinct positive integers $i$ such that $\theta^i$, $\theta^{-i}$, $\theta^{qi}$, and $\theta^{-qi}$ are all distinct, \\
$T_1 = \{1,...,\frac{1}{2}(q-3)\}$, \\
$T_2 = \{1,...,\frac{1}{2}(q-1)\}$, and \\
$T_3 = \{1,...,q-1\}$.
\end{flushleft}

The {\em general symplectic group of degree two over $\F_q$}, denoted $\GSp(4,q)$, is the set of all $g\in\GL(4,q)$ such that
$${}^tgJg=\lambda(g) J,\quad {\rm where}\quad  J=
\left(\begin{smallmatrix}
& & & 1\\
& & 1 & \\
& -1 & & \\
-1 & & & \\
\end{smallmatrix}\right)$$
and $\lambda(g)$ is an element of $\F_q^\times$ which depends on $g$. We will call $\lambda(g)$ the \emph{multiplier} of $g$. The set of all $g\in \GSp(4,q)$ such that $\lambda(g)=1$ is the subgroup $\Sp(4,q)$. Let $Z=\{\gamma^k\cdot I_4:k\in T_3\}$ denote the center of $\GSp(4,q)$.

For any $g\in G$, there exists a unique $\lambda(g)$ and $g'\in\Sp(4,q)$ such that
\begin{displaymath}
g=\left(\begin{smallmatrix}
1 & & & \\
& 1 & & \\
& & \lambda(g) & \\
& & & \lambda(g) \\
\end{smallmatrix}\right)\,\cdot\, g'.
\end{displaymath}

\vspace{0.05in}
The order of $\Sp(4,q)$, as computed by Wall in \cite{Wa}, is $q^4(q^4-1)(q^2-1)$. So the order of $\GSp(4,q)$ is $q^4(q^4-1)(q^2-1)(q-1)$.

The {\em (standard) Borel subgroup} $B$ of $\GSp(4,q)$ is the set of all of the upper triangular matrices,
\begin{displaymath}
B = \left\{ \left(\begin{smallmatrix}
* & * & * & * \\
& * & * & * \\
& & * & * \\
& & & * \\
\end{smallmatrix}\right)\in \GSp(4,q) \right\}.
\end{displaymath}
Since $\GSp(4,q)$ is a finite group of Lie type, any other Borel subgroup of $\GSp(4,q)$ is conjugate to $B$. Using the Levi decomposition, every element $g\in B$ can be written uniquely as
\begin{equation}\label{Bdecomp}
g = \left(\begin{smallmatrix}
a & & & \\
& b & & \\
& & cb^{-1} & \\
& & & ca^{-1} \\
\end{smallmatrix}\right)\cdot
\left(\begin{smallmatrix}
1 & & & \\
& 1 & x & \\
& & 1 & \\
& & & 1 \\
\end{smallmatrix}\right)\cdot
\left(\begin{smallmatrix}
1 & \lambda & \mu & y \\
& 1 & & \mu \\
& & 1 & -\lambda \\
& & & 1 \\
\end{smallmatrix}\right),
\end{equation}
with $a, b, c\in \F_q^\times$ and $x, y, \lambda, \mu \in \F_q$. The first factor in (\ref{Bdecomp}) is an element of the split torus of $\GSp(4,q)$. The subgroup of $B$ of elements which have $1$ on every entry on the main diagonal is the unipotent radical of the Borel and will be denoted by $N_{\rm GSp(4)}$. The order of $B$ is $q^4(q-1)^3$. The multiplier of the element $g$ in (\ref{Bdecomp}) is $c$. 

The {\em (standard) Siegel parabolic subgroup} $P$ of $\GSp(4,q)$ is 
\begin{displaymath}
P = \left\{ \left(\begin{smallmatrix}
* & * & * & * \\
* & * & * & * \\
& & * & * \\
& & * & * \\
\end{smallmatrix}\right)\in \GSp(4,q) \right\}.
\end{displaymath}
Any other Siegel parabolic subgroup of $\GSp(4,q)$ is conjugate to $P$. Using the Levi decomposition, every element $p\in P$ can be written uniquely as
\begin{equation}\label{Pdecomp}
p = \left(\begin{smallmatrix}
a & b & & \\
c & d & & \\
& & {\lambda}a/{\Delta} & -{\lambda}b/{\Delta} \\
& & -{\lambda}c/{\Delta} & {\lambda}d/{\Delta} \\
\end{smallmatrix}\right)\cdot
\left(\begin{smallmatrix}
1 & & \mu & y \\
& 1 & x & \mu \\
& & 1 & \\
& & & 1 \\
\end{smallmatrix}\right),
\end{equation}
with $\Delta = ad - bc \in \F_q^\times$, $\lambda \in \F_q^\times$ and $x, y, \mu \in \F_q$. The first factor in (\ref{Pdecomp}) is contained in a Levi subgroup of $\GSp(4,q)$ and the second factor is in the unipotent radical of $P$. The order of $P$ is $q^4(q^2-1)(q-1)^2$. The multiplier of the element $p$ in (\ref{Pdecomp}) is $\lambda$. We also define
\begin{displaymath}
A'=\left(\begin{smallmatrix}
& 1 \\
1 & \\
\end{smallmatrix}\right)
{}^tA^{-1}
\left(\begin{smallmatrix}
& 1 \\
1 & \\
\end{smallmatrix}\right)
\end{displaymath}
for any $A\in \GL(2, q)$. Then the Levi subgroup of $P$ is the set of all matrices of the form
$$\left(\begin{smallmatrix}
A & \\
& \lambda\cdot A'\\
\end{smallmatrix}\right)$$
where $A\in\GL(2,q)$ and $\lambda\in\F_q^\times$.

The {\em (standard) Klingen parabolic subgroup $Q$ of $\GSp(4,q)$} is
\begin{displaymath}
Q = \left\{ \left(\begin{smallmatrix}
* & * & * & * \\
& * & * & * \\
& * & * & * \\
& & & * \\
\end{smallmatrix}\right)\in \GSp(4,q) \right\}.
\end{displaymath}
Any other Klingen parabolic subgroup of $\GSp(4,q)$ is conjugate to $Q$. Using the Levi decomposition, every element $g\in Q$ can be written uniquely as
\begin{equation}\label{Qdecomp}
g = \left(\begin{smallmatrix}
t & & & \\
& a & b & \\
& c & d & \\
& & & \Delta t^{-1} \\
\end{smallmatrix}\right)\cdot
\left(\begin{smallmatrix}
1 & \lambda & \mu & y \\
& 1 & & \mu \\
& & 1 & -\lambda \\
& & & 1 \\
\end{smallmatrix}\right),
\end{equation}
with $\Delta = ad-bc \in \F_q^\times$, $t\in \F_q^\times,$ and $y, \lambda, \mu \in \F_q$. The first factor in (\ref{Qdecomp}) is contained in a Levi subgroup of $\GSp(4,q)$ and the second factor is in the unipotent radical of $Q$. The order of $Q$ is $q^4(q^2-1)(q-1)^2$. The multiplier of the element $g$ in (\ref{Qdecomp}) is $\Delta$.

We remind the reader of the following results in the representation theory of finite groups. The reader may wish to consult \cite{Serre} for more details.

Let $G$ be a finite group and let $\chi_1$ and $\chi_2$ be characters of $G$. The \emph{inner product} of $\chi_1$ and $\chi_2$ is
$$(\chi_1, \chi_2) = \frac{1}{|G|} \sum_{g\in G}\chi_1(g)\overline{\chi_2(g)}.$$
This inner product is equal to
$${\rm dim\, Hom}(\rho_1,\rho_2),$$
where $\rho_i$ is a representation with character $\chi_i$ for $i=1,2$.

Let $H$ be a subgroup of $G$ and let $(\pi,V)$ be a representation of $H$. The {\em induced representation} of $(\pi,V)$ on $G$, denoted by ${\rm Ind}_H^G V$ or by ${\rm Ind}_H^G \pi$, is the space of functions $f:G\rightarrow V$ satisfying
$$f(hg) = \pi(h)f(g), \text{for } h\in H,g\in G,$$
with group action by right translation. If $\pi$ has character $\chi$, we denote the character of ${\rm Ind}_H^G \pi$ by $\chi_H^G$ or by ${\rm Ind}_H^G(\chi)$ and call it the {\em induced character of $\pi$}.

Let $C$ be a conjugacy class of $G$. Then the conjugacy class $C$ either has an empty intersection with $H$ or it splits into finitely many distinct conjugacy classes of the subgroup $H$, say $C = D_1 \sqcup \ldots \sqcup D_r$. One may compute the values of the induced character $\chi^G_H$ on the conjugacy class $C$ using the following formula from \cite{Fulton}:
\begin{displaymath}
 \chi^G_H(C) = \frac{|G|}{|H|}\sum_{i=1}^r \frac{|D_i|}{|C|} \chi (D_i).
\end{displaymath}

We also recall some results from the representation theory of $\GL(2, q)$. These will be useful in our computations. More details can be found in \cite{Bump} and in \cite{Fulton}.

Let $B_{\GL(2)}$ be the {\em (standard) Borel subgroup} of $\GL(2, q)$, i.e., $B_{\GL(2)}$ is the subgroup of all upper triangular matrices,
\begin{displaymath}
B_{\GL(2)} = \left\{ \left(\begin{smallmatrix}
* & * \\
& * \\
\end{smallmatrix}\right) \in \GL(2, q) \right\}.
\end{displaymath}
By the Levi decomposition, every element in $B_{\GL(2)}$ can be written uniquely as
$$\left(\begin{smallmatrix}
y_1 &  \\
& y_2 \\
\end{smallmatrix}\right)\cdot\left(\begin{smallmatrix}
1 & a \\
& 1 \\
\end{smallmatrix}\right),$$
where the first factor is an element of the split torus of $\GL(2,q)$, denoted by $T_{\GL(2)}$. The subgroup of $B_{\GL(2)}$ of matrices which have 1 on every entry on the main diagonal is the unipotent radical of $B_{\GL(2)}$ and will be denoted by $N_{\GL(2)}$.

Let $\chi_1, \chi_2$ be characters of $\F_q^\times$. Define a representation $\chi$ of $B_{\GL(2)}$ by
\begin{displaymath}
\chi \left(\begin{smallmatrix}
y_1 & x \\
& y_2 \\
\end{smallmatrix}\right) = \chi_1 (y_1) \chi_2(y_2).
\end{displaymath}
Denote the representation of $\GL(2, q)$ induced from $\chi$ by $\chi_1 \times \chi_2$. The irreducible representations $\chi_1 \times \chi_2$ are called \emph{representations of the principal series}. The following two results can be found in \cite{Bump}.

\begin{thm} Let $\chi_1, \chi_2, \mu_1$ and $\mu_2$ be characters of $\F_q^\times$. Then $\chi_1 \times \chi_2$ is an irreducible representation of degree $q + 1$ of $\GL(2, q)$ unless $\chi_1 = \chi_2$, in which case it is the direct sum of two irreducible representations having degrees $1$ and $q$. We have
\begin{displaymath}
\chi_1 \times \chi_2 \cong \mu_1 \times \mu_2
\end{displaymath}
if and only if either (1) $\chi_1 = \mu_1 \, \, {\rm and} \, \, \chi_2 = \mu_2$ or (2) $\chi_1 = \mu_2 \, \, {\rm and} \, \, \chi_2 = \mu_1.$
\end{thm}

Let $\triv_{\F_q^\times}$ denote the trivial character of $\F_q^\times$. Then, $\triv_{\F_q^\times}\times\triv_{\F_q^\times}$ is not irreducible. It decomposes into the sum of two irreducible representations: $\triv_{\F_q^\times}\times\triv_{\F_q^\times}=\triv_{\GL(2)}+\St_{\GL(2)},$
where $\triv_{\GL(2)}$ is the trivial representation of $\GL(2,q)$ and $\St_{\GL(2)}$ is a $q$-dimensional representation of $\GL(2,q)$ called its {\em Steinberg representation}. Then, for any character $\chi$ of $\F_q^\times$, the irreducible representation of dimension 1 contained in $\chi \times \chi$ is the character $g \mapsto \chi({\rm det}(g))$. The other irreducible constituent is the $q$-dimensional representation which can be obtained by taking the tensor of $\chi$ with $\St_{\GL(2)}$.

\begin{prop}
Every irreducible cuspidal representation of $\GL(2, q)$ has dimension $q-1$. Every cuspidal representation is a direct sum of irreducible cuspidal representations.
\end{prop}

Let $\alpha, \beta$ be characters of $\F_q^\times$. Let $\phi$ be a character of $\F_{q^2}^\times$ such that $\phi \neq \phi^q$ and $X_\phi\cong X_{\phi^q}$ is an irreducible cuspidal representation of $\GL(2,q)$. Define $M(x):=\left(\begin{smallmatrix}x&1\\&x\end{smallmatrix}\right)$ and let $d=x+y\sqrt{\gamma},$ where $\sqrt{\gamma}=\gamma^{1/2}$ is a fixed element $t\in\F_{q^2}$ such that $t^2=\gamma$. The complete character table of $\GL(2, q)$, taken from \cite{Fulton}, is the table below, no matter what choice for $\sqrt{\gamma}$ was made.

\begin{center}
\begin{small}
\begin{tabular}{|c||c|c|c|c|}
\hline
Class & ${\rm diag}(x,x)$ & $M(x)$ & ${\rm diag}(x,y), x\neq y$ & $d\leftrightarrow\left(\begin{smallmatrix} x&y\gamma\\ y&x\end{smallmatrix}\right)$\\
Class order & 1 & $q^2-1$ & $q^2+q$ & $q^2-q$\\
\hline
\hline
$\alpha\triv_{\GL(2)}$ & $\alpha(x^2)$ & $\alpha(x^2)$ & $\alpha(xy)$ & $\alpha(d^{q+1})$\\
\hline
$\alpha\St_{\GL(2)}$ & $q\alpha(x^2)$ & 0 & $\alpha(xy)$ & $-\alpha(d^{q+1})$\\
\hline
$\alpha\times\beta, \alpha\neq\beta$ & $(q+1)\alpha(x)\beta(x)$ & $\alpha(x)\beta(x)$ & $\alpha(x)\beta(y)+\alpha(y)\beta(x)$ & 0\\
\hline
$X_\phi$ & $(q-1)\phi(x)$ & $-\phi(x)$ & 0 & $-\left(\phi(d) + \phi(d^q)\right)$\\
\hline
\end{tabular}
\end{small}
\end{center}

\section{Conjugacy classes}
\label{conjugacyclasses}
Our first step in the determination of the irreducible characters of $\GSp(4,q)$ is the computation of its conjugacy classes. This list of classes can then be used to compute the conjugacy classes of the (standard) Borel, Siegel parabolic, and Klingen parabolic subgroups. One does this by determining which conjugacy classes have a non-empty intersection with the subgroup, determining how each class splits, and computing the order of the centralizer of each conjugacy class of the subgroup.

The main tool we used for computing the conjugacy classes of $\GSp(4,q)$ is a paper of Wall, \cite{Wa}. This paper, incidentally, was also used by Srinivasan in \cite{Sr} to determine the conjugacy classes of the symplectic group $\Sp(4,q)$. Unfortunately, Wall's results cannot be immediately used in our case. We can, however, use these results to find the classes of $\SOrth(5,q)$. This is particularly useful because $\SOrth(5,q)$ is isomorphic to $\PGSp(4,q) := \GSp(4,q) / Z$. The reader can consult Appendix A.7 in \cite{RS} for a description of an isomorphism $\rho_5:\PGSp(4,F)\rightarrow\SOrth(5,F)$, where $F$ is a field of characteristic not equal to 2. The computation of the classes and their orders of this special orthogonal group are then used to determine the classes of $\GSp(4,q)$. One does this by computing the corresponding class of $\PGSp(4,q)$ under an isomorphism and then pulling the class back to $\GSp(4,q)$.

Consider the natural projection map from $\GSp(4,q)$ to $\PGSp(4,q)$ given by
\begin{displaymath}
\GSp(4,q) \longrightarrow \PGSp(4, q), \text{   } g \mapsto \overline{g}.
\end{displaymath}
Let $g,h\in \GSp(4, q)$. If $g=xhx^{-1}$ for some $x\in\GSp(4,q)$, then it is clear that $\lambda(g)=\lambda(h)$. Moreover, under the projection map, $$\overline{g}=\overline{xhx^{-1}}=\overline{x}\cdot\overline{h}\cdot \overline{x^{-1}}.$$
So if two elements are conjugate in $\GSp(4,q)$, they must be conjugate in $\PGSp(4,q)$. The list of class representatives in $\PGSp(4,q)$, when pulled back to $\GSp(4,q)$, hit class representatives of all the conjugacy classes of $\GSp(4,q)$. Suppose now that two elements $g,h\in \GSp(4, q)$ are conjugate in $\PGSp(4,q)$, i.e. $\overline{g}=\overline{x}\cdot\overline{h}\cdot\overline{x^{-1}}$, for some $\overline{x}\in \PGSp(4, q)$. Then, for some $\gamma^i\in \F_q^\times$,
\begin{displaymath}
g = \gamma^i\cdot I_4\cdot x h x^{-1}.
\end{displaymath}
Taking multipliers on both sides, we have $\lambda(g)=\gamma^{2i}\lambda(h)$. So the multiplier of $g$ is a square if and only if the multiplier of $h$ is a square. Write $g$ and $h$ as
\begin{displaymath}
g = \left(\begin{smallmatrix}
1 & & & \\
& 1 & & \\
& & \gamma^{i_g} & \\
& & & \gamma^{i_g} \\
\end{smallmatrix}\right)\cdot z_{j_g}\cdot g'\quad{\rm and}\quad
h = \left(\begin{smallmatrix}
1 & & & \\
& 1 & & \\
& & \gamma^{i_h} & \\
& & & \gamma^{i_h} \\
\end{smallmatrix}\right)\cdot z_{j_h}\cdot h',
\end{displaymath}
with $g',h'\in \Sp(4, F_q)$, $i_g,i_h\in\{0,1\}$, and $j_g,j_h\in T_3$. So the multipliers of $g$ and $h$ respectively are $\lambda(g)=\gamma^{i_g+2j_g}$ and $\lambda(h)=\gamma^{i_h+2j_h}$. If $g$ and $h$ are conjugate, then $i_g=i_h$. So $\gamma^{2j_g} = \gamma^{2j_h}$, i.e., $({\gamma^{j_g}})^2 = ({\gamma^{j_h}})^2$. So $\lambda(h)=\pm\lambda(g)$.

It is possible that an element $g$ of $\GSp(4,q)$ is conjugate to $-g$. For example,
\begin{displaymath}
\left(\begin{smallmatrix}
1 & & & \\
& 1 & & \\
& & -1 & \\
& & & -1 \\
\end{smallmatrix}\right)\cdot\left(\begin{smallmatrix}
& & 1 & \\
& & & 1 \\
1 & & & \\
& 1 & & \\
\end{smallmatrix}\right)\cdot\left(\begin{smallmatrix}
1 & & & \\
& 1 & & \\
& & -1 & \\
& & & -1 \\
\end{smallmatrix}\right)=\left(\begin{smallmatrix}
& & -1 & \\
& & & -1 \\
-1 & & & \\
& -1 & & \\
\end{smallmatrix}\right).
\end{displaymath}

The centralizers are sometimes affected when we pull back our representatives in $\PGSp(4,q)$ to $\GSp(4,q)$. There are two types of pullbacks. The first type consists of elements $g$ such that $g\neq x(-g)x^{-1}$ for any $x\in \GSp(4, q)$. The second type consists of elements $g$ such that $g=x(-g)x^{-1}$ for some $x\in \GSp(4, q)$.

{\em Type 1}. Let $g\in \GSp(4,q)$ be of the first type, i.e., $g$ is not conjugate to $-g$. Let $\overline{h}\in {\rm C}_{\PGSp(4,q)}(\overline{g})$. Then $\overline{g} = \overline{h}\cdot \overline{g}\cdot \overline{h^{-1}}$. When pulled back to $\GSp(4,q)$, $g = z_0hgh^{-1}$, with $z_0=\pm I$. We have $z_0\neq -I$ since $g$ is not conjugate to $-g$. So $z_0=I$, $g=hgh^{-1}$, and $h\in {\rm C}_{\GSp(4,q)}(g)$. We get a short exact sequence
\begin{displaymath}
1 \longrightarrow Z \longrightarrow {\rm C}_{\GSp(4,q)}(g) \longrightarrow {\rm C}_{\PGSp(4,q)}(\overline{g}) \longrightarrow 1.
\end{displaymath}
Therefore $\# {\rm C}_{\GSp(4,q)}(g) = (q-1)\cdot \# {\rm C}_{\PGSp(4,q)}(\overline{g})$.

{\em Type 2}. Let $g\in \GSp(4,q)$ be of the second type, i.e., $g$ is conjugate to $-g$. Define the set $S_g=\{h\in \GSp(4, F_q) : hgh^{-1}=-g\}$. Fix $s_0\in S_g$. The set $S_g$ is not a group, but there is a bijection of sets $S_g \longrightarrow {\rm C}_{\GSp(4,q)}(g)$, given by the map $h\mapsto s_0h$. Given $\overline{h}\in {\rm C}_{\PGSp(4,q)}(\overline{g})$, either $h\in {\rm C}_{\GSp(4,q)}(g)$ or $h\in S_g$. The set $S_g\sqcup {\rm C}_{\GSp(4,q)}(g)$ maps onto ${\rm C}_{\PGSp(4,q)}(\overline{g})$ via the projection map. Moreover, ${\rm C}'_{\GSp(4,q)}(g) := S_g\sqcup {\rm C}_{\GSp(4,q)}(g)$ is a group with respect to matrix multiplication and the projection map is a group homomorphism. ${\rm C}_{\GSp(4,q)}(g)$ is a subgroup of ${\rm C}'_{\GSp(4,q)}(g)$ of index 2. We get a short exact sequence
\begin{displaymath}
1 \longrightarrow Z \longrightarrow {\rm C}'_{\GSp(4,q)}(g) \longrightarrow {\rm C}_{\PGSp(4,q)}(\overline{g}) \longrightarrow 1.
\end{displaymath}
Then $$\# {\rm C}'_{\GSp(4,q)}(g) = (q-1) \cdot \# {\rm C}_{\PGSp(4,q)}(\overline{g}).$$ Also, $$2\cdot \# {\rm C}_{\GSp(4,q)}(g) = \# {\rm C}'_{\GSp(4,q)}(g).$$ So $\#{\rm C}_{\GSp(4,q)}(g) = \dfrac{q-1}{2} \cdot \# {\rm C}_{\PGSp(4,q)}(\overline{g})$.

Thus, given a class representative $\overline{g}\in \PGSp(4, q)$, we pull it back to $\GSp(4,q)$. Then, we determine if $g$ is of Type 1 or Type 2. If it is of Type 1, then there are $q-1$ distinct conjugacy classes $zg$, for $z\in Z$, each of order 
$$({\# \GSp(4, q)})/({(q-1)\cdot \# {\rm C}_{\PGSp(4,q)}(g)}).$$ If the pullback is of Type 2, then there are $({q-1})/{2}$ distinct conjugacy classes $z_i g$, with $i\in T_2.$

We list the conjugacy classes in the following table. The Type 1 conjugacy classes are $A_1(k), A_2(k), A_{31}(k), A_{32}(k), A_5(k), B_3(k), C_1(i,k),$ $C_3(i,k), $ $C_4(i,k), C_5(i,k)$, $C_6(i,j,k), D_1(i,k), D_2(i,k), D_4(i,j,k), D_5(i,k), D_6(i,k), D_7(i,j,k), D_8(i,k),$ and\\
$D_9(i,k).$ 

The Type 2 conjugacy classes are $B_{11}(k), B_{12}(k), B_{21}(k), B_{22}(k), B_{41}(k),$ $B_{42}(k)$, 
$B_{43}(k), B_{44}(k), B_{51}(k), B_{52}(k), C_{21}(i,k), C_{22}(i,k), D_{31}(i,k),$ and $D_{32}(i,k).$

Fix $a,b\in\F_q^\times$ such that $-a^2+b^2\gamma$ is a square and let $c=\gamma+1$. The conjugacy classes of our groups can now be given. The reader should note that there are $(q^2+2q+4)(q-1)$ conjugacy classes of $\GSp(4,q)$. 

The reader may notice that some of the class representatives given in the tables below are in a form with entries not in $\F_q$, but in $\F_{q^2}$ or $\F_{q^4}$. These are given as representatives of the conjugacy class in $\GSp(4,\overline{\F_q})$ that is stable under the action of the Frobenius, which raises each entry in the matrix to the $q$-th power.

\begin{center}
\begin{small}
\begin{longtable}{|c|c|c|}
\caption[Conjugacy classes of $\GSp(4,q)$]{\small{Conjugacy classes of $\GSp(4,q)$}} \label{GSpConjugacy} \\

\hline
   \multicolumn{1}{|c|}{\parbox{1.5cm}{\strut\centering\small\textbf{Notation}\strut}} &
   \multicolumn{1}{|c|}{\parbox{4.6cm}{\strut\centering\small\textbf{Class representative}\strut}} &
   \multicolumn{1}{|c|}{\parbox{2.4cm}{\strut\centering\small\textbf{Order of centralizer}\strut}} \\
\hline
\hline
\endfirsthead

\multicolumn{3}{c}{{\tablename} \thetable{} -- Continued} \\[0.5ex]
\hline
   \multicolumn{1}{|c|}{\parbox{1.5cm}{\strut\centering\small\textbf{Notation}\strut}} &
   \multicolumn{1}{|c|}{\parbox{4.6cm}{\strut\centering\small\textbf{Class representative}\strut}} &
   \multicolumn{1}{|c|}{\parbox{2.4cm}{\strut\centering\small\textbf{Order of centralizer}\strut}} \\
\hline
\hline
\endhead

 \hline
\endlastfoot

\begin{minipage}{15ex}
\begin{center}
\vspace{0.05in}
$A_1(k)$,\\
\scriptsize$k\in T_3$
\vspace{0.05in}
\end{center}
\end{minipage}
&
$\gamma^k\cdot I_4$
& $\# \GSp(4, q)$ \\
\hline
\begin{minipage}{15ex}
\begin{center}
\vspace{0.05in}
$A_2(k)$,\\
\scriptsize$k\in T_3$
\vspace{0.05in}
\end{center}
\end{minipage}
&
$\left(\begin{smallmatrix}
\gamma^k&&&\\
&\gamma^k&\gamma^k&\\
&&\gamma^k&\\
&&&\gamma^k
\end{smallmatrix}\right)$
& $q^4(q^2-1)(q-1)$ \\
\hline
\begin{minipage}{15ex}
\begin{center}
\vspace{0.05in}
$A_{31}(k)$,\\
\scriptsize$k\in T_3$
\vspace{0.05in}
\end{center}
\end{minipage}
&
$\left(\begin{smallmatrix}
\gamma^k & & & \gamma^k \\
& \gamma^k & -\gamma^k & \\
& & \gamma^k & \\
& & & \gamma^k \\
\end{smallmatrix}\right)$
& $2q^3(q-1)^2$ \\
\hline
\begin{minipage}{15ex}
\begin{center}
\vspace{0.05in}
$A_{32}(k)$,\\
\scriptsize$k\in T_3$
\vspace{0.05in}
\end{center}
\end{minipage}
&
$\left(\begin{smallmatrix}
\gamma^k & & & \gamma^{k+1} \\
& \gamma^k & -\gamma^k & \\
& & \gamma^k & \\
& & & \gamma^k \\
\end{smallmatrix}\right)$
& $2q^3(q^2-1)$ \\
\hline
\begin{minipage}{15ex}
\begin{center}
\vspace{0.05in}
$A_5(k)$,\\
\scriptsize$k\in T_3$
\vspace{0.05in}
\end{center}
\end{minipage}
&
$\left(\begin{smallmatrix}
\gamma^k & \gamma^k & -\gamma^k & \\
& \gamma^k & -\gamma^k & \\
& & \gamma^k & -\gamma^k \\
& & & \gamma^k \\
\end{smallmatrix}\right)$
& $q^2(q-1)$ \\
\hline
\begin{minipage}{15ex}
\begin{center}
\vspace{0.05in}
$B_{11}(k)$,\\
\scriptsize$k\in T_2$
\vspace{0.05in}
\end{center}
\end{minipage}
&
${\rm diag}(-\gamma^k,\gamma^k,\gamma^k,-\gamma^k)$
& $q^2(q^2-1)^2(q-1)$ \\
\hline
\begin{minipage}{15ex}
\begin{center}
\vspace{0.05in}
$B_{12}(k)$,\\
\scriptsize$k\in T_2$
\vspace{0.05in}
\end{center}
\end{minipage}
&
$\left(\begin{smallmatrix}
& \gamma^k & & \\
\gamma^{k+1} & & & \\
& & & \gamma^k \\
& & \gamma^{k+1} & \\
\end{smallmatrix}\right)$
& $q^2(q^4-1)(q-1)$ \\
\hline
\begin{minipage}{15ex}
\begin{center}
\vspace{0.05in}
$B_{21}(k)$,\\
\scriptsize$k\in T_2$
\vspace{0.05in}
\end{center}
\end{minipage}
&
${\rm diag}(\gamma^k,\gamma^k,-\gamma^k,-\gamma^k)$
& $q(q^2-1)(q-1)^2$ \\
\hline
\begin{minipage}{15ex}
\begin{center}
\vspace{0.05in}
$B_{22}(k)$,\\
\scriptsize$k\in T_2$
\vspace{0.05in}
\end{center}
\end{minipage}
&
$\left(\begin{smallmatrix}
& \gamma^k & & \\
\gamma^{k+1} & & & \\
& & & -\gamma^k \\
& & -\gamma^{k+1} & \\
\end{smallmatrix}\right)$
& $q(q^2-1)^2$ \\
\hline
\begin{minipage}{15ex}
\begin{center}
$B_{3}(k)$,\\
\scriptsize$k\in T_3$
\end{center}
\end{minipage}
&
$\left(\begin{smallmatrix}
-\gamma^k & & & \gamma^k \\
& \gamma^k & & \\
& & \gamma^k & \\
& & & -\gamma^k \\
\end{smallmatrix}\right)$
& $q^2(q^2-1)(q-1)$ \\
\hline
\begin{minipage}{15ex}
\begin{center}
$B_{41}(k)$,\\
\scriptsize$k\in T_2$
\end{center}
\end{minipage}
&
$\left(\begin{smallmatrix}
-\gamma^k & & & \gamma^k \\
& \gamma^k & -\gamma^k & \\
& & \gamma^k & \\
& & & -\gamma^k \\
\end{smallmatrix}\right)$
& $2q^2(q-1)$ \\
\hline
\begin{minipage}{15ex}
\begin{center}
$B_{42}(k)$,\\
\scriptsize$k\in T_2$
\end{center}
\end{minipage}
&
$\left(\begin{smallmatrix}
-\gamma^k & & & \gamma^k \\
& \gamma^k & -\gamma^{k+1} & \\
& & \gamma^k & \\
& & & -\gamma^k \\
\end{smallmatrix}\right)$
& $2q^2(q-1)$ \\
\hline
\begin{minipage}{15ex}
\begin{center}
$B_{43}(k)$,\\
\scriptsize$k\in T_2$
\end{center}
\end{minipage}
&
$\left(\begin{smallmatrix}
& \gamma^k & & -\frac{a}{2}\gamma^k \\
\gamma^{k+1} & & -\frac{a}{2}\gamma^{k+1} & b\gamma^{k+1} \\
& & & \gamma^k\\
& & \gamma^{k+1} & \\
\end{smallmatrix}\right)$
& $2q^2(q-1)$ \\
\hline

\begin{minipage}{15ex}
\begin{center}
$B_{44}(k)$,\\
\scriptsize$k\in T_2$
\end{center}
\end{minipage}
&
$\left(\begin{smallmatrix}
& \gamma^k & c\gamma^k & \\
\gamma^{k+1} & & & 2\gamma^{k+1} \\
& & & \gamma^k \\
& & \gamma^{k+1} & \\
\end{smallmatrix}\right)$
& $2q^2(q-1)$ \\
\hline
\begin{minipage}{15ex}
\begin{center}
$B_{51}(k)$,\\
\scriptsize$k\in T_2$
\end{center}
\end{minipage}
&
$\left(\begin{smallmatrix}
\gamma^k & -\gamma^k & & \\
& \gamma^k & & \\
& & -\gamma^k & -\gamma^k \\
& & & -\gamma^k \\
\end{smallmatrix}\right)$
& $q(q-1)^2$ \\
\hline
\begin{minipage}{15ex}
\begin{center}
$B_{52}(k)$,\\
\scriptsize$k\in T_2$
\end{center}
\end{minipage}
&
$\left(\begin{smallmatrix}
& \gamma^{k+1} & \gamma^k & \\
\gamma^k & & & \\
& & & -\gamma^{k+1} \\
& & -\gamma^k & \\
\end{smallmatrix}\right)$
& $q(q^2-1)$ \\
\hline
\begin{minipage}{15ex}
\vspace{0.05in}
\begin{center}
$C_1(i,k)$,\\
\scriptsize$i\in T_1, k\in T_3$
\vspace{0.02in}
\end{center}
\end{minipage}
&
${\rm diag}(\gamma^k,\gamma^k,\gamma^{k+i},\gamma^{k+i})$
& $q(q^2-1)(q-1)^2$ \\
\hline
\begin{minipage}{15ex}
\vspace{0.05in}
\begin{center}
$C_{21}(i,k)$,\\
\scriptsize$i\in T_1, k\in T_2$
\end{center}
\vspace{0.01in}
\end{minipage}
&
${\rm diag}(\gamma^k,-\gamma^k,\gamma^{k+i},-\gamma^{k+i})$
& $(q-1)^3$ \\
\hline
\begin{minipage}{15ex}
\vspace{0.05in}
\begin{center}
$C_{22}(i,k)$,\\
\scriptsize$i\in T_1, k\in T_2$
\vspace{0.01in}
\end{center}
\end{minipage}
&
$\left(\begin{smallmatrix}
& \gamma^k & & \\
\gamma^{k+1} & & & \\
& & & -\gamma^{k+i} \\
& & -\gamma^{k+i+1} & \\
\end{smallmatrix}\right)$
& $(q^2-1)(q-1)$ \\
\hline
\begin{minipage}{15ex}
\vspace{0.05in}
\begin{center}
$C_3(i,k)$,\\
\scriptsize$i\in T_1, k\in T_3$
\end{center}
\vspace{0.01in}
\end{minipage}
&
$\left(\begin{smallmatrix}
\gamma^k & -\gamma^k & & \\
& \gamma^k & & \\
& & \gamma^{k+i} & \gamma^{k+i} \\
& & & \gamma^{k+i} \\
\end{smallmatrix}\right)$
& $q(q-1)^2$ \\
\hline

\begin{minipage}{15ex}
\vspace{0.05in}
\begin{center}
$C_4(i,k)$,\\
\scriptsize$i\in T_1, k\in T_3$
\end{center}
\vspace{0.01in}
\end{minipage}
&
$\left(\begin{smallmatrix}
\gamma^k & & & -\gamma^k \\
& \gamma^{k+i} & & \\
& & \gamma^{k-i} & \\
& & & \gamma^k \\
\end{smallmatrix}\right)$
& $q(q-1)^2$ \\
\hline
\begin{minipage}{15ex}
\vspace{0.05in}
\begin{center}
$C_5(i,k)$,\\
\scriptsize$i\in T_1, k\in T_3$
\vspace{0.02in}
\end{center}
\end{minipage}
&
${\rm diag}(\gamma^{k+i},\gamma^k,\gamma^k,\gamma^{k-i})$
& $q(q^2-1)(q-1)^2$ \\
\hline
\begin{minipage}{16ex}
\vspace{0.05in}
\begin{center}
$C_6(i,j,k)$,\\
\scriptsize$i,j\in T_1, i<j,$\\
\scriptsize$k\in T_3$
\end{center}
\vspace{0.01in}
\end{minipage}
&
${\rm diag}(\gamma^k,\gamma^{k+i},\gamma^{k+j},\gamma^{k+i+j})$
& $(q-1)^3$ \\
\hline
\begin{minipage}{15ex}
\vspace{0.05in}
\begin{center}
$D_1(i,k)$,\\
\scriptsize$i\in R_2, k\in T_3$
\end{center}
\vspace{0.01in}
\end{minipage}
&
${\rm diag}(\gamma^k,\gamma^k\theta^i,\gamma^k\theta^{qi},\gamma^{k+i})$
& $(q^2-1)(q-1)$ \\
\hline
\begin{minipage}{15ex}
\vspace{0.05in}
\begin{center}
$D_2(i,k)$,\\
\scriptsize$i\in T_2, k\in T_3$
\end{center}
\vspace{0.01in}
\end{minipage}
&
${\rm diag}(\gamma^k\theta^i,\gamma^k\theta^{qi},\gamma^k\theta^i,\gamma^k\theta^{qi})$
& $q(q^2-1)^2$ \\
\hline
\begin{minipage}{15ex}
\vspace{0.05in}
\begin{center}
$D_{31}(i,k)$,\\
\scriptsize$i,k\in T_2$
\end{center}
\vspace{0.01in}
\end{minipage}
&
$\left(\begin{smallmatrix}
-\gamma^k\theta^i & & & \\
& \gamma^k\theta^{qi} & & \\
& & -\gamma^k\theta^{qi} & \\
& & & \gamma^k\theta^i \\
\end{smallmatrix}\right)$
& $(q^2-1)(q-1)$ \\
\hline
\begin{minipage}{15ex}
\vspace{0.05in}
\begin{center}
$D_{32}(i,k)$,\\
\scriptsize$i,k\in T_2$
\end{center}
\vspace{0.01in}
\end{minipage}
&
$\gamma^k\cdot{\rm diag}(\gamma^{1/2}\eta^i, -\gamma^{1/2}\eta^i, -\gamma^{1/2}\eta^{-i}, \gamma^{1/2}\eta^{-i})$
& $(q^2-1)(q+1)$ \\
\hline
\begin{minipage}{15ex}
\vspace{0.05in}
\begin{center}
$D_4(i,j,k)$,\\
\scriptsize$i\in T_2, j\in T_1,$\\
\scriptsize$k\in T_3$
\end{center}
\vspace{0.01in}
\end{minipage}
&
${\rm diag}(\gamma^k\theta^{qi},\gamma^k\theta^i,\gamma^{k+j}\theta^{qi},\gamma^{k+j}\theta^i)$
& $(q^2-1)(q-1)$ \\
\hline
\begin{minipage}{15ex}
\vspace{0.05in}
\begin{center}
$D_5(i,k)$,\\
\scriptsize$i\in T_2, k\in T_3$
\end{center}
\vspace{0.01in}
\end{minipage}
&
$\left(\begin{smallmatrix}
\gamma^k\theta^i & -\gamma^k\theta^i & & \\
& \gamma^k\theta^i & & \\
& & \gamma^k\theta^{qi} & \gamma^k\theta^{qi} \\
& & & \gamma^k\theta^{qi} \\
\end{smallmatrix}\right)$
& $q(q^2-1)$ \\
\hline
\begin{minipage}{15ex}
\vspace{0.05in}
\begin{center}
$D_6(i,k)$,\\
\scriptsize$i\in T_2, k\in T_3$
\end{center}
\vspace{0.01in}
\end{minipage}
&
${\rm diag}(\gamma^k,\gamma^k\eta^i,\gamma^k\eta^{-i},\gamma^k)$
& $q(q^2-1)^2$ \\
\hline
\begin{minipage}{16ex}
\vspace{0.05in}
\begin{center}
$D_7(i,j,k)$,\\
\scriptsize$i,j\in T_2, i<j,$\\
\scriptsize$k\in T_3$
\end{center}
\vspace{0.01in}
\end{minipage}
&
${\rm diag}(\gamma^k\eta^i, \gamma^k\eta^j, \gamma^k\eta^{-j}, \gamma^k\eta^{-i})$
& $(q^2-1)(q+1)$ \\
\hline
\begin{minipage}{15ex}
\vspace{0.05in}
\begin{center}
$D_8(i,k)$,\\
\scriptsize$i\in T_2, k\in T_3$
\end{center}
\vspace{0.01in}
\end{minipage}
&
$\left(\begin{smallmatrix}
\gamma^k & & & -\gamma^k \\
& \gamma^k\eta^i & & \\
& & \gamma^k\eta^{-i} & \\
& & & \gamma^k \\
\end{smallmatrix}\right)$
& $q(q^2-1)$ \\
\hline
\begin{minipage}{15ex}
\vspace{0.05in}
\begin{center}
$D_9(i,k)$,\\
\scriptsize$i\in R_1, k\in T_3$
\end{center}
\vspace{0.01in}
\end{minipage}
&
${\rm diag}(\gamma^k\zeta^{qi}, \gamma^k\zeta^i, \gamma^k\zeta^{q^2i}, \gamma^k\zeta^{q^3i})$
& $(q^2+1)(q-1)$ \\
\end{longtable}
\end{small}
\end{center}

As noted before, the conjugacy classes of $\GSp(4,q)$ were previously determined by Shinoda \cite{ShinodaC} using a different approach. We give a correspondence between our different notations. The classes $D_{32}(i,k)$ and $D_7(i,j,k)$ turn out to be different members of the class type $L_0$ from \cite{Shinoda}.
 
\begin{center}
\renewcommand{\arraystretch}{1.3}
\begin{small}
\begin{longtable}{|l|l||l|l|}
\caption[Notations for conjugacy classes of $\GSp(4,q)$]{\small{Notations for conjugacy classes of $\GSp(4,q)$}} \label{GSpConjugacyNotations} \\

\hline
   \multicolumn{1}{|c|}{\parbox{2cm}{\strut\centering\small\textbf{Class}\strut}} &
   \multicolumn{1}{|c|}{\parbox{2cm}{\strut\centering\small\textbf{Shinoda}\strut}} &
   \multicolumn{1}{|c|}{\parbox{2cm}{\strut\centering\small\textbf{Class}\strut}} &
   \multicolumn{1}{|c|}{\parbox{2cm}{\strut\centering\small\textbf{Shinoda}\strut}} \\
\hline
\hline
\endfirsthead

\multicolumn{3}{c}{{\tablename} \thetable{} -- Continued} \\[0.5ex]
\hline
   \multicolumn{1}{|c|}{\parbox{2cm}{\strut\centering\small\textbf{Class}\strut}} &
   \multicolumn{1}{|c|}{\parbox{2cm}{\strut\centering\small\textbf{Shinoda}\strut}} &
   \multicolumn{1}{|c|}{\parbox{2cm}{\strut\centering\small\textbf{Class}\strut}} &
   \multicolumn{1}{|c|}{\parbox{2cm}{\strut\centering\small\textbf{Shinoda}\strut}} \\
\hline
\hline
\endhead

 \hline
\endlastfoot

$A_1(k), k\in T_3$ & $A_0, a=\gamma^k$ & $C_1(i,k),$ & $D_0, a=\gamma^k, b=\gamma^{k+i}$\\
& & $i\in T_1, k\in T_3$ & \\
\hline
$A_2(k), k\in T_3$ & $A_1, a=\gamma^k$ & $C_{21}(i,k)$ & $H_0, a_1=\gamma^k, $\\
& & $i\in T_1, k\in T_2$ & $a_2=\gamma^{k+i}, b_2=-\gamma^k$ \\
\hline
$A_{31}(k), k\in T_3$ & $A_{21}, a=\gamma^k$ & $C_{22}(i,k)$ & $I_0, a=\gamma^{k+1/2},$\\
& & $i\in T_1, k\in T_2$ & $b=\gamma^{2k+i}$\\
\hline
$A_{32}(k), k\in T_3$ & $A_{22}, a=\gamma^k$ & $C_3(i,k),$ & $D_1, a=\gamma^k, b=\gamma^{k+i}$\\
&& $i\in T_1, k\in T_3$ &\\
\hline
$A_5(k), k\in T_3$ & $A_3, a=\gamma^k$ & $C_4(i,k),$ & $E_1, a=\gamma^i, c=\gamma^k$ \\
& & $i\in T_1, k\in T_3$ &\\
\hline
$B_{11}(k), k\in T_2$ & $B_0, a=\gamma^k$ & $C_5(i,k),$ & $E_0, a=\gamma^i, c=\gamma^k$ \\
& & $i\in T_1, k\in T_3$ & \\
\hline
$B_{12}(k), k\in T_2$ & $C_0, a=\gamma^{k+1/2}$ & $C_6(i,j,k),$ & $H_0, a_1=\gamma^k, $ \\
& & $i,j\in T_1, i<j,$ & $a_2=\gamma^{k+i},$\\
&& $k\in T_3$ &$b_2=\gamma^{k+j}$\\
\hline
$B_{21}(k), k\in T_2$ & $D_0, a=\gamma^k, b=-\gamma^k$ & $D_1(i,k),$& $J_0, a=\gamma^k\theta^i, c=\gamma^k$ \\
&&$i\in R_2, k\in T_3$ & \\
\hline
$B_{22}(k), k\in T_2$ & $F_0, a=\gamma^{k+1/2}$ & $D_2(i,k),$ & $F_0, a=\gamma^k\theta^i$\\
&& $i\in T_2, k\in T_3$ &\\
\hline
$B_3(k), k\in T_2$ & $B_2, a=\gamma^k$ & $D_{31}(i,k), i,k\in T_2$ & $I_0, a=\gamma^k\theta^i, b=-\gamma^{2k}$\\
\hline
$B_3\left(k+\frac{q-1}{2}\right), $ & $B_1, a=\gamma^k$ & $D_{32}(i,k), i,k\in T_2$ & $L_0, a_1=\gamma^{k+1/2}\eta^i,$\\
$k\in T_2$ & & & $a_2=-\gamma^{k+1/2}\eta^i$\\
\hline
\newpage
$B_{41}(k), k\in T_2$ & $B_{31}, a=\gamma^k$ & $D_4(i,j,k),$ & $I_0, a=\gamma^k\theta^i, b=\gamma^{2k+j}$\\
&& $i\in T_2, j\in T_1,$ & \\
&& $k\in T_3$ & \\
\hline
$B_{42}(k), k\in T_2$ & $B_{32}, a=\gamma^k$ & $D_5(i,k),$ & $F_1, a=\gamma^k\theta^i$\\
&& $i\in T_2, k\in T_3$ & \\
\hline
$B_{43}(k), k\in T_2$ & $C_{12}, a=\gamma^{k+1/2}$ & $D_6(i,k),$ & $G_0, a=\gamma^k, u=\eta^i$ \\
&& $i\in T_2, k\in T_3$ & \\
\hline
$B_{44}(k), k\in T_2$ & $C_{11}, a=\gamma^{k+1/2}$ & $D_7(i,j,k)$ & $L_0, a_1=\gamma^k\eta^i,$\\ 
&& $i,j\in T_2, i<j,$ & $a_2=\gamma^k\eta^j$\\
&&$k\in T_3$&\\
\hline
$B_{51}(k), k\in T_2$ & $D_1, a=\gamma^k, b=-\gamma^k$ & $D_8(i,k),$ & $G_1, a=\gamma^k, u=\eta^i$\\
&& $i\in T_2, k\in T_3$ & \\
\hline
$B_{52}(k), k\in T_2$ & $F_1, a=\gamma^{k+1/2}$ & $D_9(i,k),$ & $K_0, a=\gamma^k\zeta^i$\\
&& $i\in R_1, k\in T_3$ & \\
\end{longtable}
\end{small}
\end{center}

Now we compute the conjugacy classes of the Borel subgroup, the Siegel parabolic subgroup, the Klingen parabolic subgroup, and the unipotent radical of the Borel. In each of the following tables of conjugacy classes of our important subgroups, the notation will indicate which conjugacy classes of $\GSp(4,q)$ occur in the subgroup and how many components the splitting has if the class splits into multiple classes in the subgroup. For example, the Borel subgroup is denoted by $B$. The class $A_2(k)$ has a non-empty intersection with $B$, splitting into two conjugacy classes, denoted by $BA_2^1(k)$ and $BA_2^2(k)$. 

\begin{longtable}{|c|c|c|}
\caption[Conjugacy classes of the Borel subgroup]{\small{Conjugacy classes of the Borel subgroup}} \label{BConjugacy} \\

\hline
   \multicolumn{1}{|c|}{\parbox{1.5cm}{\strut\centering\small\textbf{Notation}\strut}} &
   \multicolumn{1}{|c|}{\parbox{4.6cm}{\strut\centering\small\textbf{Class representative}\strut}} &
   \multicolumn{1}{|c|}{\parbox{2.4cm}{\strut\centering\small\textbf{Order of centralizer in $B$}\strut}} \\
\hline
\hline
\endfirsthead

\multicolumn{3}{c}{{\tablename} \thetable{} -- Continued} \\[0.5ex]
\hline
   \multicolumn{1}{|c|}{\parbox{1.5cm}{\strut\centering\small\textbf{Notation}\strut}} &
   \multicolumn{1}{|c|}{\parbox{4.6cm}{\strut\centering\small\textbf{Class representative}\strut}} &
   \multicolumn{1}{|c|}{\parbox{2.4cm}{\strut\centering\small\textbf{Order of centralizer in $B$}\strut}} \\
\hline
\hline
\endhead

 \hline
\endlastfoot

\begin{minipage}{15ex}
\vspace{0.05in}
\begin{center}
$BA_1(k)$,\\
\scriptsize$k\in T_3$
\end{center}
\vspace{0.01in}
\end{minipage}
&
$\gamma^k\cdot I_4$
& $q^4(q-1)^3$ \\
\hline
\begin{minipage}{15ex}
\vspace{0.05in}
\begin{center}
$BA^1_2(k)$,\\
\scriptsize$k\in T_3$
\end{center}
\vspace{0.01in}
\end{minipage}
&
$\left(\begin{smallmatrix}
\gamma^k & & & \\
& \gamma^k & \gamma^k & \\
& & \gamma^k & \\
& & & \gamma^k \\
\end{smallmatrix}\right)$
& $q^3(q-1)^2$ \\
\hline
\begin{minipage}{15ex}
\vspace{0.05in}
\begin{center}
$BA^2_2(k)$,\\
\scriptsize$k\in T_3$
\end{center}
\vspace{0.01in}
\end{minipage}
&
$\left(\begin{smallmatrix}
\gamma^k & & & \gamma^k \\
& \gamma^k & & \\
& & \gamma^k & \\
& & & \gamma^k \\
\end{smallmatrix}\right)$
& $q^4(q-1)^2$ \\
\hline
\begin{minipage}{15ex}
\vspace{0.05in}
\begin{center}
$BA^1_{31}(k)$,\\
\scriptsize$k\in T_3$
\end{center}
\vspace{0.01in}
\end{minipage}
&
$\left(\begin{smallmatrix}
\gamma^k & & & \gamma^k \\
& \gamma^k & -\gamma^k & \\
& & \gamma^k & \\
& & & \gamma^k \\
\end{smallmatrix}\right)$
& $2q^3(q-1)$ \\
\hline

\begin{minipage}{15ex}
\vspace{0.05in}
\begin{center}
$BA^2_{31}(k)$,\\
\scriptsize$k\in T_3$
\end{center}
\vspace{0.01in}
\end{minipage}
&
$\left(\begin{smallmatrix}
\gamma^k & & \gamma^k &  \\
& \gamma^k & & \gamma^k \\
& & \gamma^k & \\
& & & \gamma^k \\
\end{smallmatrix}\right)$
& $q^3(q-1)^2$ \\
\end{longtable}

\begin{longtable}{|c|c|c|}
\caption[Conjugacy classes of the Siegel parabolic subgroup]{\footnotesize{Conjugacy classes of the Siegel parabolic subgroup}} \label{PConjugacy} \\

\hline
   \multicolumn{1}{|c|}{\parbox{1.5cm}{\strut\centering\small\textbf{Notation}\strut}} &
   \multicolumn{1}{|c|}{\parbox{4.6cm}{\strut\centering\small\textbf{Class representative}\strut}} &
   \multicolumn{1}{|c|}{\parbox{2.4cm}{\strut\centering\small\textbf{Order of centralizer in $P$}\strut}} \\
\hline
\hline
\endfirsthead

\multicolumn{3}{c}{{\tablename} \thetable{} -- Continued} \\[0.5ex]
\hline
   \multicolumn{1}{|c|}{\parbox{1.5cm}{\strut\centering\small\textbf{Notation}\strut}} &
   \multicolumn{1}{|c|}{\parbox{4.6cm}{\strut\centering\small\textbf{Class representative}\strut}} &
   \multicolumn{1}{|c|}{\parbox{2.4cm}{\strut\centering\small\textbf{Order of centralizer in $P$}\strut}} \\
\hline
\hline
\endhead

 \hline
\endlastfoot

\begin{minipage}{15ex}
\vspace{0.05in}
\begin{center}
$PA_1(k)$,\\
\scriptsize$k\in T_3$
\end{center}
\vspace{0.01in}
\end{minipage}
&
$\gamma^k\cdot I_4$
& $q^4(q^2-1)(q-1)^2$ \\
\hline
\begin{minipage}{15ex}
\vspace{0.05in}
\begin{center}
$PA_2(k)$,\\
\scriptsize$k\in T_3$
\end{center}
\vspace{0.01in}
\end{minipage}
&
$\left(\begin{smallmatrix}
\gamma^k & & & \\
& \gamma^k & \gamma^k & \\
& & \gamma^k & \\
& & & \gamma^k \\
\end{smallmatrix}\right)$
& $q^4(q-1)^2$ \\
\hline
\begin{minipage}{15ex}
\vspace{0.05in}
\begin{center}
$PA^1_{31}(k)$,\\
\scriptsize$k\in T_3$
\end{center}
\vspace{0.01in}
\end{minipage}
&
$\left(\begin{smallmatrix}
\gamma^k & & & \gamma^k \\
& \gamma^k & -\gamma^k & \\
& & \gamma^k & \\
& & & \gamma^k \\
\end{smallmatrix}\right)$
& $2q^3(q-1)^2$ \\
\hline

\begin{minipage}{15ex}
\vspace{0.05in}
\begin{center}
$PA^2_{31}(k)$,\\
\scriptsize$k\in T_3$
\end{center}
\vspace{0.01in}
\end{minipage}
&
$\left(\begin{smallmatrix}
\gamma^k & \gamma^k & & \\
& \gamma^k & & \\
& & \gamma^k & -\gamma^k \\
& & & \gamma^k \\
\end{smallmatrix}\right)$
& $q^2(q-1)^2$ \\
\hline
\begin{minipage}{15ex}
\vspace{0.05in}
\begin{center}
$PA_{32}(k)$,\\
\scriptsize$k\in T_3$
\end{center}
\vspace{0.01in}
\end{minipage}
&
$\left(\begin{smallmatrix}
\gamma^k & & & \gamma^{k+1} \\
& \gamma^k & -\gamma^k & \\
& & \gamma^k & \\
& & & \gamma^k \\
\end{smallmatrix}\right)$
& $2q^3(q^2-1)$ \\
\end{longtable}

\begin{longtable}{|c|c|c|c|}
\caption[Conjugacy classes of the Klingen parabolic subgroup]{\footnotesize{Conjugacy classes of the Klingen parabolic subgroup}} \label{QConjugacy} \\

\hline
   \multicolumn{1}{|c|}{\parbox{1.5cm}{\strut\centering\small\textbf{Notation}\strut}} &
   \multicolumn{1}{|c|}{\parbox{4.6cm}{\strut\centering\small\textbf{Class representative}\strut}} &
   \multicolumn{1}{|c|}{\parbox{2.4cm}{\strut\centering\small\textbf{Order of centralizer in $Q$}\strut}} \\
\hline
\hline
\endfirsthead

\multicolumn{3}{c}{{\tablename} \thetable{} -- Continued} \\[0.5ex]
\hline
   \multicolumn{1}{|c|}{\parbox{1.5cm}{\strut\centering\small\textbf{Notation}\strut}} &
   \multicolumn{1}{|c|}{\parbox{4.6cm}{\strut\centering\small\textbf{Class representative}\strut}} &
   \multicolumn{1}{|c|}{\parbox{2.4cm}{\strut\centering\small\textbf{Order of centralizer in $Q$}\strut}} \\
\hline
\hline
\endhead

 \hline
\endlastfoot

\begin{minipage}{15ex}
\vspace{0.05in}
\begin{center}
$QA_1(k)$,\\
\scriptsize$k\in T_3$
\end{center}
\vspace{0.01in}
\end{minipage}
&
$\gamma^k\cdot I_4$
& $q^4(q^2-1)(q-1)^2$ \\
\hline
\begin{minipage}{15ex}
\vspace{0.05in}
\begin{center}
$QA^1_2(k)$,\\
\scriptsize$k\in T_3$
\end{center}
\vspace{0.01in}
\end{minipage}
&
$\left(\begin{smallmatrix}
\gamma^k & & & \\
& \gamma^k & \gamma^k & \\
& & \gamma^k & \\
& & & \gamma^k \\
\end{smallmatrix}\right)$
& $q^3(q-1)^2$ \\
\hline
\begin{minipage}{15ex}
\vspace{0.05in}
\begin{center}
$QA^2_2(k)$,\\
\scriptsize$k\in T_3$
\end{center}
\vspace{0.01in}
\end{minipage}
&
$\left(\begin{smallmatrix}
\gamma^k & & & \gamma^k \\
& \gamma^k & & \\
& & \gamma^k & \\
& & & \gamma^k \\
\end{smallmatrix}\right)$
& $q^4(q^2-1)(q-1)$ \\
\hline
\begin{minipage}{15ex}
\vspace{0.05in}
\begin{center}
$QA^1_{31}(k)$,\\
\scriptsize$k\in T_3$
\end{center}
\vspace{0.01in}
\end{minipage}
&
$\left(\begin{smallmatrix}
\gamma^k & & & \gamma^k \\
& \gamma^k & -\gamma^k & \\
& & \gamma^k & \\
& & & \gamma^k \\
\end{smallmatrix}\right)$
& $2q^3(q-1)$ \\
\hline

\begin{minipage}{15ex}
\vspace{0.05in}
\begin{center}
$QA^2_{31}(k)$,\\
\scriptsize$k\in T_3$
\end{center}
\vspace{0.01in}
\end{minipage}
&
$\left(\begin{smallmatrix}
\gamma^k & & \gamma^k & \\
& \gamma^k & & \gamma^k \\
& & \gamma^k & \\
& & & \gamma^k \\
\end{smallmatrix}\right)$
& $q^3(q-1)^2$ \\
\end{longtable}

\begin{longtable}{|c|c|c|}
\caption[Conjugacy classes of $N_{\rm GSp(4)}$]{\footnotesize{Conjugacy classes of $N_{\rm GSp(4)}$}} \label{NConjugacy} \\

\hline
   \multicolumn{1}{|c|}{\parbox{1.5cm}{\strut\centering\small\textbf{Notation}\strut}} &
   \multicolumn{1}{|c|}{\parbox{4.6cm}{\strut\centering\small\textbf{Class representative}\strut}} &
   \multicolumn{1}{|c|}{\parbox{2.4cm}{\strut\centering\small\textbf{Order of centralizer in $N_{\rm GSp(4)}$}\strut}} \\
\hline
\hline
\endfirsthead

\multicolumn{3}{c}{{\tablename} \thetable{} -- Continued} \\[0.5ex]
\hline
   \multicolumn{1}{|c|}{\parbox{1.5cm}{\strut\centering\small\textbf{Notation}\strut}} &
   \multicolumn{1}{|c|}{\parbox{4.6cm}{\strut\centering\small\textbf{Class representative}\strut}} &
   \multicolumn{1}{|c|}{\parbox{2.4cm}{\strut\centering\small\textbf{Order of centralizer in $N_{\rm GSp(4)}$}\strut}} \\
\hline
\hline
\endhead

 \hline
\endlastfoot

$NA_1$
&
$I_4$
& $q^4$ \\
\hline
\begin{minipage}{15ex}
\vspace{0.05in}
\begin{center}
$NA_2^1(k)$,\\
\scriptsize$k\in T_3$
\end{center}
\vspace{0.01in}
\end{minipage}
&
$\left(\begin{smallmatrix}
1 & & & \\
& 1 & \gamma^k & \\
& & 1 & \\
& & & 1 \\
\end{smallmatrix}\right)$
& $q^3$ \\
\hline
\begin{minipage}{15ex}
\vspace{0.05in}
\begin{center}
$NA_2^2(k)$,\\
\scriptsize$k\in T_3$
\end{center}
\vspace{0.01in}
\end{minipage}
&
$\left(\begin{smallmatrix}
1 & & & \gamma^k \\
& 1 & & \\
& & 1 & \\
& & & 1 \\
\end{smallmatrix}\right)$
& $q^4$ \\
\hline
\begin{minipage}{13ex}
\begin{center}
$NA_{31}^1(i ,j , y)$,\\
\scriptsize$i,j\in T_3,$\\
{$\gamma^{2i}-\gamma^jy = \gamma^{2n}$}\\
for some $n\in T_3$
\end{center}
\vspace{0.01in}
\end{minipage}
&
$\left(\begin{smallmatrix}
1 & & \gamma^i & y \\
& 1 & \gamma^j & \gamma^i \\
& & 1 & \\
& & & 1 \\
\end{smallmatrix}\right)$
& $q^3$ \\
\hline
\begin{minipage}{15ex}
\vspace{0.05in}
\begin{center}
$NA_{31}^2(k)$,\\
\scriptsize$k\in T_3$
\end{center}
\vspace{0.01in}
\end{minipage}
&
$\left(\begin{smallmatrix}
1 & \gamma^k & & \\
& 1 & & \\
& & 1 & -\gamma^k \\
& & & 1 \\
\end{smallmatrix}\right)$
& $q^2$ \\
\hline
\begin{minipage}{15ex}
\vspace{0.05in}
\begin{center}
$NA_{31}^3(K)$,\\
\scriptsize$k\in T_3$
\end{center}
\vspace{0.01in}
\end{minipage}
&
$\left(\begin{smallmatrix}
1 & & \gamma^k & \\
& 1 & & \gamma^k \\
& & 1 & \\
& & & 1 \\
\end{smallmatrix}\right)$
& $q^3$ \\
\end{longtable}

\section{Induced characters}\label{inducedchars}

This section contains the character tables of representations induced from representations of $N_{\GSp(4)}$, $B$, $P$, and $Q$. If a conjugacy class of $\GSp(4,q)$ is not listed in the table, the character takes the value 0 on that conjugacy class.

\subsection{Generic characters}\label{Genericreps}
Let $\psi_1$ and $\psi_2$ be non-trivial characters of $\F_q$ and let $\psi_N$ be the character of $N_{\GSp(4)}$ defined by
\begin{displaymath}
\psi_N \left(%
\right) = \psi_1(x)\psi_2(y).
\end{displaymath}
This defines a representation of $N_{\GSp(4)}$. The representation ${\rm Ind}_{N_{\GSp(4)}}^G(\psi_N)$ will be denoted by $\mathcal{G}$ and its character by $\chi_\mathcal{G}$. If $(\pi, V)$ is an irreducible representation that can be embedded into $\mathcal{G}$, we call its image a \emph{Whittaker model of} $\pi$ and say that $\pi$ is \emph{generic}. The character table of $\mathcal{G}$ is Table \ref{GCharacterValues}.

\begin{center}
\begin{small}
\begin{longtable}{|c|l|}
\caption[$\mathcal{G}$ character values]{\small{$\mathcal{G}$ character values}} \label{GCharacterValues} \\

\hline
   \multicolumn{1}{|c|}{\parbox{1.5cm}{\strut\centering\small\textbf{Class}\strut}} &
   \multicolumn{1}{|c|}{\parbox{3.5cm}{\strut\centering\small\textbf{$\mathcal{G}$ character value}\strut}} \\
\hline
\hline
\endfirsthead

\multicolumn{2}{c}{{\tablename} \thetable{} -- Continued} \\[0.5ex]
\hline
   \multicolumn{1}{|c|}{\parbox{1.5cm}{\strut\centering\small\textbf{Class}\strut}} &
   \multicolumn{1}{|c|}{\parbox{3.5cm}{\strut\centering\small\textbf{$\mathcal{G}$ character value}\strut}} \\
\hline
\hline
\endhead

 \hline
\endlastfoot

$A_1(q-1)$ & $(q^4-1)(q^2-1)(q-1)$ \\
\hline
$A_2(q-1)$ & $-(q^2-1)(q-1)$ \\
\hline
$A_{31}(q-1)$ & $-(q^2-1)(q-1)$ \\
\hline
$A_{32}(q-1)$ & $-(q^2-1)(q-1)$ \\
\hline
$A_5(q-1)$ & $q-1$ \\
\end{longtable}
\end{small}
\end{center}
We compute $(\chi_\mathcal{G}, \chi_\mathcal{G}) = q^2(q-1)$. So there are precisely $q^2(q-1)$ irreducible generic representations of $\GSp(4,q)$ since Gelfand-Graev representations of finite groups of Lie type are multiplicity free,  \cite{DM}, \cite{GG}, \cite{SteinGG}.

\subsection{Borel}\label{Borelreps}
Let $\chi_1, \chi_2,$ and $\sigma$ be characters of the multiplicative group $\F_q^\times$. Define a representation $\pi$ on the Borel subgroup $B$ by
\begin{displaymath}
\left(%
\right) \longmapsto \chi_1(a)\chi_2(b)\sigma(c).
\end{displaymath}
The character of this representation is given by $\chi_1(a)\chi_2(b)\sigma(c)$. The representation ${\rm Ind}_B^G(\pi)$ will be denoted by $\chi_1 \times \chi_2 \rtimes \sigma$. The central character of $\chi_1 \times \chi_2 \rtimes \sigma$ is $\chi_1 \chi_2 \sigma^2$. The character table of $\chi_1 \times \chi_2 \rtimes \sigma$ is Table \ref{BCharacterValues}.

\begin{center}
\begin{small}
\begin{longtable}{|c|l|}
\caption[$\chi_1 \times \chi_2 \rtimes \sigma$ character values]{\small{$\chi_1 \times \chi_2 \rtimes \sigma$ character values}} \label{BCharacterValues} \\

\hline
   \multicolumn{1}{|c|}{\parbox{1cm}{\strut\centering\small\textbf{Class}\strut}} &
   \multicolumn{1}{|c|}{\parbox{6cm}{\strut\centering\small\textbf{$\chi_1 \times \chi_2 \rtimes \sigma$ character value}\strut}} \\
\hline
\hline
\endfirsthead

\multicolumn{2}{c}{{\tablename} \thetable{} -- Continued} \\[0.5ex]
\hline
   \multicolumn{1}{|c|}{\parbox{1cm}{\strut\centering\small\textbf{Class}\strut}} &
   \multicolumn{1}{|c|}{\parbox{6cm}{\strut\centering\small\textbf{$\chi_1 \times \chi_2 \rtimes \sigma$ character value}\strut}} \\
\hline
\hline
\endhead

 \hline
\endlastfoot

$A_1(k)$ & $(q^2+1)(q+1)^2 \chi_1(\gamma^k)\chi_2(\gamma^k)\sigma(\gamma^{2k})$ \\
\hline
$A_2(k)$ & $(q+1)^2\chi_1(\gamma^k)\chi_2(\gamma^k)\sigma(\gamma^{2k})$ \\
\hline
$A_{31}(k)$ & $(3q+1)\chi_1(\gamma^k)\chi_2(\gamma^k)\sigma(\gamma^{2k})$ \\
\hline
$A_{32}(k)$ & $(q+1)\chi_1(\gamma^k)\chi_2(\gamma^k)\sigma(\gamma^{2k})$ \\
\hline
$A_5(k)$ & $\chi_1(\gamma^k)\chi_2(\gamma^k)\sigma(\gamma^{2k})$ \\
\hline
$B_{11}(k)$ & $(q+1)^2\chi_1(\gamma^k)\chi_2(\gamma^k)\sigma(\gamma^{2k})\Bigl(\chi_1(-1)+\chi_2(-1)\Bigr)$ \\
\hline
$B_{21}(k)$ & \small$(q+1)\chi_1(\gamma^k)\chi_2(\gamma^k)\sigma(-\gamma^{2k})\Bigl(1+\chi_1(-1)\chi_2(-1)+\chi_1(-1)+\chi_2(-1)\Bigr)$\\
\hline
$B_3(k)$ & $(q+1)\chi_1(\gamma^k)\chi_2(\gamma^k)\sigma(\gamma^{2k})\Bigl(\chi_1(-1)+\chi_2(-1)\Bigr)$ \\
\end{longtable}
\end{small}
\end{center}

\subsection{Siegel}\label{Siegelreps}
Let $(\pi, V)$ be an irreducible representation of $\GL(2, q)$, and let $\sigma$ be a character of $\F_q^\times$. Define a representation $\pi'$ of the Siegel parabolic subgroup $P$ on $V$ by
\begin{displaymath}
\left(%
\right) \longmapsto \sigma(\lambda)\pi(A).
\end{displaymath}
The character of this representation is given by $\sigma(\lambda)\chi_\pi(A)$, where $\chi_\pi$ is the character of $\pi$. The representation ${\rm Ind}_P^G(\pi')$ will be denoted by $\pi \rtimes \sigma$. If $\pi$ has central character $\omega_\pi$, then the central character of $\pi \rtimes \sigma$ is $\omega_\pi \sigma^2$. The character table of $\pi \rtimes \sigma$ is Table \ref{PCharacterValues}.
\begin{center}
\begin{small}
\begin{longtable}{|c|l|}
\caption[$\pi \rtimes \sigma$ character values]{\small{$\pi \rtimes \sigma$ character values}} \label{PCharacterValues} \\

\hline
   \multicolumn{1}{|c|}{\parbox{1cm}{\strut\centering\small\textbf{Class}\strut}} &
   \multicolumn{1}{|c|}{\parbox{6cm}{\strut\centering\small\textbf{$\pi \rtimes \sigma$ character value}\strut}} \\
\hline
\hline
\endfirsthead

\multicolumn{2}{c}{{\tablename} \thetable{} -- Continued} \\[0.5ex]
\hline
   \multicolumn{1}{|c|}{\parbox{1cm}{\strut\centering\small\textbf{Class}\strut}} &
   \multicolumn{1}{|c|}{\parbox{6cm}{\strut\centering\small\textbf{$\pi \rtimes \sigma$ character value}\strut}} \\
\hline
\hline
\endhead

 \hline
\endlastfoot

$A_1(k)$ & $(q^2+1)(q+1)\sigma(\gamma^{2k})\omega_\pi(\gamma^k)$ \\
\hline
$A_2(k)$ & $(q+1)\sigma(\gamma^{2k})\omega_\pi(\gamma^k)$ \\
\hline
$A_{31}(k)$ & $\sigma(\gamma^{2k})\Bigl(\omega_\pi(\gamma^k) + 2q\chi_\pi(M(\gamma^k))\Bigr)$ \\
\hline
$A_{32}(k)$ & $\sigma(\gamma^{2k})\omega_\pi(\gamma^k)$ \\
\hline
$A_5(k)$ & $\sigma(\gamma^{2k})\chi_\pi(M(\gamma^k))$ \\
\hline
$B_{11}(k)$ & $(q+1)^2 \sigma(\gamma^{2k})\chi_\pi({\rm diag}(\gamma^k,-\gamma^k))$ \\
\hline
$B_{12}(k)$ & $(q^2+1) \sigma(\gamma^{2k+1})\chi_\pi({\rm diag}(\gamma^{k+1/2},-\gamma^{k+1/2}))$ \\
\hline
$B_{21}(k)$ & $\sigma(-\gamma^{2k})\Bigl(\omega_\pi(\gamma^k) + \omega_\pi(-\gamma^k)+ (q+1)\chi_\pi({\rm diag}(\gamma^k,-\gamma^k))\Bigr)$ \\
\end{longtable}
\end{small}
\end{center}

\subsection{Klingen}\label{Klingenreps}
Let $\chi$ be a character of $\F_q^\times$, and let $(\pi, V)$ be an irreducible representation of $\GL(2, q)$. Define a representation $\pi'$ of the Klingen parabolic subgroup $Q$ on $V$ by
\begin{displaymath}
\left(%
\right)), \, \, \, {\rm where \,} \Delta = ad-bc.
\end{displaymath}
The character of this representation is given by $\chi(t)\chi_\pi(\left(%
\right))$, where $\chi_\pi$ is the character of $\pi$. The representation ${\rm Ind}_K^G(\pi')$ will be denoted by $\chi \rtimes \pi$. If $\pi$ has central character $\omega_\pi$, then the central character of $\chi \rtimes \pi$ is $\chi \omega_\pi$. The character table of $\chi \rtimes \pi$ is Table \ref{QCharacterValues}.

\begin{center}
\begin{small}
\begin{longtable}{|c|l|}
\caption[$\chi \rtimes \pi$ character values]{\small{$\chi \rtimes \pi$ character values}} \label{QCharacterValues} \\

\hline
   \multicolumn{1}{|c|}{\parbox{1cm}{\strut\centering\small\textbf{Class}\strut}} &
   \multicolumn{1}{|c|}{\parbox{6cm}{\strut\centering\small\textbf{$\chi \rtimes \pi$ character value}\strut}} \\
\hline
\hline
\endfirsthead

\multicolumn{2}{c}{{\tablename} \thetable{} -- Continued} \\[0.5ex]
\hline
   \multicolumn{1}{|c|}{\parbox{1cm}{\strut\centering\small\textbf{Class}\strut}} &
   \multicolumn{1}{|c|}{\parbox{6cm}{\strut\centering\small\textbf{$\chi \rtimes \pi$ character value}\strut}} \\
\hline
\hline
\endhead

 \hline
\endlastfoot

$A_1(k)$ & $(q^2+1)(q+1)\chi(\gamma^k)\omega_\pi(\gamma^k)$ \\
\hline
$A_2(k)$ & 
$\chi(\gamma^k)\Bigl(\omega_\pi(\gamma^k) + q(q+1)\chi_\pi(M(\gamma^k))\Bigr)$ \\
\hline
$A_{31}(k)$ & $2\chi(\gamma^k)\Bigl(\omega_\pi(\gamma^k) + \dfrac{(q-1)}{2}\chi_\pi(M(\gamma^k))\Bigr)$ \\
\hline
$A_{32}(k)$ & $(q+1)\chi(\gamma^k)\chi_\pi(M(\gamma^k))$ \\
\hline
$A_5(k)$ & $\chi(\gamma^k)\chi_\pi(M(\gamma^k))$ \\
\hline
$B_{11}(k)$ & $(q+1)\chi(\gamma^k)\Bigl(\chi(-1)\omega_\pi(\gamma^k) + \omega_\pi(-\gamma^k)\Bigr)$ \\
\hline
$B_{21}(k)$ & $(q+1)\chi(\gamma^k)\Bigl(1 + \chi(-1)\Bigr)\chi_\pi({\rm diag}(\gamma^k,-\gamma^k))$ \\
\hline
$B_3(k)$ & $\chi(\gamma^k)\Bigl(\chi(-1)\omega_\pi(\gamma^k)+ (q+1)\chi_\pi(M(-\gamma^k))\Bigr)$ \\
\end{longtable}
\end{small}
\end{center}

\section{Irreducible characters}
\label{Irredchar}
All of the non-trivial irreducible characters of $\Sp(4,q)$ were determined by Srinivasan, \cite{Sr}. Her list of characters can be used to determine the list of all of the irreducible characters of $\GSp(4,q)$ and some of their character tables. This list will help to determine the irreducible constituents of representations which are induced from the Borel, Siegel parabolic, and Klingen parabolic subgroups. Such constituents are precisely the \emph{non-cuspidal} representations of $\GSp(4,q)$.

The list of all of the nontrivial irreducible characters of $\GSp(4,q)$ is given below in terms of Srinivasan's list of irreducible characters of $\Sp(4,q)$. Note that we use the same monomorphism from $\F_{q^4}^\times$ into $\C^\times$ that we chose at the beginning of this paper to determine values of characters of $\Sp(4,q)$. See \cite{Sr} for the definitions\footnote{There is a misprint in \cite{Sr}, page 523. Note that $\dim(\theta_1)$ is indeed $\frac{1}{2}q^2(q^2+1).$}  of the characters $\chi_1,\dots, \chi_9, \xi_1,\dots,\xi_{41}',\Phi_1,\dots,\Phi_9,$ and $\theta_1,\dots,\theta_{13}$. 

\subsection{Constructing characters of $\GSp(4,q)$}
Define the group $\GSp(4,q)^+:=Z\, \cdot \Sp(4,q)$. Let $\chi$ be an irreducible character of $\Sp(4,q)$ and let $\alpha$ be a character of $\F_q^\times\cong Z$ such that $\alpha(-1)=\chi(-I_4)$. We extend the character $\chi$ to an irreducible character of $\GSp(4,q)^+$ by defining $\alpha\chi(zg):=\alpha(z)\chi(g)$, for $z\in \F_q\cong Z$ and $g\in\Sp(4,q)$. Since $\alpha(-1)=\chi(-I_4)$, we have that $\alpha\chi$ is well-defined. The character of $\GSp(4,q)$ induced from $\alpha\chi$ is denoted by ${\rm Ind}(\alpha\chi)$ or simply by ${\rm Ind}(\chi)$ if $\alpha$ is the trivial character.

\begin{thm}\label{IrredCharacterThm}
All irreducible characters of $\GSp(4,q)$ are listed in Table \ref{GSpIrredCharacters}. The characters in the table are pairwise inequivalent.
\end{thm}

In each row in the table, $\alpha$ is any character of $\F_q^\times\cong Z$ with the specified value of $\alpha(-1)$. Genericity of a character is indicated by a bullet, $\bullet$, in the ``g" column. The abbreviation $t=\frac{1}{2}(q-1)$ is also used. In the cases where the induced character decomposes, say into $\chi_a$ and $\chi_b$, we have that $\chi_b = \xi\chi_a$, where $\xi:\F_q^\times\longrightarrow\C^\times$ is the unique non-trivial quadratic character defined by $\xi(\gamma)=-1$. Note that there are more characters of $\Sp(4,q)$ than those given explicitly in Srinivasan's tables, \cite{Sr}. The missing characters are the trivial character, $\xi_{22}, \xi'_{22}, \xi_{42}, \xi'_{42}, \phi_i,$ and $\theta_i$ for $i=2,4,6,8$. Srinivasan explains how to compute the missing non-trivial characters using the tables in her paper. Note that we do not obtain any new characters of $\GSp(4,q)$ with these characters, because ${\rm Ind}(\alpha\xi_{21})={\rm Ind}(\alpha\xi_{22})$, ${\rm Ind}(\alpha\xi'_{21})={\rm Ind}(\alpha\xi'_{22})$, ${\rm Ind}(\alpha\xi_{41})={\rm Ind}(\alpha\xi_{42})$, ${\rm Ind}(\alpha\xi'_{41})={\rm Ind}(\alpha\xi'_{42})$, ${\rm Ind}(\alpha\phi_i)={\rm Ind}(\alpha\phi_{i+1})$, and ${\rm Ind}(\alpha\theta_{i})={\rm Ind}(\alpha\theta_{i+1})$  for $i=2,4,6,8.$

\begin{center}
\begin{small}
\begin{longtable}{|c|c|l|l|c|}
\caption[Irreducible characters of $\GSp(4,q)$]{\small{Irreducible characters of $\GSp(4,q)$}} \label{GSpIrredCharacters} \\

\hline
   \multicolumn{1}{|c|}{\parbox{2cm}{\strut\centering\small\textbf{Character}\strut}} &
   \multicolumn{1}{|c|}{\parbox{1cm}{\strut\centering\small\textbf{$\alpha(-1)$}\strut}} &
   \multicolumn{1}{|c|}{\parbox{2.5cm}{\strut\centering\small\textbf{Constituents}\strut}} &
   \multicolumn{1}{|c|}{\parbox{2.5cm}{\strut\centering\small\textbf{Dimension}\strut}} &
   \multicolumn{1}{|c|}{\parbox{.3cm}{\strut\centering\small\textbf{g}\strut}} \\
\hline
\hline
\endfirsthead

\multicolumn{5}{c}{{\tablename} \thetable{} -- Continued} \\[0.5ex]
\hline
   \multicolumn{1}{|c|}{\parbox{2cm}{\strut\centering\small\textbf{Character}\strut}} &
   \multicolumn{1}{|c|}{\parbox{1cm}{\strut\centering\small\textbf{$\alpha(-1)$}\strut}} &
   \multicolumn{1}{|c|}{\parbox{2.5cm}{\strut\centering\small\textbf{Constituents}\strut}} &
   \multicolumn{1}{|c|}{\parbox{2.5cm}{\strut\centering\small\textbf{Dimension}\strut}} &
   \multicolumn{1}{|c|}{\parbox{.3cm}{\strut\centering\small\textbf{g}\strut}} \\
\hline
\hline
\endhead

 \hline
\endlastfoot

$\alpha\triv_{\GSp(4,q)}$ & no condition & $\alpha\triv_{\GSp(4,q)}$ & $1$ & \\
\hline
${\rm Ind}(\alpha\chi_1(n))$ & & ${\rm Ind}(\alpha\chi_1(n))_a$ & $(q^2-1)^2$ & $\bullet$ \\
\cline{3-5}
\raisebox{2ex}{\scriptsize{$n\in R_1$}} &\raisebox{2ex}{$(-1)^n$} & ${\rm Ind}(\alpha\chi_1(n))_b$ & $(q^2-1)^2$ & $\bullet$ \\
\hline
${\rm Ind}(\alpha\chi_2(n))$ & & ${\rm Ind}(\alpha\chi_2(n))_a$ & $q^4-1$ & $\bullet$ \\
\cline{3-5}
\raisebox{2ex}{\scriptsize{$n\in R_2$}} &\raisebox{2ex}{$(-1)^n$} & ${\rm Ind}(\alpha\chi_2(n))_b$ & $q^4-1$ & $\bullet$ \\
\hline
${\rm Ind}(\alpha\chi_3(n,m))$ & & ${\rm Ind}(\alpha\chi_3(n,m))_a$ & $(q^2+1)(q+1)^2$ & $\bullet$ \\
\cline{3-5}
\raisebox{2ex}{\scriptsize{$n,m\in T_1, n<m$}} &\raisebox{2ex}{$(-1)^{n+m}$} & ${\rm Ind}(\alpha\chi_3(n,m))_b$ & $(q^2+1)(q+1)^2$ & $\bullet$ \\
\hline
${\rm Ind}(\alpha\chi_4(n,m))$ & & ${\rm Ind}(\alpha\chi_4(n,m))_a$ & $(q^2+1)(q-1)^2$ & $\bullet$ \\
\cline{3-5}
\raisebox{2ex}{\scriptsize{$n,m\in T_2, n\neq m$}} &\raisebox{2ex}{$(-1)^{n+m}$} & ${\rm Ind}(\alpha\chi_4(n,m))_b$ & $(q^2+1)(q-1)^2$ & $\bullet$ \\
\hline
${\rm Ind}(\alpha\chi_5(n,m))$ & & ${\rm Ind}(\alpha\chi_5(n,m))_a$ & $q^4-1$ & $\bullet$ \\
\cline{3-5}
\raisebox{2ex}{\scriptsize{$n\in T_2, m\in T_1$}} &\raisebox{2ex}{$(-1)^{n+m}$} & ${\rm Ind}(\alpha\chi_5(n,m))_b$ & $q^4-1$ & $\bullet$ \\
\hline
${\rm Ind}(\alpha\chi_6(n))$ & & ${\rm Ind}(\alpha\chi_6(n))_a$ & $(q^2+1)(q-1)$ & \\
\cline{3-5}
\raisebox{2ex}{\scriptsize{$n\in T_2$}} &\raisebox{2ex}{$1$} & ${\rm Ind}(\alpha\chi_6(n))_b$ & $(q^2+1)(q-1)$ & \\
\hline
${\rm Ind}(\alpha\chi_7(n))$ & & ${\rm Ind}(\alpha\chi_7(n))_a$ & $q(q^2+1)(q-1)$ & $\bullet$ \\
\cline{3-5}
\raisebox{2ex}{\scriptsize{$n\in T_2$}} &\raisebox{2ex}{$1$} & ${\rm Ind}(\alpha\chi_7(n))_b$ & $q(q^2+1)(q-1)$ & $\bullet$ \\
\hline
${\rm Ind}(\alpha\chi_8(n))$ & & ${\rm Ind}(\alpha\chi_8(n))_a$ & $(q^2+1)(q+1)$ & \\
\cline{3-5}
\raisebox{2ex}{\scriptsize{$n\in T_1$}} &\raisebox{2ex}{$1$} & ${\rm Ind}(\alpha\chi_8(n))_b$ & $(q^2+1)(q+1)$ & \\
\hline
${\rm Ind}(\alpha\chi_9(n))$ & & ${\rm Ind}(\alpha\chi_9(n))_a$ & $q(q^2+1)(q+1)$ & $\bullet$ \\
\cline{3-5}
\raisebox{2ex}{\scriptsize{$n\in T_1$}} &\raisebox{2ex}{$1$} & ${\rm Ind}(\alpha\chi_9(n))_b$ & $q(q^2+1)(q+1)$ & $\bullet$ \\
\hline

${\rm Ind}(\alpha\xi_1(n))$ & & ${\rm Ind}(\alpha\xi_1(n))_a$ & $(q^2+1)(q-1)$ & \\
\cline{3-5}
\raisebox{2ex}{\scriptsize{$n\in T_2$}} &\raisebox{2ex}{$(-1)^n$} & ${\rm Ind}(\alpha\xi_1(n))_b$ & $(q^2+1)(q-1)$ & \\
\hline
${\rm Ind}(\alpha\xi_1'(n))$ & & ${\rm Ind}(\alpha\xi_1'(n))_a$ & $q(q^2+1)(q-1)$ & $\bullet$ \\
\cline{3-5}
\raisebox{2ex}{\scriptsize{$n\in T_2$}} &\raisebox{2ex}{$(-1)^n$} & ${\rm Ind}(\alpha\xi_1'(n))_b$ & $q(q^2+1)(q-1)$ & $\bullet$ \\
\hline
\begin{minipage}{11ex}
\begin{center}
\vspace{0.05in}
${\rm Ind}(\alpha\xi_{21}(n))$\\
\scriptsize{$n\in T_2$}
\vspace{0.05in}
\end{center}
\end{minipage} & $(-1)^{n+t}$ & ${\rm Ind}(\alpha\xi_{21}(n))$ & $q^4-1$ & $\bullet$ \\
\hline
\begin{minipage}{11ex}
\begin{center}
\vspace{0.05in}
${\rm Ind}(\alpha\xi_{21}'(n))$\\
\scriptsize{$n\in T_2$}
\vspace{0.05in}
\end{center}
\end{minipage} & $(-1)^{n+t+1}$ & ${\rm Ind}(\alpha\xi_{21}'(n))$ & $(q^2+1)(q-1)^2$ & $\bullet$ \\
\hline

${\rm Ind}(\alpha\xi_3(n))$ & & ${\rm Ind}(\alpha\xi_3(n))_a$ & $(q^2+1)(q+1)$ & \\
\cline{3-5}
\raisebox{2ex}{\scriptsize{$n\in T_1$}} &\raisebox{2ex}{$(-1)^n$} & ${\rm Ind}(\alpha\xi_3(n))_b$ & $(q^2+1)(q+1)$ & \\
\hline
${\rm Ind}(\alpha\xi_3'(n))$ & & ${\rm Ind}(\alpha\xi_3'(n))_a$ & $q(q^2+1)(q+1)$ & $\bullet$ \\
\cline{3-5}
\raisebox{2ex}{\scriptsize{$n\in T_1$}} &\raisebox{2ex}{$(-1)^n$} & ${\rm Ind}(\alpha\xi_3'(n))_b$ & $q(q^2+1)(q+1)$ & $\bullet$ \\
\hline
\begin{minipage}{11ex}
\begin{center}
\vspace{0.05in}
${\rm Ind}(\alpha\xi_{41}(n))$\\
\scriptsize{$n\in T_1$}
\vspace{0.05in}
\end{center}
\end{minipage} & $(-1)^{n+t}$ & ${\rm Ind}(\alpha\xi_{41}(n))$ & $(q^2+1)(q+1)^2$ & $\bullet$ \\
\hline
\begin{minipage}{11ex}
\begin{center}
\vspace{0.05in}
${\rm Ind}(\alpha\xi_{41}'(n))$\\
\scriptsize{$n\in T_1$}
\vspace{0.05in}
\end{center}
\end{minipage} & $(-1)^{n+t+1}$ & ${\rm Ind}(\alpha\xi_{41}'(n))$ & $q^4-1$ & $\bullet$ \\
\hline
${\rm Ind}(\alpha\Phi_1)$ & $(-1)^{t+1}$ & ${\rm Ind}(\alpha\Phi_1)$ & $(q^2+1)(q-1)$ & \\
\hline
${\rm Ind}(\alpha\Phi_3)$ & $(-1)^{t+1}$ & ${\rm Ind}(\alpha\Phi_3)$ & $q(q^2+1)(q-1)$ & $\bullet$ \\
\hline
${\rm Ind}(\alpha\Phi_5)$ & $(-1)^t$ & ${\rm Ind}(\alpha\Phi_5)$ & $(q^2+1)(q+1)$ & \\
\hline
${\rm Ind}(\alpha\Phi_7)$ & $(-1)^t$ & ${\rm Ind}(\alpha\Phi_7)$ & $q(q^2+1)(q+1)$ & $\bullet$ \\
\hline
& & ${\rm Ind}(\alpha\Phi_9)_a$ & $q(q^2+1)$ & \\
\cline{3-5}
\raisebox{2ex}{${\rm Ind}(\alpha\Phi_9)$}&\raisebox{2ex}{$1$} & ${\rm Ind}(\alpha\Phi_9)_b$ & $q(q^2+1)$ & \\
\hline
\hline
${\rm Ind}(\alpha\theta_1)$ & $1$ & ${\rm Ind}(\alpha\theta_1)$ & $q^2(q^2+1)$ & $\bullet$ \\
\hline
${\rm Ind}(\alpha\theta_3)$ & $1$ & ${\rm Ind}(\alpha\theta_3)$ & $q^2+1$ & \\
\hline
${\rm Ind}(\alpha\theta_5)$ & $-1$ & ${\rm Ind}(\alpha\theta_5)$ & $q^2(q^2-1)$ & $\bullet$ \\
\hline
${\rm Ind}(\alpha\theta_7)$ & $-1$ & ${\rm Ind}(\alpha\theta_7)$ & $q^2-1$ & \\
\hline

& & ${\rm Ind}(\alpha\theta_9)_a$ & $\frac{1}{2}q(q+1)^2$ & \\
\cline{3-5}
\raisebox{2ex}{${\rm Ind}(\alpha\theta_9)$}&\raisebox{2ex}{$1$} & ${\rm Ind}(\alpha\theta_9)_b$ & $\frac{1}{2}q(q+1)^2$ & \\
\hline
& & ${\rm Ind}(\alpha\theta_{10})_a$ & $\frac{1}{2}q(q-1)^2$ & \\
\cline{3-5}
\raisebox{2ex}{${\rm Ind}(\alpha\theta_{10})$}&\raisebox{2ex}{$1$} & ${\rm Ind}(\alpha\theta_{10})_b$ & $\frac{1}{2}q(q-1)^2$ & \\
\hline
& & ${\rm Ind}(\alpha\theta_{11})_a$ & $\frac{1}{2}q(q^2+1)$ & \\
\cline{3-5}
\raisebox{2ex}{${\rm Ind}(\alpha\theta_{11})$}&\raisebox{2ex}{$1$} & ${\rm Ind}(\alpha\theta_{11})_b$ & $\frac{1}{2}q(q^2+1)$ & \\
\hline
& & ${\rm Ind}(\alpha\theta_{12})_a$ & $\frac{1}{2}q(q^2+1)$ & \\
\cline{3-5}
\raisebox{2ex}{${\rm Ind}(\alpha\theta_{12})$}&\raisebox{2ex}{$1$} & ${\rm Ind}(\alpha\theta_{12})_b$ & $\frac{1}{2}q(q^2+1)$ & \\
\hline
& & ${\rm Ind}(\alpha\theta_{13})_a$ & $q^4$ & $\bullet$ \\
\cline{3-5}
\raisebox{2ex}{${\rm Ind}(\alpha\theta_{13})$}&\raisebox{2ex}{$1$} & ${\rm Ind}(\alpha\theta_{13})_b$ & $q^4$ & $\bullet$ \\
\end{longtable}
\end{small}
\end{center}

Before we begin the proof of Theorem \ref{IrredCharacterThm}, we note the following useful result from Fulton-Harris, \cite{Fulton}. Let $H$ be a subgroup of index 2 in a group $G$. Then $H$ is a normal subgroup of $G$ and $G/H$ is a group of order 2. Let $U$ and $U'$ denote the trivial and non-trivial representations of $G$, respectively, obtained from the two irreducible representations of $G/H$. For any representation $V$ of $G$, let $V' = V \otimes U'$. The character of $V'$ is the same as the character of $V$ on elements of $H$, but takes opposite values on elements not in $H$. In particular, ${\rm Res}_H^G V' = {\rm Res}_H^G V$. If $W$ is any representation of $H$, there is a \emph{conjugate} representation defined by conjugating by any element $t\in G$ such that $t\notin H$; if $\psi$ is the character of $W$, the character of the conjugate is $h\mapsto \psi(tht^{-1})$. Since $t$ is unique up to multiplication by an element of $H$, the conjugate representation is unique up to isomorphism.

\begin{prop}[\cite{Fulton}, Proposition 5.1]\label{index2prop}
Let $V$ be an irreducible representation of $G$. Let $W={\rm Res}_H^G V$ be the restriction of $V$ to $H$. Then exactly one of the following holds:

{\rm (1)} $V$ is not isomorphic to $V'$; $W$ is irreducible and isomorphic to its conjugate; ${\rm Ind}_H^G W \cong V\oplus V'$.

{\rm (2)} $V\cong V'$; $W=W'\oplus W''$, where $W'$ and $W''$ are irreducible and conjugate but not isomorphic; ${\rm Ind}_H^G W'\cong {\rm Ind}_H^G W''\cong V$.

Each irreducible representation of $H$ arises uniquely in this way, noting that in case {\rm (1)} the restrictions of $V$ and $V'$ to $H$ determine the same representation.
\end{prop}

\begin{proof}
{\em of Theorem \ref{IrredCharacterThm}.} Let $\pi$ be an irreducible representation of $\Sp(4,q)$. Let $\alpha$ be a character of $Z$. Define a representation $\pi^+$ of $\GSp(4,q)^+$ by 
$$\pi^+(z\cdot g) := \alpha(z)\pi(g).$$
By Schur's Lemma, elements of $Z$ act as scalars on vectors in the space of $\pi$. It is clear that this representation is well-defined if and only if $\alpha(\pm 1)$ acts as $\pi(\pm I)$ on the space of $\pi$. The character of this new representation is $\alpha(z)\chi_\pi(g)$, where $\chi_\pi$ is the character of $\pi$. By taking the inner product of the character with itself, one sees that this representation is irreducible.

Consider the induced representation ${\rm Ind}_{\GSp(4,q)^+}^{\GSp(4,q)}(\pi^+).$ By Proposition \ref{index2prop}, the induced representation is either irreducible or it has precisely two irreducible constituents, since $\GSp(4,q)^+$ is an index two subgroup of $\GSp(4,q)$. Moreover, the induced representation has precisely two irreducible constituents if and only if
$$\chi_\pi(g)=\chi_\pi(\left(\begin{smallmatrix}
1&&&\\
&1&&\\
&&\gamma&\\
&&&\gamma\\
\end{smallmatrix}\right)
g
\left(\begin{smallmatrix}
1&&&\\
&1&&\\
&&\gamma^{-1}&\\
&&&\gamma^{-1}\\
\end{smallmatrix}\right))=\overline{\chi_\pi(g)}$$
for all $g\in \Sp(4,q)$. Otherwise, it is irreducible. So one only needs to check which characters of $\Sp(4,q)$ are equal to their conjugate. For those characters $\chi$ such that $\chi\neq\overline\chi$, both $\chi$ and $\overline{\chi}$ yield the same irreducible representation of $\GSp(4,q)$ via the induction of the extension of these characters to $\GSp(4,q)^+$. In total, $(q^2+2q+4)(q-1)$ unique irreducible characters are obtained in this manner, which is the same amount as the number of conjugacy classes of $\GSp(4,q)$. This completes the proof.
\qed
\end{proof}

The computation of the character values was performed using the values from Srinivasan's character tables, \cite{Sr}. The conjugacy classes we obtained for $\GSp(4,q)$ are, due to the manner of their computation, organized in a way which is quite different from the classes of $\Sp(4,q)$. This naturally presented some problems in bookkeeping. Moreover, Srinivasan defined $\Sp(4,q)$ using a different $J$ matrix, which she calls $A$. The classes of $\Sp(4,q)$ with respect to our $J$ matrix will be denoted using Srinivasan's notation: $A_n,B_n,B_n',C_n(i),C_6(i,j),D_n(i)$. These are the images of Srinivasan's conjugacy classes under the isomorphism that sends an element $g$ to
$$\left(\begin{smallmatrix}
 1 & & & \\
 & & & 1 \\
 & & -1 & \\
 & 1 & & \\
\end{smallmatrix}\right)\cdot g\cdot \left(\begin{smallmatrix}
 1 & & & \\
 & & & 1 \\
 & & -1 & \\
 & 1 & & \\
\end{smallmatrix}\right).$$

In Table \ref{Indalphachivalues}, we describe how to compute character values of the induced characters ${\rm Ind}_{\GSp(4,q)^+}^{\GSp(4,q)}(\chi^+),$ where $\chi^+=\alpha\chi$ with $\chi$ a character of $\Sp(4,q)$ and $\alpha$ a character of the center $Z$ of $\GSp(4,q)$.

\begin{center}
\begin{small}
\begin{longtable}{|c|l|}
\caption[${\rm Ind}(\alpha\chi)$ character values]{${\rm Ind}(\alpha\chi)$ character values} \label{Indalphachivalues} \\

\hline
   \multicolumn{1}{|c|}{\parbox{2cm}{\strut\centering\small\textbf{Class}\strut}} &
   \multicolumn{1}{|c|}{\parbox{8cm}{\strut\centering\small\textbf{${\rm Ind}(\alpha\chi)$ character value}\strut}} \\
\hline
\endfirsthead

\multicolumn{2}{c}{{\tablename} \thetable{} -- Continued} \\[0.5ex]
\hline
   \multicolumn{1}{|c|}{\parbox{2cm}{\strut\centering\small\textbf{Class}\strut}} &
   \multicolumn{1}{|c|}{\parbox{8cm}{\strut\centering\small\textbf{${\rm Ind}(\alpha\chi)$ character value}\strut}} \\
\hline
\endhead

 \hline
\endlastfoot

$A_1(k)$ & $2\chi(A_1)\alpha(\gamma^k)$\\
\hline
$A_2(k)$ & $(\chi(A_{21})+\chi(A_{22}))\alpha(\gamma^k)$\\
\hline
$A_{31}(k)$ & $2\chi(A_{31})\alpha(\gamma^k) $\\
\hline
$A_{32}(k)$ & $2\chi(A_{32})\alpha(\gamma^k) $\\
\hline
$A_5(k)$ & $2(\chi(A_{41})+\chi(A_{42}))\alpha(\gamma^k)$\\
\hline
$B_{11}(k)$ & $2\chi(D_1)\alpha(\gamma^k) $\\
\hline
$B_{12}(k)$ & 0\\
\hline
$B_{21}(k)$ & \begin{minipage}{30ex}
\vspace{0.05in}
$\begin{cases}
2\chi(B_8(\frac{1-q}{4}))\alpha(\gamma^{k+\frac{q-1}{4}}) & {\rm if}\, q\equiv 1 (4)\\
0 & {\rm if}\, q\equiv 3 (4)\\
\end{cases}$
\vspace{0.01in}
\end{minipage}\\
\hline
$B_{22}(k)$ &  \begin{minipage}{30ex}
\vspace{0.05in}
$\begin{cases}
0 & {\rm if}\, q\equiv 1 (4)\\
2\chi(B_6(\frac{q+1}{4}))\alpha(\gamma^{k+\frac{q+1}{4}}) & {\rm if}\, q\equiv 3 (4)\\
\end{cases}$
\vspace{0.01in}
\end{minipage}\\
\hline
$B_3(k)$ & $(\chi(D_{21})+\chi(D_{22}))\alpha(\gamma^k)$\\
\hline
$B_{41}(k)$ & $(\chi(D_{31})+\chi(D_{34}))\alpha(\gamma^k)$\\
\hline
$B_{42}(k)$ & $(\chi(D_{32})+\chi(D_{33}))\alpha(\gamma^k)$\\
\hline
$B_{43}(k)$ & 0\\
\hline
$B_{44}(k)$ & 0\\
\hline
$B_{51}(k)$ & \begin{minipage}{30ex}
\vspace{0.05in}
$\begin{cases}
2\chi(B_9(\frac{1-q}{4}))\alpha(\gamma^{k+\frac{q-1}{4}}) & {\rm if}\, q\equiv 1 (4)\\
0 & {\rm if}\, q\equiv 3 (4)\\
\end{cases}$
\vspace{0.01in}
\end{minipage}\\
\hline
$B_{52}(k)$ & \begin{minipage}{30ex}
\vspace{0.05in}
$\begin{cases}
0 & {\rm if}\, q\equiv 1 (4)\\
2\chi(B_7(\frac{q+1}{4}))\alpha(\gamma^{k+\frac{q+1}{4}}) & {\rm if}\, q\equiv 3 (4)\\
\end{cases}$
\vspace{0.01in}
\end{minipage}\\
\hline
$C_1(i,k)$ & \begin{minipage}{30ex}
\vspace{0.05in}
$\begin{cases}
0 & {\rm if}\, i\, {\rm is\, odd}\\
2\chi(B_8(-\frac{i}{2}))\alpha(\gamma^{k+\frac{i}{2}}) & {\rm if}\, i\, {\rm is\, even}\\
\end{cases}$
\vspace{0.01in}
\end{minipage}\\
\hline
$C_{21}(i,k)$ & \begin{minipage}{30ex}
\vspace{0.05in}
$\begin{cases}
2\chi(B_3(\frac{-i-t}{2},\frac{t-i}{2}))\alpha(\gamma^{\frac{2k+i+t}{2}}) & {\rm if}\, q\equiv 1 (4)\, {\rm and}\,  i\, {\rm is\, even}\\
2\chi(B_3(\frac{-i-t}{2},\frac{t-i}{2}))\alpha(\gamma^{\frac{2k+i+t}{2}}) & {\rm if}\, q\equiv 3 (4)\, {\rm and}\, i\, {\rm is\, odd}\\
0 & {\rm otherwise}\\
\end{cases}$
\vspace{0.01in}
\end{minipage}\\
\hline
$C_{22}(i,k)$ & \begin{minipage}{30ex}
\vspace{0.05in}
$\begin{cases}
2\chi(B_2((q+1)(\frac{t-i}{2})))\alpha(\gamma^{k+\frac{q+1}{4}+\frac{i}{2}}) & {\rm if}\, q\equiv 1 (4)\\
& {\rm and}\, i\, {\rm odd}\\
2\chi(B_2((q+1)(\frac{t-i}{2})))\alpha(\gamma^{k+\frac{q+1}{4}+\frac{i}{2}}) & {\rm if}\, q\equiv 3 (4)\\
& {\rm and}\, i\, {\rm even}\\
0 & {\rm otherwise}\\
\end{cases}$
\vspace{0.01in}
\end{minipage}\\
\hline
$C_3(i,k)$ & \begin{minipage}{30ex}
\vspace{0.05in}
$\begin{cases}
0 & {\rm if}\, i\, {\rm is\, odd}\\
2\chi(B_9(-\frac{i}{2}))\alpha(\gamma^{k+\frac{i}{2}}) & {\rm if}\, i\, {\rm is\, even}\\
\end{cases}$
\vspace{0.01in}
\end{minipage}\\
\hline
$C_4(i,k)$ & $(\chi(C_{41}(i))+\chi(C_{42}(i)))\alpha(\gamma^k)$\\
\hline
$C_5(i,k)$ & $2\chi(C_3(i))\alpha(\gamma^{k})$\\
\hline
$C_6(i,j,k)$ & \begin{minipage}{30ex}
\vspace{0.05in}
$\begin{cases}
0 & {\rm if}\, i+j\, {\rm is\, odd}\\
2\chi(B_3(\frac{-i-j}{2},\frac{j-i}{2}))\alpha(\gamma^{k+\frac{i+j}{2}}) & {\rm if}\, i+j\, {\rm is\, even}\\
\end{cases}$
\vspace{0.01in}
\end{minipage}\\
\hline
$D_1(i,k)$ & \begin{minipage}{30ex}
\vspace{0.05in}
$\begin{cases}
0 & {\rm if}\, i\, {\rm is\, odd}\\
2\chi(B_5(-\frac{i}{2},-\frac{i}{2}))\alpha(\gamma^{k+\frac{i}{2}}) & {\rm if}\, i\, {\rm is\, even}\\
\end{cases}$
\vspace{0.01in}
\end{minipage}\\
\hline
$D_2(i,k)$ & $\chi(B_6(i))\alpha(\gamma^k)$\\
\hline
$D_{32}(i,k)$ & 0\\
\hline
$D_5(i,k)$ & $\chi(B_7(i))\alpha(\gamma^k)$\\
\hline
$D_6(i,k)$ & $2\chi(C_1(i))\alpha(\gamma^k)$\\
\hline
$D_8(i,k)$ & $(\chi(C_{21}(i))+\chi(C_{22}(i)))\alpha(\gamma^k)$\\
\end{longtable}
\end{small}
\end{center}

We omit the values of ${\rm Ind}(\alpha\chi)$ on the classes $D_{31}(i,k), D_4(i,j,k)$, $D_7(i,j,k),$ and $D_9(i,k)$ since they turn out to not be necessary for our determination of the irreducible non-cuspidal characters.

Note that the group $\GSp(4,q)$$^+$ only has elements with square multipliers and so the induced character takes the value 0 on the non-square multiplier classes of $\GSp(4,q)$ by the induced character formula. If the induced character decomposes into two constituents, then we know that the sum of the values of the constituent characters on the non-square multiplier classes is 0, but the particular values each constituent takes on these classes are unknown. However, we can say that the values of the constituent characters on the square multiplier classes are half the values of the induced character on those classes.

We now give a correspondence between our notations and Shinoda's notations for the irreducible characters of $\GSp(4,q)$. Here, $\omega_0$ is the non-trivial quadratic character on $<\eta>$ and $\Lambda_1$ is the character $\Lambda_1:\F_{q^2}^\times\to\C$ such that $\Lambda_1(\theta)=\omega_0(\theta^{q-1})=\omega_0(\eta)$.

\begin{center}
\renewcommand{\arraystretch}{1.3}
\begin{small}
\begin{longtable}{|l|l||l|l|}
\caption[Notations for irreducible characters of $\GSp(4,q)$]{\small{Notations for irreducible characters of $\GSp(4,q)$}} \label{GSpIrreducibleCharactersShinoda} \\

\hline
   \multicolumn{1}{|c|}{\parbox{1.2cm}{\strut\centering\small\textbf{Shinoda}\strut}} &
   \multicolumn{1}{|c|}{\parbox{2cm}{\strut\centering\small\textbf{Our notation}\strut}} &
   \multicolumn{1}{|c|}{\parbox{1.2cm}{\strut\centering\small\textbf{Shinoda}\strut}} &
   \multicolumn{1}{|c|}{\parbox{2cm}{\strut\centering\small\textbf{Our notation}\strut}} \\
\hline
\hline
\endfirsthead

\multicolumn{3}{c}{{\tablename} \thetable{} -- Continued} \\[0.5ex]
\hline
   \multicolumn{1}{|c|}{\parbox{1.2cm}{\strut\centering\small\textbf{Shinoda}\strut}} &
   \multicolumn{1}{|c|}{\parbox{2cm}{\strut\centering\small\textbf{Our notation}\strut}} &
   \multicolumn{1}{|c|}{\parbox{1.2cm}{\strut\centering\small\textbf{Shinoda}\strut}} &
   \multicolumn{1}{|c|}{\parbox{2cm}{\strut\centering\small\textbf{Our notation}\strut}} \\
\hline
\hline
\endhead

 \hline
\endlastfoot

$X_1(\lambda,\mu),$ & ${\rm Ind}(\lambda\mu\cdot\chi_3(n,m))_a,$ & $\chi_6(\omega),$ & ${\rm Ind}(\chi_7(n))_a$\\
$X_1(\lambda,\mu)\neq \xi X_1(\lambda,\mu)$ & & $\omega(-1)=(-1)^n$ & \\
\hline
$X_1(\lambda,\mu),$ & ${\rm Ind}(\lambda\mu\cdot\xi_{41}(n))$ & $\chi_7(\Lambda),$ & ${\rm Ind}(\Lambda|_{\F_q^\times}\xi_1(n))_a$ \\
$X_1(\lambda,\mu)=\xi X_1(\lambda,\mu)$ && $\Lambda\neq\Lambda_1$ & \\
\hline
$X_2(\Lambda,1)$ & ${\rm Ind}(\Lambda|_{\F_q^\times}\cdot\chi_2(n))_a,$ & $\chi_7(\Lambda_1)$ & ${\rm Ind}(\Lambda_1|_{\F_q^\times}\cdot\Phi_1)$ \\

\hline
$X_3(\Lambda,\nu), \nu\neq\xi$ & ${\rm Ind}(\Lambda|_{\F_q^\times}\nu\cdot\chi_5(n,m))_a$ & $\chi_8(\Lambda),$ & ${\rm Ind}(\Lambda|_{\F_q^\times}\xi_1'(n))_a$ \\
 & & $\Lambda\neq\Lambda_1$ &\\\hline
$X_3(\Lambda,\xi), \Lambda\neq \Lambda_1$ & ${\rm Ind}(\Lambda|_{\F_q^\times}\xi\cdot\xi_{21}(n))$ & $\chi_8(\Lambda_1)$ & ${\rm Ind}(\Lambda_1|_{\F_q^\times}\Phi_3)$ \\
\hline
$X_3(\Lambda_1,\xi)$ & ${\rm Ind}(\Lambda_1|_{\F_q^\times}\xi\cdot\xi_{41}'(n))$ & $\tau_1$ & ${\rm Ind}(\theta_3)$ \\
\hline
$X_4(\Theta)$ & ${\rm Ind}(\Theta|_{\F_q^\times}\chi_1(n))_a$ & $\tau_2$ & ${\rm Ind}(\Phi_9)_a$ \\
$\omega(-1)=(-1)^n$ and &&&\\
$\Lambda$ not quadratic & & &\\
\hline
$X_5(\Lambda,\omega),$ & ${\rm Ind}(\Lambda|_{\F_q^\times}\chi_4(n,m))_a$ & $\tau_3$ & ${\rm Ind}(\theta_1)$  \\
$\omega(-1)=(-1)^n$ and &&&\\
$1\neq\Lambda$ quadratic & & &\\
\hline
$X_5(\Lambda,\omega),$ & ${\rm Ind}(\Lambda|_{\F_q^\times}\xi_{21}'(n))$ & $\tau_4(\lambda')$ &${\rm Ind}(\lambda'|_{\F_q^\times}\theta_7)$ \\\hline
$\chi_1(\lambda), \lambda\neq\xi$ & ${\rm Ind}(\lambda\cdot\xi_3(n))_a$ & $\tau_5(\lambda')$ & ${\rm Ind}(\lambda'|_{\F_q^\times}\theta_5)$ \\
\hline
$\chi_1(\xi)$ & ${\rm Ind}(\xi\cdot\Phi_5)$ & $\theta_1$ & ${\rm Ind}(\theta_9)_a$ \\
\hline
$\chi_2(\lambda), \lambda\neq\xi$ & ${\rm Ind}(\lambda\cdot\xi_3'(n))_a$ & $\theta_2$ & ${\rm Ind}(\theta_{10})_a$ \\
\hline
$\chi_2(\xi)$ & ${\rm Ind}(\xi\cdot\Phi_7) $ & $\theta_3$ & ${\rm Ind}(\theta_{11})_a$ \\
\hline
$\chi_3(\lambda)$ & ${\rm Ind}(\lambda^2\cdot\chi_8(n))_a$ & $\theta_4$ & ${\rm Ind}(\theta_{12})_a$ \\
\hline
$\chi_4(\lambda)$ & ${\rm Ind}(\lambda^2\cdot\chi_9(n))_a$ & $\theta_5$ & ${\rm Ind}(\theta_{13})_a$ \\
\hline
$\chi_5(\omega), \omega(-1)=(-1)^n$ & ${\rm Ind}(\chi_6(n))_a$ & $\theta_0(\lambda)$ & $\lambda\triv_{\GSp(4,q)}$ \\
\end{longtable}
\end{small}
\end{center}

For our families of characters with parameters $n$ or $(n,m)$, we note that the corresponding character from \cite{Shinoda} depends on our fixed embedding of $\F_{q^4}^\times$ into $\C^\times$. For example, the $n$ such that ${\rm Ind}(\lambda^2\cdot\chi_8(n))_a=\chi_3(\lambda)$ must, in particular, satisfy
$\tilde{\gamma}^{in}+\tilde{\gamma}^{-in}=\lambda(\gamma^i)+\lambda(\gamma^{-i}).$ The conditions on the parameters are suppressed in Table \ref{GSpIrreducibleCharactersShinoda}.

\subsection{Irreducible non-cuspidal representations}\label{irrednc}

The irreducible non-cuspidal representations can now be determined by decomposing the parabolically induced representations of the (standard) Borel, Siegel, and Klingen subgroups defined above into irreducible constituents. In the theorem, we also describe our characters in terms of the irreducible characters from Table \ref{GSpIrredCharacters} where possible. Precisely which character from Table \ref{GSpIrredCharacters} a particular irreducible non-cuspidal character is sometimes depends on our embedding of $\F_{q^4}^\times$ into $\C^\times$.

\begin{thm}\label{cusprep} Let $\sigma$ be a character of $\F_q^\times$. Every irreducible non-cuspidal representation of $\GSp(4,q)$ is one of the following types.

\emph{Type I.} These are the irreducible representations obtained by parabolic induction from the Borel subgroup. More precisely, they are irreducible representations of the form $\chi_1 \times \chi_2 \rtimes \sigma$, where $\chi_1$ and $\chi_2,$ are characters of $\F_q^\times$. Representations of this form are irreducible if and only if $\chi_1 \neq \triv_{\F_q^\times}, \chi_2 \neq \triv_{\F_q^\times}$ and $\chi_1 \neq \chi_2^{\pm 1}$, where $\triv_{\F_q^\times}$ is the trivial character of $\F_q^\times$. Let $\alpha=\chi_1\chi_2\sigma^2$. In terms of the irreducible characters that appear in Table \ref{GSpIrredCharacters}, these are the characters ${\rm Ind}(\alpha\chi_3(n,m))_a, {\rm Ind}(\alpha\chi_3(n,m))_b, {\rm Ind}(\alpha\xi_{41}(n))$, where $n,m\in T_1$, $n<m$.

\emph{Type II.} Let $\chi$ be a character of $\F_q^\times$ with $\chi^2 \neq \triv_{\F_q^\times}$. Then the induced representation $\chi \times \chi \rtimes \sigma$ decomposes into two irreducible constituents
\begin{center}
IIa : $\chi\St_{\GL(2)} \rtimes \sigma$ \, \, and \, \, IIb: $\chi \triv_{\GL(2)} \rtimes \sigma$.
\end{center}
Let $\alpha=(\chi\sigma)^2$. Characters of Type IIa are the characters ${\rm Ind}(\alpha\chi_9(n))_a$ or ${\rm Ind}(\alpha\chi_9(n))_b$, where $n\in T_1$.  Characters of Type IIb are the characters\\ ${\rm Ind}(\alpha\chi_8(n))_a$ or ${\rm Ind}(\alpha\chi_8(n))_b$, where $n\in T_1$.

\emph{Type III.} Let $\chi$ be a character of $\F_q^\times$ such that $\chi \neq \triv_{\F_q^\times}$. Then $\chi \times \triv_{\F_q^\times} \rtimes \sigma$ decomposes into two irreducible constituents
\begin{center}
IIIa : $\chi \rtimes \sigma\St_{\rm GSp(2)}$ \, \, and \, \, IIIb: $\chi \rtimes \sigma \triv_{\rm GSp(2)}$.
\end{center}
Let $\alpha=\chi\sigma^2$. Characters of Type IIIa are the characters ${\rm Ind}(\alpha\xi_3'(n))_a,$\\ ${\rm Ind}(\alpha\xi_3'(n))_b,$ or ${\rm Ind}(\alpha\Phi_7)$, where $n\in T_1$. Characters of Type IIIb are the characters ${\rm Ind}(\alpha\xi_3(n))_a, {\rm Ind}(\alpha\xi_3(n))_b,$ or ${\rm Ind}(\alpha\Phi_5)$, where $n\in T_1$.

\emph{Type V.} $\xi \times \xi \rtimes \sigma$ decomposes into four irreducible constituents
\begin{center}
Va: $\sigma{\rm Ind}(\theta_1)$ \, \, Vb: $\sigma{\rm Ind}(\Phi_9)_a$ \\
Vc: $\sigma{\rm Ind}(\Phi_9)_b$ \, \, Vd: $\sigma{\rm Ind}(\theta_3).$
\end{center}

\emph{Type VI*.} $\triv_{\F_q^\times} \times \triv_{\F_q^\times} \rtimes \sigma$ decomposes into six irreducible constituents
\begin{center}
\begin{tabular}{ll}
VI*a: $\sigma \St_{\rm GSp(4)}$ & VI*b: $\sigma {\rm Ind}(\theta_9)_a$ \\
VI*c: $\sigma {\rm Ind}(\theta_9)_a$ & VI*d: $\sigma {\rm Ind}(\theta_{11})_a$ \\
VI*e: $\sigma {\rm Ind}(\theta_{12})_a$ & VI*f: $\sigma \triv_{\rm GSp(4)},$ \\
\end{tabular}
\end{center}
where $\St_{\rm GSp(4)}:={\rm Ind}(\theta_{13})_a$ is the Steinberg representation of $\GSp(4,q)$ and $\triv_{\rm GSp(4)}$ is the trivial representation of $\GSp(4,q)$.

\emph{Type VII.} These are the irreducible representations of the form $\chi \rtimes \pi$, where $\pi$ is an irreducible cuspidal representation of $\GL(2, q)$ and $\chi$ is a character of $\F_q^\times$. Representations of this form are irreducible if and only if $\chi \neq \triv_{\F_q^\times}$ and $\chi \neq \xi$, where $\xi$ is the non-trivial quadratic character such that $\xi\pi \cong \pi$. Let $\omega_\pi$ be the central character of $\pi$ and let $\alpha=\chi\omega_\pi$. Characters of Type VII are the characters ${\rm Ind}(\alpha\chi_5(n,m))$, ${\rm Ind}(\alpha\xi_{21}(n))$, and ${\rm Ind}(\alpha\xi_{41}'(m))$, where $n\in T_2, m\in T_1$.

\emph{Type VIII.} Let $\pi$ be an irreducible cuspidal representation of $\GL(2, q)$ with central character $\omega_\pi$. Then $\triv_{\F_q^\times} \rtimes \pi$ decomposes into two irreducible constituents
\begin{center}
VIIIa: ${\rm Ind}(\omega_\pi\Phi_3)$ \, \, and \, \, VIIIb: ${\rm Ind}(\omega_\pi\Phi_1).$
\end{center}
or
\begin{center}
VIIIa: ${\rm Ind}(\omega_\pi\xi_1'(n))$ \, \, and \, \, VIIIb: ${\rm Ind}(\omega_\pi\xi_1(n)).$
\end{center}
for some $n\in T_2$.

\emph{Type IX.} Let $\pi$ be an irreducible cuspidal representation of $\GL(2, q)$ with central character $\omega_\pi$ such that $\xi\pi \cong \pi$. Then $\xi \rtimes \pi$ decomposes into two irreducible constituents
\begin{center}
IXa: ${\rm Ind}(\xi\omega_\pi\theta_5)$ \, \, and \, \, IXb: ${\rm Ind}(\xi\omega_\pi\theta_7).$
\end{center}

\emph{Type X.} These are the irreducible representations of the form $\pi \rtimes \sigma$, where $\pi$ is an irreducible cuspidal representation of $\GL(2, q)$. Representations of this form are irreducible if and only if the central character $\omega_\pi$ of $\pi$ is not trivial. Let $\alpha=\omega_\pi\sigma^2$. Characters of Type X are the characters ${\rm Ind}(\alpha\chi_2(n))_a$ or ${\rm Ind}(\alpha\chi_2(n))_b$, where $n\in R_2$.

\emph{Type XI.} Let $\pi$ be an irreducible cuspidal representation of $\GL(2, q)$ with trivial central character. Then $\pi \rtimes \sigma$ decomposes into two irreducible constituents
\begin{center}
XIa: $\sigma{\rm Ind}(\chi_{7}(n))_{a}$ \, \, and \, \, XIb: $\sigma{\rm Ind}(\chi_{6}(n))_{a}.$
\end{center}
\end{thm}

Note that the irreducible constituents in Theorem \ref{cusprep} are sometimes written as twists of an irreducible character rather than in the form that is in Table \ref{GSpIrredCharacters}. 

The reader may have noticed that the type notation IV is not present. The representations in our type VI* turn out to be the irreducible representations that appear in the image of the non-supercuspidal types IV and VI under the functor $\mathcal{F}$ we will define in Section \ref{dimension}. Type VI* representations also occur as the images of non-supercuspidal representations of type V in the case that the quadratic character $\xi$ is unramified. See \cite{RS} and \cite{ST} for a description of the different types of non-supercuspidal representations of $\GSp(4,F)$, where $F$ is a non-archimedean local field of characteristic zero. 

\begin{cor}
The irreducible cuspidal characters of $\GSp(4,q)$ are the following characters in Table \ref{GSpIrredCharacters}: ${\rm Ind}(\alpha\chi_1(n))_a,$ ${\rm Ind}(\alpha\chi_1(n))_b,$ ${\rm Ind}(\alpha\chi_4(n,m))_a,$\\ ${\rm Ind}(\alpha\chi_4(n,m))_b,$ ${\rm Ind}(\alpha\xi_{21}'(n)),$ ${\rm Ind}(\alpha\theta_{10})_a,$ and ${\rm Ind}(\alpha\theta_{10})_b$.
\end{cor}

We remark that for our purposes of determining cuspidality of characters we only needed the portion of the character table of ${\rm Ind}(\alpha\chi)$ that is given in Table \ref{Indalphachivalues}. Since, in many cases, only the values on the classes with square multiplier are found with our method, we could not verify irreducible constituents, except for certain types, using the inner product. On the other hand, because we had the complete character list, the character values we could compute were enough to show that a particular irreducible character was a constituent. This can be done by listing the irreducible characters with the correct dimension and comparing enough character values.

The complete character tables of the irreducible non--cuspidal characters\\ $\sigma{\rm Ind}(\chi_6(n))_a$, $\sigma{\rm Ind}(\chi_7(n))_a$, ${\rm Ind}(\omega_\pi\Phi_1)$, ${\rm Ind}(\omega_\pi\Phi_3)$, $\sigma{\rm Ind}(\Phi_9)$, $\sigma{\rm Ind}(\theta_1)$,\\ $\sigma{\rm Ind}(\theta_3)$,  ${\rm Ind}(\xi\omega_\pi\theta_5)$, ${\rm Ind}(\xi\omega_\pi\theta_7)$, $\sigma{\rm Ind}(\theta_9)_a$, $\sigma{\rm Ind}(\theta_{11})_a$, $\sigma{\rm Ind}(\theta_{12})_a$, $\sigma\St_{\rm GSp(4)}$, ${\rm Ind}(\omega_\pi\xi_{1}(n))_a$, ${\rm Ind}(\omega_\pi\xi_{1}'(n))_a$, ${\rm Ind}(\chi\omega_\pi\xi_{21}(n))$, and ${\rm Ind}(\chi\omega_\pi\xi_{41}'(n))$ are given at the end of this paper. These tables were completed using both \cite{Sr} and \cite{Shinoda}. 

{\it Proof of Theorem \ref{cusprep}} The irreducible non-cuspidal representations are supported in the Borel, the Siegel parabolic, or the Klingen parabolic. We first consider those supported in the Borel.

\underline{Borel:} Let $\chi_1$, $\chi_2$, and $\sigma$ be characters of $\F_q^\times$. As in \ref{Borelreps}, these characters are used to define a representation of the Borel subgroup and induced to $\GSp(4,q)$ to obtain the representation $\chi_1\times\chi_2\rtimes\sigma$. From its character table, we have $$\chi_1\times\chi_2\rtimes\sigma\cong\chi_2\times\chi_1\rtimes\sigma.$$ We also have $$(\chi_1\times\chi_2\rtimes\sigma,\, \chi_\mathcal{G})=1,$$ indicating that exactly one irreducible constituent of $\chi_1\times\chi_2\rtimes\sigma$ is generic. We also compute that the possible values of
$$(\chi_1\times\chi_2\rtimes\sigma,\chi_1\times\chi_2\rtimes\sigma)$$
are $1,2,4,$ and $8$ according to the conditions in the following table.
\begin{center}
\begin{tabular}{|l|l|}
\hline
Value & Conditions\\
\hline
\hline
1 & $\chi_1\neq\triv_{\F_q^\times}, \chi_2\neq\triv_{\F_q^\times}, \chi_2\neq\chi_1^{\pm1}$\\
\hline
2 & $\chi_2=\chi_1^{\pm 1}$ and $\chi_1^2\neq\triv_{\F_q^\times}$\\
\hline
2 & $\chi_i=\triv_{\F_q^\times}$, $\chi_j\neq \triv_{\F_q^\times}$ for $i\neq j$\\
\hline
4 & $\chi_1=\chi_2=\xi$\\
\hline
8 & $\chi_1=\chi_2=\triv_{\F_q^\times}$\\
\hline
\end{tabular}
\end{center}
Using either the character inner product or by adding character values, the constituents in types II--VI* are verified.

\underline{Klingen:} Let $\chi$ be a character of $\F_q^\times$ and let $\pi$ be an irreducible cuspidal representation of $\GL(2,q)$ with central character $\omega_\pi$.  As in \ref{Klingenreps}, define a representation of the Klingen parabolic subgroup and induce to $\GSp(4,q)$ to obtain the representation $\chi\rtimes\pi$. Let $\rho_{\chi\rtimes\pi}$ be the character of $\chi\rtimes\pi$. Then $(\rho_{\chi\rtimes\pi},\, \chi_\mathcal{G})=1,$ indicating that exactly one irreducible constituent of $\chi\rtimes\pi$ is generic. We also compute that the possible values of $(\rho_{\chi\rtimes\pi},\rho_{\chi\rtimes\pi})$ are 1 and 2. By adding character values, if $\chi\rtimes\pi$ is reducible then either (1) $\chi=\triv_{\F_q^\times}$ and $\triv_{\F_q^\times}\rtimes\pi={\rm Ind}(\omega_\pi\Phi_3)+{\rm Ind}(\omega_\pi\Phi_1)$ or $\triv_{\F_q^\times}\rtimes\pi={\rm Ind}(\omega_\pi\xi_1'(n))+{\rm Ind}(\omega_\pi\xi_1(n))$ or (2) $\chi=\xi$ and $\pi$ is such that $\xi\pi=\pi$ and we have $\xi\rtimes\pi={\rm Ind}(\xi\omega_\pi\theta_5)+{\rm Ind}(\xi\omega_\pi\theta_7).$

\underline{Siegel:} Let $\sigma$ be a character of $\F_q^\times$ and let $\pi$ be an irreducible cuspidal representation of $\GL(2,q)$ with central character $\omega_\pi$. As in \ref{Siegelreps}, define a representation of the Siegel parabolic subgroup and induce to $\GSp(4,q)$ to obtain the representation $\pi\rtimes\sigma$. We have $(\pi\rtimes\sigma,\, \chi_\mathcal{G})=1$ for all such representations $\pi\rtimes\sigma$, indicating that exactly one irreducible constituent of $\pi\rtimes\sigma$ is generic. Let $\chi_{\pi\rtimes\sigma}$ be the character of $\pi\rtimes\sigma$. Then the possible values of $(\chi_{\pi\rtimes\sigma},\chi_{\pi\rtimes\sigma})$ are $1$ and $2$. We have $(\chi_{\pi\rtimes\sigma},\chi_{\pi\rtimes\sigma})=2$ precisely when $\omega_\pi$ is trivial. By adding character values, if $\pi\rtimes\sigma$ is reducible then $$\pi\rtimes\sigma=\sigma{\rm Ind}(\chi_7(n))_a+\sigma{\rm Ind}(\chi_6(n))_a$$ for some $n\in T_2$, depending on $\pi$ and on our embedding of $\F_{q^4}^\times$ into $\C^\times$.
\qed

\subsubsection{Decompositions for types V and VI*}\label{decomp}
We give more information on the decompositions of the non-cuspidal representations supported in the Borel subgroup for types V and VI*.  These decompositions can be verified with the character tables provided in sections \ref{inducedchars} and \ref{noncuspchartables}.

\underline{Type V}: Constituents of $\xi\times\xi\rtimes\sigma$.
\begin{align*}
 \xi\times\xi\rtimes\sigma
  &=\xi\,\St_{\GL(2)}\rtimes\sigma
  +\xi\,\triv_{\GL(2)}
  \rtimes\sigma\\
 &=\xi\,\St_{\GL(2)}
  \rtimes\xi\sigma
  +\xi\,\triv_{\GL(2)}
  \rtimes\xi\sigma.
\end{align*}
Each of the four representations on the right side is reducible and has two constituents which are shown in the following table.
\begin{center}
\begin{small}
\begin{longtable}{|c||c|c|}
\caption[Type V constituents]{\small{Type V constituents}} \label{Vconstituents} \\
\hline
   {\parbox{2cm}{\strut\centering\small{}\strut}} &
   {\parbox{2.5cm}{\strut\centering\small{$\xi\,\St_{\GL(2)}\rtimes\xi\sigma$}\strut}} &
   {\parbox{2.5cm}{\strut\centering\small{$\xi\,\triv_{\GL(2)}\rtimes\xi\sigma$}\strut}} \\
\hline
\hline
\endfirsthead
\multicolumn{3}{c}{{\tablename} \thetable{} -- Continued} \\[0.5ex]
\hline
   {\parbox{2cm}{\strut\centering\small{}\strut}} &
   {\parbox{2.5cm}{\strut\centering\small{$\xi\,\St_{\GL(2)}\rtimes\xi\sigma$}\strut}} &
   {\parbox{2.5cm}{\strut\centering\small{$\xi\,\triv_{\GL(2)}\rtimes\xi\sigma$}\strut}} \\
\hline
\hline
\endhead
\hline
\endlastfoot
$\xi\,\St_{\GL(2)}\rtimes\sigma$ & $\sigma{\rm Ind}(\theta_1)$ & $\sigma{\rm Ind}(\Phi_9)_a$ \\
\hline
$\xi\,\triv_{\GL(2)}\rtimes\sigma$ & $\sigma{\rm Ind}(\Phi_9)_b$ & $\sigma{\rm Ind}(\theta_3)$ \\
\end{longtable}
\end{small}
\end{center}

\underline{Type VI*}: Constituents of $\triv_{\F_q^\times}\times\triv_{\F_q^\times}\rtimes\sigma$.
\begin{align*}
\triv_{\F_q^\times}\times\triv_{\F_q^\times}\rtimes\sigma&=\St_{\GL(2)}\rtimes\sigma + \triv_{\GL(2)}\rtimes\sigma\\
 &=\triv_{\F_q^\times}\rtimes\sigma\St_{\rm GSp(2)} + \triv_{\F_q^\times}\rtimes\sigma \triv_{\rm GSp(2)}.
\end{align*}
Each of the four representations given on the right is reducible and has three irreducible constituents, which are shown in the following table. The common factor $\sigma{\rm Ind}(\theta_9)_a$ occurs as a constituent of each of the four representations $\St_{\GL(2)}\rtimes\sigma,\, \triv_{\GL(2)}\rtimes\sigma,$ $\triv_{\F_q^\times}\rtimes\sigma\St_{\rm GSp(2)},$ and $\triv_{\F_q^\times}\rtimes\sigma \triv_{\rm GSp(2)}$.

\begin{center}
\begin{small}
\begin{longtable}{|c||c|c|c|}
\caption[Type VI* constituents]{\small{Type VI* constituents}} \label{VIconstituents} \\

\hline
   {\parbox{2cm}{\strut\centering\small{}\strut}} &
   {\parbox{2cm}{\strut\centering\small{$\St_{\GL(2)}\rtimes\sigma$}\strut}} &
   {\parbox{2cm}{\strut\centering\small{(common)}\strut}} &
   {\parbox{2cm}{\strut\centering\small{$\triv_{\GL(2)}\rtimes\sigma$}\strut}} \\
\hline
\hline
\endfirsthead

\multicolumn{4}{c}{{\tablename} \thetable{} -- Continued} \\[0.5ex]
\hline
   {\parbox{2cm}{\strut\centering\small{}\strut}} &
   {\parbox{2cm}{\strut\centering\small{$\St_{\GL(2)}\rtimes\sigma$}\strut}} &
   {\parbox{2cm}{\strut\centering\small{(common)}\strut}} &
   {\parbox{2cm}{\strut\centering\small{$\triv_{\GL(2)}\rtimes\sigma$}\strut}} \\
\hline
\hline
\endhead

 \hline
\endlastfoot

$\triv_{\F_q^\times}\rtimes\sigma\St_{\rm GSp(2)}$ & $\sigma\St_{\rm GSp(4)}$ & $\sigma{\rm Ind}(\theta_9)_a$ & $\sigma{\rm Ind}(\theta_{11})_a$ \\
\hline
(common factor) & $\sigma{\rm Ind}(\theta_9)_a$ & $-$ & $\sigma{\rm Ind}(\theta_9)_a$ \\
\hline
$\triv_{\F_q^\times}\rtimes\sigma \triv_{\rm GSp(2)}$ & $\sigma{\rm Ind}(\theta_{12})_a$ & $\sigma{\rm Ind}(\theta_9)_a$ & $\sigma \triv_{\rm GSp(4)}$ \\
\end{longtable}
\end{small}
\end{center}

\section{Dimension formulas}
\label{dimension}

A representation $(\pi, V )$ of a group $G$ defined over a non-archimedean local field is called {\em smooth} if every vector in $V$ is fixed by an open-compact subgroup $K$ of $G$. A representation $(\pi, V )$ is called {\em admissible} if the space $V^K$ is finite dimensional for any open compact subgroup $K$.

A fundamental problem for representations of groups defined over local fields is that they are not always semi-simple, i.e., they do not decompose into a direct sum of irreducible constituents. For such reducible representations, one must determine how they can be constructed from irreducible representations as extensions of a quotient by a subrepresentation.

Consider the group $G=\GSp(4, F)$, where $F$ is a non-archimedean local field of characteristic zero with ring of integers $\mathfrak{o}$ and maximal ideal $\mathfrak{p}$ such that $\mathfrak{o}/\mathfrak{p}$ is a finite field with $q$ elements. Fix a generator $\varpi$ of $\mathfrak{p}$. If $x$ is in $F^\times$, then define $v(x)$ to be the unique integer such that $x = u\varpi^{v(x)}$ for some unit $u$ in $\mathfrak{o}^\times$. Write $\nu(x)$ or $|x|$ for the normalized absolute value of $x$ such that $\nu(\varpi) = q^{-1}$.

The {\em principal congruence subgroup of level $\mathfrak{p}^n$} of $G$, denoted by $\Gamma(\mathfrak{p}^n)$, is 
\begin{displaymath}
\Gamma(\mathfrak{p}^n) = \{ g\in G : g\equiv I_4\, ({\rm mod}\, \mathfrak{p}^n) \}.
\end{displaymath}
We have the following short exact sequence
\begin{displaymath}
1 \longrightarrow \Gamma(\mathfrak{p}) \longrightarrow K \longrightarrow \GSp(4,\mathfrak{o}/\mathfrak{p}) \longrightarrow 1,
\end{displaymath}
where $K = \GSp(4,\mathfrak{o})$ is a maximal compact subgroup of $G$. 

Define the functor $\mathcal{F}$ which maps admissible representations of $\GSp(4, F)$ to representations of $\GSp(4,q)$ as follows. Let $(\pi, V)$ be an admissible representation of $\GSp(4, F)$. Consider the action of $K$ on the subspace $V^{\Gamma(\mathfrak{p})}$ obtained via the usual restriction functor. Since $\Gamma(\mathfrak{p})$ is a normal subgroup of $K$ and $\Gamma(\mathfrak{p})$ acts trivially, there is an action of $K/\Gamma(\mathfrak{p})$ on this subspace. Indeed, $K/\Gamma(\mathfrak{p})$ acts on vectors $v\in V^{\Gamma(\mathfrak{p})}$ by the following action:
$$k\Gamma(\mathfrak{p})\cdot v = \pi(k)v,$$
where $k\in K$. This action is well-defined since if $k_1\Gamma(\mathfrak{p})=k_2\Gamma(\mathfrak{p})$, where $k_1,k_2\in K$, then $k_1=k_2g$ for some $g\in\Gamma(\mathfrak{p})$. It follows that
$$\pi(k_1)v=k_1\Gamma(\mathfrak{p})\cdot v=k_2g\Gamma(\mathfrak{p})\cdot v=\pi(k_2g)v=\pi(k_2)\pi(g)v=\pi(k_2).$$

By our short exact sequence, we have $K/\Gamma(\mathfrak{p})\cong \GSp(4,\mathfrak{o}/\mathfrak{p})\cong \GSp(4,q)$. So we define an action of $\GSp(4,q)$ on $V^{\Gamma(\mathfrak{p})}$ as the action of $K/\Gamma(\mathfrak{p})$ on $V^{\Gamma(\mathfrak{p})}$, yielding a representation $\tilde{\pi}$ of $\GSp(4,q)$ on the space $V^{\Gamma(\mathfrak{p})}$. Hence, we define $\mathcal{F}((\pi,V))=\left(\tilde{\pi},V^{\Gamma(\mathfrak{p})}\right)$. The functor $\mathcal{F}$ is a functor of invariants. So it is left exact. It is right exact because $\Gamma(\mathfrak{p})$ is compact, \cite{Bernstein}. Similarly, one defines a functor which maps admissible representations $\pi$ of a split connected reductive subgroup $G$ of $\GL(n,F)$ to representations $\tilde{\pi}$ of $G(\F_q)$ and a functor which maps smooth characters $\chi$ of $F$ to characters $\tilde{\chi}$ of $\F_q$. These functors will also be denoted by $\mathcal{F}.$

For certain admissible non-supercuspidal representations $(\pi,V)$ of $G$, the dimension of the subspace $V^{\Gamma(\mathfrak{p})}$ can be computed by determining the representation of $\GSp(4,q)$ obtained by applying the restriction functor $\mathcal{F}$ to $\pi$. To this end, the following lemmas are useful. We use the notations $K(1) = \Gamma(\mathfrak{p}),$ $K_M(1) = \Gamma(\mathfrak{p}) \cap M,$ $K_{\GL(2)}=\GL(2,\mathfrak{o}),$ and $$K_{\GL(2)}(1)=\Gamma_{\GL(2)}(\mathfrak{p})=\{g\in K_{\GL(2)}: g\equiv I_2 ({\rm mod}\, \mathfrak{p})\}.$$ 

\begin{lemma}\label{KM1fixedvectors}
Let $G$ be a split connected reductive subgroup of $\GL(n,F)$ with center $Z$. Let $K_G=G(\mathfrak{o}), K_G(1)=\{g\in K_G:g\equiv I\, {\rm mod}\, \mathfrak{p}\}.$ Let $P$ be a parabolic subgroup of $G$. Write $P=MN$, where $M$ is the Levi subgroup and $N$ is the unipotent radical of $P$. Let $K_M(1) = K_G(1) \cap M.$ Let $(\rho,V)$ be a subquotient of a representation which is parabolically induced from a representation $\sigma$ of $P$. If $\sigma$ has a non-zero subspace of $K_M(1)$-fixed vectors, then $\rho$ has a non-zero subspace of $K(1)$-fixed vectors.
\end{lemma}
\begin{proof}
Since $K(1)$ is an open compact subgroup of $G$ which has the Iwahori decomposition with respect to $H$, the Jacquet functor
$$J_N: V\rightarrow J_N(V)$$
gives, by Theorem 2.2 of \cite{MoyPrasad}, a surjection
$$V^{K(1)}\rightarrow\left(J_N(V)\right)^{K_M(1)}.$$
Since the space of $K_M(1)$-fixed vectors in the Jacquet module is non-zero, $\rho$ has a non-zero subspace of $K(1)$-fixed vectors.
\qed
\end{proof}

\begin{lemma}\label{Borelfixedvectors}
Let $\chi_1,\chi_2,\sigma$ be smooth characters of $F^\times$ that are at most tamely ramified. Then we have
\begin{equation}
\mathcal{F}(\chi_1\times\chi_2\rtimes\sigma) \cong \tilde{\chi}_1\times \tilde{\chi}_2\rtimes \tilde{\sigma}.
\end{equation}
where $\tilde{\chi}_i:=\mathcal{F}(\chi_i), \tilde{\sigma}:=\mathcal{F}(\sigma)$.
\end{lemma}
\begin{proof}
The standard model for this induced representation is the space $V$ of smooth functions $f: G\rightarrow \C$ that satisfy
\begin{displaymath}
f(\left(\begin{smallmatrix}
a & * & * & * \\
& b & * & * \\
& & cb^{-1} & * \\
& & & ca^{-1} \\
\end{smallmatrix}\right)g) = |a^2 b||c|^{-3/2}\chi_1(a)\chi_2(b)\sigma(c)f(g)
\end{displaymath}
for all $g\in G$, with group action by right translation. Restricting the characters $\chi_1, \chi_2,$ and $\sigma$ to $\mathfrak{o}^\times$ yields characters on $\F_q^\times$, i.e., 
$$\tilde{\chi}_i=\chi_i|_{\mathfrak{o}^\times}, \tilde{\sigma}=\sigma|_{\mathfrak{o}^\times}:\F_q^\times\cong (\mathfrak{o}/\mathfrak{p})^\times=\mathfrak{o}^\times/(1+\mathfrak{p})\rightarrow\C^\times.$$

Let $g\in G$. By the Iwasawa decomposition of $G$, we can write $g=bk$, where $k\in K=\GSp(4,\mathfrak{o})$ and $b\in B(F)$, where $B(F)$ is the standard Borel subgroup of $G$. Note that the Iwasawa decomposition is not unique and $K(1)\subset\GSp(4,\mathfrak{o})$. Consider the subspace $V^{K(1)}$ of $K(1)$-fixed vectors in the representation $\chi_1\times\chi_2\rtimes\sigma$,
\begin{displaymath}
V^{K(1)} = \{ f: f(g\gamma) = f(g)\, {\rm for\, all}\, \gamma\in K(1), g\in G\}.
\end{displaymath}
Let $f\in V^{K(1)}$. Then for all $\gamma\in K(1)$,
\begin{equation}\label{fbkprop}
f(bk)=f(g)=f(g\gamma)=f(bk\gamma).
\end{equation}
This space is then isomorphic to the space $\tilde{V}^{K(1)}$ of functions $f: K\rightarrow \C$ satisfying Equation \ref{fbkprop} via the restriction map: $f\mapsto\overline{f},$ where $\overline{f}(k)=f(k)$ for $k\in K.$ This map is well-defined. These functions can then be considered functions on $$K/K(1) \cong \GSp(4,\mathfrak{o}/\mathfrak{p}) \cong \GSp(4,q).$$ That is, $f$ can be considered a function on the finite group $\GSp(4,q)$ satisfying
\begin{displaymath}
f(\left(\begin{smallmatrix}
a & * & * & * \\
& b & * & * \\
& & cb^{-1} & * \\
& & & ca^{-1} \\
\end{smallmatrix}\right)g) = \tilde{\chi}_1(a)\tilde{\chi}_2(b)\tilde{\sigma}(c)f(g).
\end{displaymath}
So
\begin{displaymath}
\mathcal{F}(\chi_1\times\chi_2\rtimes\sigma) \cong  \tilde{\chi}_1\times \tilde{\chi}_2\rtimes \tilde{\sigma}.
\end{displaymath}
\qed
\end{proof}

\begin{cor}
Let $\chi_1,\chi_2,\sigma$ be smooth characters of $F^\times$ that are at most tamely ramified. Then the dimension of $(\chi_1\times\chi_2\rtimes\sigma)^{K(1)}$ is $(q^2+1)(q+1)^2$. 
\end{cor}

\begin{lemma}\label{supercuspreps}
Let $G$ be a split connected reductive subgroup of $\GL(n,F)$ with center $Z$. Let $K_G=G(\mathfrak{o}), K_G(1)=\{g\in K_G:g\equiv I\, {\rm mod}\, \mathfrak{p}\}.$ Suppose $(\pi,V)$ is an irreducible smooth supercuspidal representation of $G$ which admits a non-zero space of $K_{G}(1)$-fixed vectors. Then $ \mathcal{F}(\pi)\cong \tilde{\pi}$, where $\tilde{\pi}$ is an irreducible cuspidal representation of $G(\F_q)$ obtained from the restriction of $\pi$ to $K_{G}$.
\end{lemma}
\begin{proof}
Since $\pi$ is an irreducible smooth supercuspidal representation of $G$ which admits a non-zero space of $K_{G}(1)$-fixed vectors, $\pi$ has depth zero. So $\pi|_{K_{G}}$ contains an irreducible representation, say $\tilde{\pi}$, of $K_{G}/{K_{G}(1)}\cong G(\mathfrak{o}/\mathfrak{p})\cong G(\F_q)$. 

We claim that $\tilde{\pi}$ is cuspidal. Indeed, if not, then $\tilde{\pi}$ contains a cuspidal representation of a proper parabolic subgroup $\tilde{P}$ of $G(\mathfrak{o}/\mathfrak{p})$. Write $\tilde{P}=\tilde{M}\tilde{N}$, where $\tilde{M}$ is the Levi subgroup and $\tilde{N}$ is the unipotent radical of $P$. The inverse image of $\tilde{P}$ under the map $K\to G(\mathfrak{o}/\mathfrak{p})$ is a non-maximal parahoric subgroup $P\subset K_G$. The inverse image of $\tilde{U}$ is $P^+$. Then $\pi^{P^+}\neq 0$. By Proposition 6.7 of \cite{MoyPrasad}, $\pi$ then has a non-trivial Jacquet module. This contradicts our assumption that $\pi$ is supercuspidal.

The representation $\pi$ also contains an extension $\tilde{\pi}'$ of $\tilde{\pi}$ to $ZK_{G}$, where $Z$ is the center of $G$. We have that ${\rm Hom}_{ZK_{G}}(\tilde{\pi}',\pi|_{ZK_{G}})\neq 0$. By Frobenius reciprocity,
$${\rm Hom}_G({\rm Ind}_{ZK_{G}}^G \tilde{\pi}',\pi)\cong {\rm Hom}_{ZK_{G}}(\tilde{\pi}',\pi|_{ZK_{G}}).$$
Since $\tilde{\pi}'$ is irreducible and $g\in G$ intertwines $\tilde{\pi}'$ if and only if $g\in ZK$, we have that ${\rm Ind}_{ZK_{G}}^G\tilde{\pi}'$ is irreducible by Theorem 11.4 in \cite{BH}. So it must be isomorphic to $\pi$.

Consider $\left({\rm Ind}_{ZK_{G}}^G\tilde{\pi}'\right)|_{{K_{G}(1)}}$. Its decomposition using Mackey's restriction formula is 
$$\left({\rm Ind}_{ZK_{G}}^G\tilde{\pi}'\right)|_{{K_{G}(1)}}\cong\bigoplus_{x\in ZK_{G}\backslash G/{K_{G}(1)}} {\rm Ind}_{ZK^x(1)}^{{K_{G}(1)}} ((\tilde{\pi}')^x|_{ZK^x(1)}).$$
where $ZK^x(1) = ZK_{G}^x\cap{K_{G}(1)},$ $ZK_{G}^x=x^{-1}ZK_{G} x$, and $(\tilde{\pi}')^x$ is the representation of $ZK_{G}^x$ defined by $$(\tilde{\pi}')^x(x^{-1}hx) = \tilde{\pi}'(h)$$ for $h\in ZK_{G}, x\in ZK_{G}\backslash G/{K_{G}(1)}$. 

If ${\rm Ind}_{ZK^x(1)}^{{K_{G}(1)}} ((\tilde{\pi}')^x|_{ZK^x(1)})$ contains the trivial representation $\triv_{K_{G}(1)}$ of ${K_{G}(1)}$, then ${\rm Ind}_{ZK_{G}}^G\tilde{\pi}'|_{K_{G}}$ contains an irreducible representation, say $\tau$, such that $$\tau|_{{K_{G}(1)}}\supset\textbf{1}_{{K_{G}(1)}}.$$
This implies that $\tau$ is trivial on ${K_{G}(1)}$. The representations $\tilde{\pi},\tau\subset\pi$ then intertwine. That is, there exists $x\in G$ such that
$${\rm Hom}_{K_{G}^x\cap K_{G}}(\tilde{\pi}^x,\tau)\neq 0.$$
This $x$ is then necessarily an element of $ZK_{G}$. 

So, by the cuspidality of $\tilde{\pi}$ and using representatives for the Cartan decomposition $K_{G}\backslash G/K_{G}$, we have $\tilde{\pi}\cong\tau$. 
It follows that $\mathcal{F}(\pi)\cong\tilde{\pi}'|_{{K_{G}(1)}}\cong\tilde{\pi}$.
\qed
\end{proof}

\begin{lemma}\label{PQfixedvectors}
Let $G=\GSp(4, F)$ and let $H$ be either the Klingen or the Siegel parabolic subgroup of $G$. Write $H=MN$, where $M$ is its Levi subgroup and $N$ its unipotent radical. Let $\pi$ be an irreducible smooth supercuspidal representation of $M$ which admits a non-zero space of $K_M(1)$-fixed vectors. Then
\begin{equation}\label{functorofinductioneqn}
\mathcal{F}({\rm Ind}_H^G(\pi)) \cong {\rm Ind}_{H(\F_q)}^{G(\F_q)}(\mathcal{F}(\pi)).
\end{equation}
\end{lemma}
\begin{proof}
Let $V$ be the standard model of ${\rm Ind}_H^G(\pi)$. By Lemma \ref{KM1fixedvectors}, $V^{K(1)}\neq\{0\}$. Restricting $\pi$ to $K_{\GL(2)}$, we have that $\mathcal{F}(\pi)$ is isomorphic to a cuspidal representation, say $\tilde{\pi}$, of $K_{\GL(2)}/{K_{\GL(2)}(1)}\cong\GL(2,q)$ obtained from the restriction of $\pi$ to $K_{\GL(2)}$ by Lemma \ref{supercuspreps}. Let $f:G\rightarrow\C$ be a non-zero vector in $V^{K(1)}$. By the Iwasawa decomposition of $G$, every element $g\in G$ can be written as $g=hk$, for some $k\in K$ and $h\in H$. Note that the Iwasawa decomposition is not unique and $K(1)\subset\GSp(4,\mathfrak{o})$. Write $g=hk$, where $k\in K$ and $h\in H$. Then for all $\gamma\in K(1)$,
\begin{equation}\label{PQfhkprop}
f(hk)=f(g)=f(g\gamma)=f(hk\gamma).
\end{equation}
The subspace of $K(1)$-fixed vectors of ${\rm Ind}_H^G(\pi)$ is then isomorphic to the space $\tilde{V}^{K(1)}$ of functions $f: K\rightarrow \C$ satisfying Equation \ref{PQfhkprop} via the restriction map: $f\mapsto\overline{f},$ where $\overline{f}(k)=f(k)$ for $k\in K.$ This map is well-defined. These functions can then be considered functions on $$K/K(1) \cong \GSp(4,\mathfrak{o}/\mathfrak{p}) \cong \GSp(4,q).$$ That is, $f$ can be considered a function on the finite group $\GSp(4,q)$ satisfying
\begin{displaymath}
f(hg) = \tilde{\pi}(h)f(g)
\end{displaymath}
for all $h\in H(q)$, $g\in\GSp(4,q)$, where $H(q)$ is the Klingen parabolic subgroup, respectively the Siegel parabolic subgroup, of $\GSp(4,q)$. So
$$\mathcal{F}({\rm Ind}_H^G(\pi)) \cong {\rm Ind}_{H(\F_q)}^{G(\F_q)}(\mathcal{F}(\pi)).$$
\qed
\end{proof}

\begin{cor}
Let $G=\GSp(4, F)$ and let $H$ be either the Klingen or the Siegel parabolic subgroup of $G$. Write $H=MN$, where $M$ is its Levi subgroup and $N$ its unipotent radical. Let $\pi$ be an irreducible smooth supercuspidal representation of $M$ which admits a non-zero space of $K_M(1)$-fixed vectors. Then the dimension of $({\rm Ind}_H^G(\pi))^{K(1)}$ is $q^4-1$. 
\end{cor}

We are now ready to state and prove our last main result.

\begin{thm}\label{nonsupercuspidaldimensions}
Let $\chi,\chi_1,\chi_2,\xi,\xi_1,\xi_2,\sigma$ be characters of $F^\times$ that are trivial on\\ $1+\mathfrak{p}$. Suppose $\xi,\xi_1,\xi_2$ are non-trivial quadratic characters with $\xi_1$ tamely ramified and $\xi_2$ unramified. Let $\pi$ be an irreducible smooth supercuspidal representation of $\GL(2,F)$ which admits a non-zero space of $K_{\GL(2)}(1)$-fixed vectors. Then Table \ref{Fixed} gives the dimensions of $K(1)$-fixed vectors in the ``$\mathcal{F}(\rho)$ dimension'' column for the non-supercuspidal representations $\rho$ of $G=\GSp(4, F)$ supported in the Borel, the Siegel parabolic, or the Klingen parabolic subgroup. The ``$\mathcal{F}(\rho,V)$" column gives the image of the representation $\rho$ under $\mathcal{F}$.
\end{thm}

\begin{center}
\begin{small}
\begin{longtable}{|c|c|c|l|c|}
\caption[Dimensions of $K(1)$-fixed vectors]{\small{Dimensions of $K(1)$-fixed vectors}} \label{Fixed} \\
\hline
   \multicolumn{1}{|c|}{\parbox{0cm}{\strut\centering\small\textbf{}\strut}} &
   \multicolumn{1}{|c|}{\parbox{0cm}{\strut\centering\small\textbf{}\strut}} &
   \multicolumn{1}{|c|}{\parbox{2.8cm}{\strut\centering\small\textbf{Representation $\rho$}\strut}} &
   \multicolumn{1}{|c|}{\parbox{2cm}{\strut\centering\small\textbf{$\mathcal{F}(\rho)$ dimension}\strut}} &
   \multicolumn{1}{|c|}{\parbox{1.8cm}{\strut\centering\small\textbf{$\mathcal{F}(\rho,V)$}\strut}} \\
\hline
\hline
\endfirsthead
\multicolumn{5}{c}{{\tablename} \thetable{} -- Continued} \\[0.5ex]
\hline
   \multicolumn{1}{|c|}{\parbox{0cm}{\strut\centering\small\textbf{}\strut}} &
   \multicolumn{1}{|c|}{\parbox{0cm}{\strut\centering\small\textbf{}\strut}} &
   \multicolumn{1}{|c|}{\parbox{2.8cm}{\strut\centering\small\textbf{Representation $\rho$}\strut}} &
   \multicolumn{1}{|c|}{\parbox{2cm}{\strut\centering\small\textbf{$\mathcal{F}(\rho)$ dimension}\strut}} &
   \multicolumn{1}{|c|}{\parbox{1.8cm}{\strut\centering\small\textbf{$\mathcal{F}(\rho,V)$}\strut}} \\
\hline
\hline
\endhead
 \hline
\endlastfoot

\makebox{\small I} & & $\chi_1 \times \chi_2 \rtimes \sigma$\quad(irreducible) & $(q^2+1)(q+1)^2$ & $\tilde{\chi}_1 \times \tilde{\chi}_2 \rtimes \tilde{\sigma}$ \\
\hline
& \small\mbox{a} & $\chi\St_{\GL(2,F)} \rtimes \sigma$ & $q(q^2+1)(q+1)$ & $\tilde{\chi}\St_{\GL(2,q)} \rtimes \tilde{\sigma}$ \\
\cline{3-5} \raisebox{2.5ex}[-2.5ex]{\small II}
& \small\mbox{\small b} & $\chi \triv_{\GL(2,F)} \rtimes \sigma$ & $(q^2+1)(q+1)$ & $\tilde{\chi} \triv_{\GL(2,q)} \rtimes \tilde{\sigma}$ \\
\hline
& \small\mbox{a} & $\chi \rtimes \sigma\St_{\GL(2)}$ & $q(q^2+1)(q+1)$ & $\tilde{\chi} \rtimes \tilde{\sigma}\St_{\GSp(2,q)}$ \\
\cline{3-5} \raisebox{2.5ex}[-2.5ex]{\small III}
& \small\mbox{b} & $\chi \rtimes \sigma \triv_{\GL(2)}$ & $(q^2+1)(q+1)$ &  $\tilde{\chi} \rtimes \tilde{\sigma} \triv_{\GSp(2,q)}$\\
\hline
& \small\mbox{a} & $\sigma\St_{\GSp(4,F)}$ & $q^4$ & $\tilde{\sigma}\St_{\GSp(4,q)}$ \\
\cline{3-5}
&\small\mbox{b} & $L(\nu^2,\nu^{-1}\sigma\St_{\GL(2)})$ & $q(q^2+q+1)$ & \begin{minipage}{15ex}
\begin{center}
$\tilde{\sigma}{\rm Ind}(\theta_9)_a+\tilde{\sigma}{\rm Ind}(\theta_{11})_a$
\end{center}
\end{minipage} \\
\cline{3-5} \raisebox{2.5ex}[-2.5ex]{\small IV} &
\small\mbox{c} & $L(\nu^{3/2}\St_{\GL(2)},\nu^{-3/2}\sigma)$ & $q(q^2+q+1)$ & \begin{minipage}{15ex}
\begin{center}
$\tilde{\sigma}{\rm Ind}(\theta_9)_a+\tilde{\sigma}{\rm Ind}(\theta_{12})_a$ 
\end{center}
\end{minipage} \\
\cline{3-5}
& \small\mbox{d} & $\sigma \triv_{\GSp(4)}$ & $1$ & $\tilde{\sigma} \triv_{\GSp(4,q)}$\\
\hline
& \small\mbox{a} & $\delta([\xi_1,\nu\xi_1],\nu^{-1/2}\sigma)$ & $q^2(q^2+1)$ & $\tilde{\sigma}{\rm Ind}(\theta_1)$ \\
\cline{3-5}
&\small\mbox{b} & $L(\nu^{1/2}\xi_1\St_{\GL(2,F)},\nu^{-1/2}\sigma)$ & $q(q^2+1)$ & $\tilde{\sigma}{\rm Ind}(\Phi_9)_a$ \\
\cline{3-5} \raisebox{2.5ex}[-2.5ex]{\small V} &
\small\mbox{c} & $L(\nu^{1/2}\xi_1\St_{\GL(2,F)},\xi_1\nu^{-1/2}\sigma)$ & $q(q^2+1)$ & $\tilde{\sigma}{\rm Ind}(\Phi_9)_b$ \\
\cline{3-5}
& \small\mbox{d} & $L(\nu\xi_1,\xi_1\rtimes\nu^{-1/2}\sigma)$ & $q^2+1$ & $\tilde{\sigma}{\rm Ind}(\theta_3)$ \\
\hline
& \small\mbox{a} & $\delta([\xi_2,\nu\xi_2],\nu^{-1/2}\sigma)$ & $q^4+\frac{1}{2}q(q^2+1)$ & \begin{minipage}{15ex}
\begin{center}
$\tilde{\sigma}\St_{\GSp(4,q)}+\tilde{\sigma}{\rm Ind}(\theta_{12})_a$
\end{center}
\end{minipage}\\
\cline{3-5}
&\small\mbox{b} & $L(\nu^{1/2}\xi_2\St_{\GL(2,F)},\nu^{-1/2}\sigma)$ & $\frac{1}{2}q(q+1)^2$ & $\tilde{\sigma}{\rm Ind}(\theta_9)_a$ \\
\cline{3-5} \raisebox{2.5ex}[-2.5ex]{\small V*} &
\small\mbox{c} & $L(\nu^{1/2}\xi_2\St_{\GL(2,F)},\xi_2\nu^{-1/2}\sigma)$ & $\frac{1}{2}q(q+1)^2$ & $\tilde{\sigma}{\rm Ind}(\theta_9)_a$ \\
\cline{3-5}
& \small\mbox{d} & $L(\nu\xi_2,\xi_2\rtimes\nu^{-1/2}\sigma)$ & $1+\frac{1}{2}q(q^2+1)$ & \begin{minipage}{15ex}
\begin{center}
$\tilde{\sigma} \triv_{\GSp(4,q)}+\tilde{\sigma}{\rm Ind}(\theta_{11})_a$
\end{center}
\end{minipage} \\
\hline
& \small\mbox{a} & $\tau(S,\nu^{-1/2}\sigma)$ & $q^4+\frac{1}{2}q(q+1)^2$ & \begin{minipage}{15ex}
\begin{center}
$\tilde{\sigma}\St_{\GSp(4,q)}+\tilde{\sigma}{\rm Ind}(\theta_9)_a$
\end{center}
\end{minipage}\\
\cline{3-5}
&\small\mbox{b} & $\tau(T,\nu^{-1/2}\sigma)$ & $\frac{1}{2}q(q^2+1)$ & $\tilde{\sigma}{\rm Ind}(\theta_{11})_a$ \\
\cline{3-5} \raisebox{2.5ex}[-2.5ex]{\small VI} &
\small\mbox{c} & $L(\nu^{1/2}\St_{\GL(2,F)},\nu^{-1/2}\sigma)$ & $\frac{1}{2}q(q^2+1)$ & $\tilde{\sigma}{\rm Ind}(\theta_{12})_a$ \\
\cline{3-5}
& \small\mbox{d} & $L(\nu,\triv_{F^\times}\rtimes\nu^{-1/2}\sigma)$ & $1+\frac{1}{2}q(q+1)^2$ & \begin{minipage}{15ex}
\begin{center}
$\tilde{\sigma} \triv_{\GSp(4,q)}+\tilde{\sigma}{\rm Ind}(\theta_9)_a$
\end{center}
\end{minipage} \\
\hline

\makebox{\small VII} & & $\chi \rtimes \pi$\quad (irreducible) & $q^4 - 1$ & $\tilde{\chi} \rtimes \tilde{\pi}$ \\
\hline
\makebox{\small VIII} & & $\tau(S,\pi)+\tau(T,\pi)$ & $q^4-1$ & \begin{minipage}{15ex}
\begin{center}
$({\rm Ind}(\omega_{\tilde{\pi}}\Phi_1)+{\rm Ind}(\omega_{\tilde{\pi}}\Phi_3))$

or

$({\rm Ind}(\omega_{\tilde{\pi}}\xi_1'(n))+{\rm Ind}(\omega_{\tilde{\pi}}\xi_1(n)))$
\end{center}
\end{minipage}\\
\hline
\makebox{\small IX} & & \begin{minipage}{15ex}
\begin{center}
$\delta(\nu\xi_1,\nu^{-1/2}\pi)+L(\nu\xi_1,\nu^{-1/2}\pi)$
\end{center}
\end{minipage}
 & $q^4-1$ & \begin{minipage}{15ex}
\begin{center}
${\rm Ind}(\tilde{\xi}\omega_{\tilde{\pi}}\theta_5)+{\rm Ind}(\tilde{\xi}\omega_{\tilde{\pi}}\theta_7)$
\end{center}
\end{minipage}\\
\hline
\makebox{\small IX*} & & \begin{minipage}{15ex}
\begin{center}
$\delta(\nu\xi_2,\nu^{-1/2}\pi)+L(\nu\xi_2,\nu^{-1/2}\pi)$
\end{center}
\end{minipage}
 & $q^4-1$ & \begin{minipage}{15ex}
\begin{center}
$({\rm Ind}(\omega_{\tilde{\pi}}\Phi_1)+{\rm Ind}(\omega_{\tilde{\pi}}\Phi_3))$

or

$({\rm Ind}(\omega_{\tilde{\pi}}\xi_1'(n))+{\rm Ind}(\omega_{\tilde{\pi}}\xi_1(n)))$
\end{center}
\end{minipage}\\
\hline
\hline
\makebox{\small X} & & $\pi \rtimes \sigma$\quad (irreducible) & $q^4 - 1$ & $\tilde{\pi} \rtimes \tilde{\sigma}$ \\
\hline
\makebox{\small XI} & & \begin{minipage}{15ex}
\begin{center}
$\delta(\nu^{1/2}\pi,\nu^{-1/2}\sigma)+L(\nu^{1/2}\pi,\nu^{-1/2}\sigma)$
\end{center}
\end{minipage}
 & $q^4-1$ & \begin{minipage}{15ex}
\begin{center}
$\tilde{\sigma}{\rm Ind}(\chi_7(n))_a+\tilde{\sigma}{\rm Ind}(\chi_6(n))_a$
\end{center}
\end{minipage}\\
\end{longtable}
\end{small}
\end{center}

\begin{proof}
By Lemma \ref{Borelfixedvectors}, the dimension of $(\chi_1\times\chi_2\rtimes\sigma)^{K(1)}$ is $(q^2+1)(q+1)^2$. If $\chi_1\times\chi_2\rtimes\sigma$ is of type II, III, IV, V, or VI, then there are at least two irreducible constituents. Let $\rho$ be an irreducible constituent in one of these types. Then the representation $\mathcal{F}(\rho)$ and the dimension of $K(1)$-fixed vectors are completely determined in all cases except for the type V* and type VI cases using the Tables \ref{Vconstituents} and \ref{VIconstituents}, Lemmas \ref{KM1fixedvectors}, \ref{Borelfixedvectors}, and \ref{PQfixedvectors}, and the decompositions determined by Sally and Tadi\'{c} \cite{ST}, summarized in Tables 2.10 and 2.11 on page 39 in \cite{RS}.

The type V* and type VI representations require some additional information to determine their image under the functor $\mathcal{F}$. The issue is where to place the common factor $\sigma{\rm Ind}(\theta_9)_a$. The representation $L(\nu^{1/2}\xi_2\St_{\GL(2)},\nu^{-1/2}\sigma)$ is a common constituent of $\nu^{1/2}\xi_2\St_{\GL(2)}\rtimes\nu^{-1/2}\sigma$ and of $\nu^{1/2}\xi_2\triv_{\GL(2)}\rtimes\xi_2\nu^{-1/2}\sigma$. So $\mathcal{F}(L(\nu^{1/2}\xi_2\St_{\GL(2)},\nu^{-1/2}\sigma))$ is a common constituent of $\St_{\GL(2,q)}\rtimes\tilde{\sigma}$ and of $\triv_{\GL(2,q)}\rtimes\tilde{\sigma}$. The only possible common constituent is $\tilde{\sigma}{\rm Ind}(\theta_9)_a$. Since\\ $\mathcal{F}(L(\nu^{1/2}\xi_2\St_{\GL(2)},\nu^{-1/2}\sigma))\neq0$ by Lemma \ref{KM1fixedvectors}, it must be $\tilde{\sigma}{\rm Ind}(\theta_9)_a$. 

It is important to note that VId has a three-dimensional subspace of Iwahori subgroup-fixed vectors, see Appendix A.10 on page 297 in \cite{RS}. So the dimension of the space of $K(1)$-fixed vectors of VId is at least three-dimensional. 

Comparing our Tables \ref{Vconstituents} and \ref{VIconstituents} with the decompositions determined by Sally and Tadi\'{c} \cite{ST} completes the argument for the dimensions for types V* and VI.

Let $H$ be either the Klingen or the Siegel parabolic subgroup and let $\rho$ be $\chi\rtimes\pi$, respectively $\pi\rtimes\sigma$. In either case $\rho$ is an admissible representation of $G$. Write $H=MN$, where $M$ is its Levi subgroup and $N$ its unipotent radical. By Lemma \ref{PQfixedvectors}, the dimension of $\rho^{K(1)}$ is $q^4-1.$

When the representation $\rho$ decomposes, it decomposes into an irreducible subrepresentation, which we will call $\rho_{\rm sub}$ and an irreducible subquotient, which we will call $\rho_{\rm quot}$, each with multiplicity one. We have that $(\rho_{\rm sub},V)$ is also an admissible representation of $G$ and so, by Lemma \ref{KM1fixedvectors}, it has a non-zero subspace of $K(1)$-fixed vectors. 

By Theorem 2.5 of \cite{MoyPrasad}, there exists an $M$-associate $MU$ such that $\rho_{\rm quot}$ is a subrepresentation of the induced representation ${\rm Ind}_{MU}^\mathcal{G}$. Then, by applying the previous argument for subrepresentations, $\rho_{\rm quot}$ has a non-zero subspace of $K(1)$-fixed vectors.

\qed
\end{proof}

\begin{remark}
For types VIII, IX, IX*, and XI, it is not clear where to place the representations of $\GSp(4,q)$ obtained by applying the functor $\mathcal{F}$ to each constituent. This problem could be solved if one could prove that for generic $\rho$, we have that $\mathcal{F}(\rho)$ is also generic, but, as far as the author knows, this is an open problem.
\end{remark}

\begin{remark}
We are grateful to the referee for pointing out that, in Theorem \ref{nonsupercuspidaldimensions}, we are considering the cases where $\tilde{\chi},\tilde{\chi_1},\tilde{\chi_2},\tilde{\xi},\tilde{\xi_1},\tilde{\xi_2},\tilde{\sigma},$ and $\tilde{\pi}\cong \pi^{K_{\GL(2)}(1)}$ are non-zero. If $\rho$ is a non-supercuspidal admissible representation that doesn't satisfy the assumptions in Theorem \ref{nonsupercuspidaldimensions}, then $\mathcal{F}(\rho)$ is (a subquotient of) zero. For example, in the Klingen case, if $\tilde{\pi}=0$ or $\tilde{\chi}=0$, then $\mathcal{F}(\chi\rtimes\pi)=0$.
\end{remark}

\section{Non-cuspidal character tables}
\label{noncuspchartables}

We now give the character tables of the non-cuspidal characters that have not been explicitly given previously in this paper. These character tables were computed using Table \ref{Indalphachivalues} and the character tables in \cite{Sr} and \cite{Shinoda}. We make use of our correspondence in Table \ref{GSpIrreducibleCharactersShinoda} for some of these values. See \cite{Shinoda} for more details on these characters. Note that, in the following tables, if a conjugacy class row is omitted, then the character takes the value 0 on that class.

\begin{landscape}
\begin{small}
\begin{center}
\begin{longtable}{|c|l|l|l|}
\caption[$\sigma{\rm Ind}(\chi_6(n))_a, \sigma{\rm Ind}(\chi_7(n))_a, {\rm Ind}(\omega_\pi\Phi_1)$]{\small{$\sigma{\rm Ind}(\chi_6(n))_a, \sigma{\rm Ind}(\chi_7(n))_a, {\rm Ind}(\omega_\pi\Phi_1)$}} \label{noncuspchartable1}\\
\hline
   \multicolumn{1}{|c|}{\parbox{1cm}{\strut\centering\small\textbf{Class}\strut}} &
   \multicolumn{1}{|c|}{\parbox{4cm}{\strut\centering\small\textbf{$\sigma{\rm Ind}(\chi_6(n))_a=\chi_5(\omega)$}\strut}} &
   \multicolumn{1}{|c|}{\parbox{4cm}{\strut\centering\small\textbf{$\sigma{\rm Ind}(\chi_7(n))_a=\chi_6(\omega)$}\strut}} &
   \multicolumn{1}{|c|}{\parbox{4cm}{\strut\centering\small\textbf{${\rm Ind}(\omega_\pi\Phi_1)=\chi_7(\Lambda_1)$}\strut}} \\
\hline
\hline
\endfirsthead
\hline
   \multicolumn{1}{|c|}{\parbox{1cm}{\strut\centering\small\textbf{Class}\strut}} &
   \multicolumn{1}{|c|}{\parbox{4cm}{\strut\centering\small\textbf{$\sigma{\rm Ind}(\chi_6(n))_a=\chi_5(\omega)$}\strut}} &
   \multicolumn{1}{|c|}{\parbox{4cm}{\strut\centering\small\textbf{$\sigma{\rm Ind}(\chi_7(n))_a=\chi_6(\omega)$}\strut}} &
   \multicolumn{1}{|c|}{\parbox{4cm}{\strut\centering\small\textbf{${\rm Ind}(\omega_\pi\Phi_1)=\chi_7(\Lambda_1)$}\strut}} \\
\hline
\hline
\endhead
 \hline
\endlastfoot
$A_1(k)$ & $(q^2+1)(q-1)\sigma(\gamma^{2k})$ & $q(q-1)(q^2+1)\sigma(\gamma^{2k})$ & $(q^2+1)(q-1)\omega_\pi(\gamma^k)$ \\
\hline
$A_2(k)$ & $(q-1)\sigma(\gamma^{2k})$ & $q(q-1)\sigma(\gamma^{2k})$ & $-(q^2-q+1)\omega_\pi(\gamma^k)$ \\
\hline
$A_{31}(k)$ & $-\sigma(\gamma^{2k})$ & $-q\sigma(\gamma^{2k})$ & $(q-1)\omega_\pi(\gamma^k)$ \\
\hline
$A_{32}(k)$ & $(2q-1)\sigma(\gamma^{2k})$ & $-q\sigma(\gamma^{2k})$ & $(q-1)\omega_\pi(\gamma^k)$ \\
\hline
$A_5(k)$ & $-\sigma(\gamma^{2k})$ & $0$ & $-\omega_\pi(\gamma^k)$ \\
\hline
$B_{11}(k)$ & $-(-1)^n (q-1)^2\sigma(\gamma^{2k})$ & $(-1)^n(q-1)^2\sigma(\gamma^{2k})$ & $(q-1)(1-(-1)^t)\omega_\pi(\gamma^k)$ \\
\hline
$B_{12}(k)$ & $-(-1)^n (q^2+1)\sigma(\gamma^{2k+1})$ & $-(-1)^n(q^2+1)\sigma(\gamma^{2k+1})$ & $0$ \\
\hline
$B_{21}(k)$ & $(q-1)\sigma(-\gamma^{2k})$ & $(q-1)\sigma(-\gamma^{2k})$ & $0$ \\
\hline
$B_{22}(k)$ & $(q-1-\tilde{\omega}(\gamma^{1/2})+\tilde{\omega}(\gamma^{q/2}))\sigma(-\gamma^{2k+1})$ & $-(q-1+q(\tilde{\omega}(\gamma^{1/2})+\tilde{\omega}(\gamma^{q/2})))\sigma(-\gamma^{2k+1})$ & $(q-1)(\Lambda_1(\gamma^{1/2})+\Lambda_1(\gamma^{q/2}))\omega_\pi(\gamma^k)$\\
\hline
$B_3(k)$ & $(-1)^n(q-1)\sigma(\gamma^{2k})$ & $(-1)^n(q-1)\sigma(\gamma^{2k})$ & $((q-1)\omega_\pi(-1)-1)\omega_\pi(\gamma^k)$ \\
\hline
$B_{41}(k)$ & $-(-1)^n \sigma(\gamma^{2k})$ & $-(-1)^n \sigma(\gamma^{2k})$ & $-(1+\omega_\pi(-1))\omega_\pi(\gamma^k)$ \\
\hline
$B_{42}(k)$ & $-(-1)^n \sigma(\gamma^{2k})$ & $-(-1)^n \sigma(\gamma^{2k})$ & $-(1+\omega_\pi(-1))\omega_\pi(\gamma^k)$ \\
\hline
$B_{43}(k)$ & $-(-1)^n\sigma(\gamma^{2k+1})$ & $-(-1)^n\sigma(\gamma^{2k+1})$  & $0$\\
\hline
$B_{44}(k)$ & $-(-1)^n\sigma(\gamma^{2k+1})$  & $-(-1)^n\sigma(\gamma^{2k+1})$  & $0$\\
\hline
$B_{51}(k)$ & $-\sigma(-\gamma^{2k})$ & $-\sigma(-\gamma^{2k})$ & $0$ \\
\hline
$B_{52}(k)$ & $(-1-\tilde{\omega}(\gamma^{1/2})-\tilde{\omega}(\gamma^{q/2}))\sigma(-\gamma^{2k+1})$ & $\sigma(-\gamma^{2k+1})$ & $-(\Lambda_1(\gamma^{1/2})+\Lambda_1(\gamma^{q/2}))\omega_\pi(\gamma^k)$\\
\hline
$C_1(i,k)$ & $(q-1)\sigma(\gamma^{2k+i})$ & $(q-1)\sigma(\gamma^{2k+i})$ & $0$ \\
\hline
$C_{21}(i,k)$ & $0$ & $0$ & $0$ \\
\hline
$C_{22}(i,k)$ & $-(\tilde{\omega}(\gamma^{1/2})+\tilde{\omega}(\gamma^{q/2}))\sigma(-\gamma^{2k+i+1})$ & $-(\tilde{\omega}(\gamma^{1/2})+\tilde{\omega}(\gamma^{q/2}))\sigma(-\gamma^{2k+i+1})$ & $0$\\
\hline
$C_3(i,k)$ & $-\sigma(\gamma^{2k+i})$ & $-\sigma(\gamma^{2k+i})$ & $0$ \\
\hline
$C_4(i,k)$ & $0$ & $0$ & $-\omega_\pi(\gamma^k)$ \\
\hline
$C_5(i,k)$ & $0$ & $0$ & $(q-1)\omega_\pi(\gamma^k)$ \\
\hline
$D_1(i,k)$ & $0$ & $0$ & $-2(-1)^i\omega_\pi(\gamma^k)$ \\
\hline
$D_2(i,k)$ & $(q-1-\omega(\eta^i)-\omega(\eta^{-i}))\sigma(\gamma^{2k+i})$ & $-(q-1+\omega(\eta^i)+\omega(\eta^{-i}))\sigma(\gamma^{2k+i})$ & $2(-1)^i(q-1)\omega_\pi(\gamma^k)$\\
\hline
$D_{31}(i,k)$ & $-(\omega(\eta^i)+\omega(\eta^{-i}))\sigma(-\gamma^{2k+i})$ & $-(\omega(\eta^i)+\omega(\eta^{-i}))\sigma(-\gamma^{2k+i})$ & $0$\\
\hline
$D_{32}(i,k)$ & $(-2(-1)^n-\omega(\eta^{2i})-\omega(\eta^{-2i}))\sigma(\gamma^{2k+1})$ & $(2(-1)^n+\omega(\eta^{2i})+\omega(\eta^{-2i}))\sigma(\gamma^{2k+1})$  & $-(-1)^i(1+\Lambda_1(-1))(\Lambda_1(\gamma^{1/2})$\\
& & & $\quad+\Lambda_1(\gamma^{q/2}))\omega_\pi(\gamma^k)$\\
\hline
$D_4(i,j,k)$ & $-(\omega(\eta^i)+\omega(\eta^{-i}))\sigma(\gamma^{2k+i+j})$ & $-(\omega(\eta^i)+\omega(\eta^{-i}))\sigma(\gamma^{2k+i+j})$ & $0$\\
\hline
$D_5(i,k)$ & $-(1+\omega(\eta^i)+\omega(\eta^{-i}))\sigma(\gamma^{2k+i})$ & $\sigma(\gamma^{2k+i})$ & $-2(-1)^i\omega_\pi(\gamma^k)$\\
\hline
$D_6(i,k)$ & $(q-1)(\omega(\eta^i)+\omega(\eta^{-i}))\sigma(\gamma^{2k})$ & $-(q-1)(\omega(\eta^i)+\omega(\eta^{-i}))\sigma(\gamma^{2k})$ & $(q-1-2(-1)^i)\omega_\pi(\gamma^k)$ \\
\hline
$D_7(i,j,k)$ & $0$ & $0$ & $-2((-1)^i+(-1)^j)\omega_\pi(\gamma^k)$\\
\hline
$D_8(i,k)$ & $-(\omega(\eta^i)+\omega(\eta^{-i}))\sigma(\gamma^{2k})$ & $(\omega(\eta^i)+\omega(\eta^{-i}))\sigma(\gamma^{2k})$ & $-(2(-1)^i+1)\omega_\pi(\gamma^k)$ \\
\end{longtable}
\end{center}

\begin{center}
\begin{longtable}{|c|l|l|l|}
\caption[${\rm Ind}(\omega_\pi\Phi_3), \sigma{\rm Ind}(\Phi_9)_a, \sigma{\rm Ind}(\theta_1)$]{\small{${\rm Ind}(\omega_\pi\Phi_3), \sigma{\rm Ind}(\Phi_9)_a, \sigma{\rm Ind}(\theta_1)$}} \label{noncuspchartable1}\\
\hline
   \multicolumn{1}{|c|}{\parbox{1cm}{\strut\centering\small\textbf{Class}\strut}} &
   \multicolumn{1}{|c|}{\parbox{4cm}{\strut\centering\small\textbf{${\rm Ind}(\omega_\pi\Phi_3)=\chi_8(\Lambda_1)$}\strut}} &
   \multicolumn{1}{|c|}{\parbox{4cm}{\strut\centering\small\textbf{$\sigma{\rm Ind}(\Phi_9)_a=\sigma\tau_2$}\strut}} &
   \multicolumn{1}{|c|}{\parbox{4cm}{\strut\centering\small\textbf{$\sigma{\rm Ind}(\theta_1)=\sigma\tau_3$}\strut}} \\
\hline
\hline
\endfirsthead
\hline
   \multicolumn{1}{|c|}{\parbox{1cm}{\strut\centering\small\textbf{Class}\strut}} &
   \multicolumn{1}{|c|}{\parbox{4cm}{\strut\centering\small\textbf{${\rm Ind}(\omega_\pi\Phi_3)=\chi_8(\Lambda_1)$}\strut}} &
   \multicolumn{1}{|c|}{\parbox{4cm}{\strut\centering\small\textbf{$\sigma{\rm Ind}(\Phi_9)_a=\sigma\tau_2$}\strut}} &
   \multicolumn{1}{|c|}{\parbox{4cm}{\strut\centering\small\textbf{$\sigma{\rm Ind}(\theta_1)=\sigma\tau_3$}\strut}} \\
\hline
\hline
\endhead
 \hline
\endlastfoot
$A_1(k)$ & $q(q^2+1)(q-1)\omega_\pi(\gamma^k)$ & $q(q^2+1)\sigma(\gamma^{2k})$ & $q^2(q^2+1)\sigma(\gamma^{2k})$ \\
\hline
$A_2(k)$ & $-q\omega_\pi(\gamma^k)$ & $q\sigma(\gamma^{2k})$ & $q^2\sigma(\gamma^{2k})$ \\
\hline
$A_{31}(k)$ & $0$ & $q\sigma(\gamma^{2k})$ & $0$ \\
\hline
$A_{32}(k)$ & $-2q\omega_\pi(\gamma^k)$ & $q\sigma(\gamma^{2k})$ & $0$ \\
\hline
$B_{11}(k)$ & $q(q-1)(1+\omega_\pi(-1))\omega_\pi(\gamma^k)$ & $(-1)^t(q^2+1)\sigma(\gamma^{2k})$ & $2(-1)^tq\sigma(\gamma^{2k})$ \\
\hline
$B_{12}(k)$ & $0$ & $-(-1)^t(q^2+1)\sigma(\gamma^{2k+1})$ & $0$\\
\hline
$B_{21}(k)$ & $0$ & $(q+(-1)^t)\sigma(-\gamma^{2k})$ & $(1+(-1)^t)q\sigma(-\gamma^{2k})$ \\
\hline
$B_{22}(k)$ & $-(q-1)(\Lambda_1(\gamma^{1/2})+\Lambda_1(\gamma^{q/2}))\omega_\pi(\gamma^k)$ & $(q+(-1)^t)\sigma(-\gamma^{2k+1})$ & $-q(1+(-1)^t)\sigma(-\gamma^{2k+1})$\\
\hline
$B_3(k)$ & $((q-1)\omega_\pi(-1)-1)\omega_\pi(\gamma^k)$ & $(-1)^t\sigma(\gamma^{2k})$ & $(-1)^t q\sigma(\gamma^{2k})$ \\
\hline
$B_{41}(k)$ & $0$ & $(-1)^t\sigma(\gamma^{2k})$ & $0$ \\
\hline
$B_{42}(k)$ & $0$ & $(-1)^t\sigma(\gamma^{2k})$ & $0$ \\
\hline
$B_{43}(k)$ & $0$ & $-(-1)^t\sigma(\gamma^{2k+1})$ & $0$\\
\hline
$B_{44}(k)$ & $0$ & $-(-1)^t\sigma(\gamma^{2k+1})$ & $0$\\
\hline
$B_{51}(k)$ & $0$ & $(1+(-1)^t)\sigma(-\gamma^{2k})$ & $0$ \\
\hline
$B_{52}(k)$ & $(\Lambda_1(\gamma^{1/2})+\Lambda_1(\gamma^{q/2}))\omega_\pi(\gamma^k)$ & $(-1)^t\sigma(-\gamma^{2k+1})$ & $0$ \\
\hline
$C_1(i,k)$ & $0$ & $(q+(-1)^i)\sigma(\gamma^{2k+i})$ & $(1+(-1)^i)q\sigma(\gamma^{2k+i})$ \\
\hline
$C_{21}(i,k)$ & $0$ & $((-1)^t+(-1)^i)\sigma(-\gamma^{2k+i})$ & $((-1)^t+(-1)^i)\sigma(-\gamma^{2k+i})$ \\
\hline
$C_{22}(i,k)$ & $0$ & $((-1)^t-(-1)^i)\sigma(-\gamma^{2k+i+1})$ & $((-1)^t+(-1)^i)\sigma(-\gamma^{2k+i+1})$ \\
\hline
$C_3(i,k)$ & $0$ & $(-1)^i\sigma(\gamma^{2k+i})$ & $0$ \\
\hline
$C_4(i,k)$ & $-\omega_\pi(\gamma^k)$ & $(-1)^i\sigma(\gamma^{2k})$ & $(-1)^i\sigma(\gamma^{2k})$ \\
\hline
$C_5(i,k)$ & $(q-1)\omega_\pi(\gamma^k)$ & $(-1)^i(q+1)\sigma(\gamma^{2k})$ & $(-1)^i(q+1)\sigma(\gamma^{2k})$ \\
\hline
$C_6(i,j,k)$ & $0$ & $((-1)^i+(-1)^j)\sigma(\gamma^{2k+i+j})$ & $((-1)^i+(-1)^j)\sigma(\gamma^{2k+i+j})$ \\
\hline
$D_1(i,k)$ & $-2(-1)^i\omega_\pi(\gamma^k)$ & $0$ & $0$ \\
\hline
$D_2(i,k)$ & $-2(-1)^i(q-1)\omega_\pi(\gamma^k)$ & $(q-(-1)^i)\sigma(\gamma^{2k+1})$ &$-q(1+(-1)^i)\sigma(\gamma^{2k+1})$ \\
\hline
$D_{31}(i,k)$ & $0$ & $(-1)^i((-1)^t-1)\sigma(-\gamma^{2k+i})$ & $-(-1)^i(1+(-1)^t)\sigma(-\gamma^{2k+i})$ \\
\hline
$D_{32}(i,k)$ & $(-1)^i(1+\Lambda_1(-1))(\Lambda_1(\gamma^{1/2})+\Lambda_1(\gamma^{q/2}))\omega_\pi(\gamma^k)$ & $-(1+\omega_0(-1))\sigma(\gamma^{2k+1})$ & $-(1+\omega_0(-1))\sigma(\gamma^{2k+1})$ \\
\hline
$D_4(i,j,k)$ & $0$ & $(-1)^i((-1)^j-1)\sigma(\gamma^{2k+i+j})$ & $(-1)^i(1+(-1)^j)\sigma(\gamma^{2k+i+j})$\\
\hline
$D_5(i,k)$ & $2(-1)^i\omega_\pi(\gamma^k)$ & $-(-1)^i\sigma(\gamma^{2k+1})$ & $0$\\
\hline
$D_6(i,k)$ & $-(q-1+2(-1)^i q)\omega_\pi(\gamma^k)$ & $(-1)^i(q-1)\sigma(\gamma^{2k})$ & $-(-1)^i(q-1)\sigma(\gamma^{2k})$\\
\hline
$D_7(i,j,k)$ & $2((-1)^i+(-1)^j)\omega_\pi(\gamma^k)$ & $-2(-1)^{i+j}\omega_0(-1)\sigma(\gamma^{2k})$ & $2(-1)^{i+j}\omega_0(-1)\sigma(\gamma^{2k})$\\
\hline
$D_8(i,k)$ & $\omega_\pi(\gamma^k)$ & $-(-1)^i\sigma(\gamma^{2k})$ & $(-1)^i\sigma(\gamma^{2k})$\\
\end{longtable}
\end{center}

\begin{center}
\begin{longtable}{|c|l|l|l|}
\caption[$\sigma{\rm Ind}(\theta_3), {\rm Ind}(\xi\omega_\pi\theta_5), {\rm Ind}(\xi\omega_\pi\theta_7)$]{\small{$\sigma{\rm Ind}(\theta_3), {\rm Ind}(\xi\omega_\pi\theta_5), {\rm Ind}(\xi\omega_\pi\theta_7)$}} \label{noncuspchartable1}\\
\hline
   \multicolumn{1}{|c|}{\parbox{1cm}{\strut\centering\small\textbf{Class}\strut}} &
   \multicolumn{1}{|c|}{\parbox{4cm}{\strut\centering\small\textbf{$\sigma{\rm Ind}(\theta_3)=\sigma\tau_1$}\strut}} &
   \multicolumn{1}{|c|}{\parbox{4cm}{\strut\centering\small\textbf{${\rm Ind}(\xi\omega_\pi\theta_5)=\tau_5(\lambda')$}\strut}} &
   \multicolumn{1}{|c|}{\parbox{4cm}{\strut\centering\small\textbf{${\rm Ind}(\xi\omega_\pi\theta_7)=\tau_4(\lambda')$}\strut}} \\
\hline
\hline
\endfirsthead
\hline
   \multicolumn{1}{|c|}{\parbox{1cm}{\strut\centering\small\textbf{Class}\strut}} &
   \multicolumn{1}{|c|}{\parbox{4cm}{\strut\centering\small\textbf{$\sigma{\rm Ind}(\theta_3)=\sigma\tau_1$}\strut}} &
   \multicolumn{1}{|c|}{\parbox{4cm}{\strut\centering\small\textbf{${\rm Ind}(\xi\omega_\pi\theta_5)=\tau_5(\lambda')$}\strut}} &
   \multicolumn{1}{|c|}{\parbox{4cm}{\strut\centering\small\textbf{${\rm Ind}(\xi\omega_\pi\theta_7)=\tau_4(\lambda')$}\strut}} \\
\hline
\hline
\endhead
 \hline
\endlastfoot
$A_1(k)$ & $(q^2+1)\sigma(\gamma^{2k})$ & $(-1)^kq^2(q^2-1)\omega_\pi(\gamma^k)$ & $(-1)^k(q^2-1)\omega_\pi(\gamma^k)$ \\
\hline
$A_2(k)$ & $\sigma(\gamma^{2k})$ & $-(-1)^k q^2\omega_\pi(\gamma^k)$ & $-(-1)^k\omega_\pi(\gamma^k)$ \\
\hline
$A_{31}(k)$ & $(q+1)\sigma(\gamma^{2k})$ & $0$ & $(-1)^k(q-1)\omega_\pi(\gamma^k)$ \\
\hline
$A_{32}(k)$ & $-(q-1)\sigma(\gamma^{2k})$ & $0$ & $-(-1)^k(q+1)\omega_\pi(\gamma^k)$ \\
\hline
$A_5(k)$ & $\sigma(\gamma^{2k})$ & $0$ & $-(-1)^k\omega_\pi(\gamma^k)$ \\
\hline
$B_{11}(k)$ & $2(-1)^tq\sigma(\gamma^{2k})$ & $0$ & $0$ \\
\hline
$B_{21}(k)$ & $(1+(-1)^t)\sigma(-\gamma^{2k})$ & $0$ & $0$ \\
\hline
$B_{22}(k)$ & $(1+(-1)^t)\sigma(-\gamma^{2k+1})$ & $0$ & $0$ \\
\hline
$B_3(k)$ & $(-1)^t q\sigma(\gamma^{2k})$ & $(-1)^{t+k}q\omega_\pi(\gamma^k)$ & $-(-1)^{t+k}q\omega_\pi(\gamma^k)$ \\
\hline
$B_{51}(k)$ & $(1+(-1)^t)\sigma(-\gamma^{2k})$ & $0$ & $0$ \\
\hline
$B_{52}(k)$ & $(1+(-1)^t)\sigma(-\gamma^{2k+1})$ & $0$ & $0$ \\
\hline
$C_1(i,k)$ & $(1+(-1)^i)\sigma(\gamma^{2k+i})$ & $0$ & $0$ \\
\hline
$C_{21}(i,k)$ & $((-1)^i+(-1)^t)\sigma(-\gamma^{2k+i})$ & $0$ & $0$ \\
\hline
$C_{22}(i,k)$ & $-(-1)^t(1+(-1)^i)\sigma(-\gamma^{2k+i+1})$ & $0$ & $0$ \\
\hline
$C_3(i,k)$ & $(1+(-1)^i)\sigma(\gamma^{2k+i})$ & $0$ & $0$ \\
\hline
$C_4(i,k)$ & $(-1)^i\sigma(\gamma^{2k})$ & $(-1)^{i+k+1}\omega_\pi(\gamma^k)$ & $(-1)^{i+k+1}\omega_\pi(\gamma^k)$ \\
\hline
$C_5(i,k)$ & $(-1)^i (q+1)\sigma(\gamma^{2k})$ & $(-1)^{i+k}(q-1)\omega_\pi(\gamma^k)$ & $(-1)^{i+k}(q-1)\omega_\pi(\gamma^k)$ \\
\hline
$C_6(i,j,k)$ & $((-1)^i+(-1)^j)\sigma(\gamma^{2k+i+j})$ & $0$ & $0$ \\
\hline
$D_1(i,k)$ & $0$ & $-(-1)^k$ & $-(-1)^k(\tilde{\lambda'}+\tilde{\lambda'}^q)(\gamma^k\theta^i)$ \\
\hline
$D_2(i,k)$ & $(1+(-1)^i)\sigma(\gamma^{2k+i})$ & $0$ & $0$ \\
\hline
$D_{31}(i,k)$ & $(-1)^i(1+(-1)^t)\sigma(-\gamma^{2k+i})$ & $0$ & $0$ \\
\hline
$D_{32}(i,k)$ & $(1+\omega_0(-1))\sigma(\gamma^{2k+1})$ & $0$ & $0$ \\
\hline
$D_4(i,j,k)$ & $(-1)^i(1+(-1)^j)\sigma(\gamma^{2k+i+j})$ & $0$ & $0$ \\
\hline
$D_5(i,k)$ & $(1+(-1)^i)\sigma(\gamma^{2k+1})$ & $0$ & $0$ \\
\hline
$D_6(i,k)$ & $-(-1)^i(q-1)\sigma(\gamma^{2k})$ & $-(-1)^{i+k}(q+1)\omega_\pi(\gamma^k)$ & $-(-1)^{i+k}(q+1)\omega_\pi(\gamma^k)$ \\
\hline
$D_7(i,j,k)$ & $2(-1)^{i+j}\omega_0(-1)\sigma(\gamma^{2k})$ & $0$ & $0$\\
\hline
$D_8(i,k)$ & $(-1)^i\sigma(\gamma^{2k})$ & $-(-1)^{i+k}\omega_\pi(\gamma^k)$ & $-(-1)^{i+k}\omega_\pi(\gamma^k)$ \\
\hline
$D_9(i,k)$ & $0$ & $(-1)^k(\lambda'(\zeta^i)+\lambda'(\zeta^{qi}))\omega_\pi(\gamma^k)$ & $-(-1)^k(\lambda'(\zeta^i)+\lambda'(\zeta^{qi}))\omega_\pi(\gamma^k)$\\
\end{longtable}
\end{center}

\begin{center}
\begin{longtable}{|c|l|l|l|l|}
\caption[$\sigma{\rm Ind}(\theta_9)_a, \sigma{\rm Ind}(\theta_{11})_a, \sigma{\rm Ind}(\theta_{12})_a, \sigma\St_{\GSp(4)}$]{\small{$\sigma{\rm Ind}(\theta_9)_a, \sigma{\rm Ind}(\theta_{11})_a, \sigma{\rm Ind}(\theta_{12})_a, \sigma\St_{\GSp(4)}$}} \label{noncuspchartable1}\\
\hline
   \multicolumn{1}{|c|}{\parbox{1cm}{\strut\centering\small\textbf{Class}\strut}} &
   \multicolumn{1}{|c|}{\parbox{4cm}{\strut\centering\small\textbf{$\sigma{\rm Ind}(\theta_9)_a=\sigma\theta_1$}\strut}} &
   \multicolumn{1}{|c|}{\parbox{4cm}{\strut\centering\small\textbf{$\sigma{\rm Ind}(\theta_{11})_a=\sigma\theta_3$}\strut}} &
   \multicolumn{1}{|c|}{\parbox{4cm}{\strut\centering\small\textbf{$\sigma{\rm Ind}(\theta_{12})_a=\sigma\theta_4$}\strut}} &
   \multicolumn{1}{|c|}{\parbox{4cm}{\strut\centering\small\textbf{$\sigma\St_{\GSp(4)}=\theta_5$}\strut}} \\
\hline
\hline
\endfirsthead
\hline
   \multicolumn{1}{|c|}{\parbox{1cm}{\strut\centering\small\textbf{Class}\strut}} &
   \multicolumn{1}{|c|}{\parbox{4cm}{\strut\centering\small\textbf{$\sigma{\rm Ind}(\theta_9)_a=\sigma\theta_1$}\strut}} &
   \multicolumn{1}{|c|}{\parbox{4cm}{\strut\centering\small\textbf{$\sigma{\rm Ind}(\theta_{11})_a=\sigma\theta_3$}\strut}} &
   \multicolumn{1}{|c|}{\parbox{4cm}{\strut\centering\small\textbf{$\sigma{\rm Ind}(\theta_{12})_a=\sigma\theta_4$}\strut}} &
   \multicolumn{1}{|c|}{\parbox{4cm}{\strut\centering\small\textbf{$\sigma\St_{\GSp(4)}=\theta_5$}\strut}} \\
\hline
\hline
\endhead
 \hline
\endlastfoot
$A_1(k)$ & $\frac{1}{2}q(q+1)^2\sigma(\gamma^{2k})$ & $\frac{1}{2}q(q^2+1)\sigma(\gamma^{2k})$ & $\frac{1}{2}q(q^2+1)\sigma(\gamma^{2k})$ & $q^4\sigma(\gamma^{2k})$ \\
\hline
$A_2(k)$ & $\frac{1}{2}q(q+1)\sigma(\gamma^{2k})$ & $\frac{1}{2}q(1-q)\sigma(\gamma^{2k})$ & $\frac{1}{2}q(q+1)\sigma(\gamma^{2k})$ & $0$ \\
\hline
$A_{31}(k)$ & $q\sigma(\gamma^{2k})$ & $q\sigma(\gamma^{2k})$ & $0$ & $0$ \\
\hline
$A_{32}(k)$ & $0$ & $0$ & $q\sigma(\gamma^{2k})$ & $0$ \\
\hline
$B_{11}(k)$ & $\frac{1}{2}(q+1)^2\sigma(\gamma^{2k})$ & $\frac{1}{2}(q^2+2q-1)\sigma(\gamma^{2k})$ & $-\frac{1}{2}(q^2-2q-1)\sigma(\gamma^{2k})$ & $q^2\sigma(\gamma^{2k})$ \\
\hline
$B_{12}(k)$ & $\frac{1}{2}(q^2-1)\sigma(\gamma^{2k+1})$ & $\frac{1}{2}(q^2+1)\sigma(\gamma^{2k+1})$ & $-\frac{1}{2}(q^2+1)\sigma(\gamma^{2k+1})$ & $-q^2\sigma(\gamma^{2k+1})$\\
\hline
$B_{21}(k)$ & $(q+1)\sigma(-\gamma^{2k})$ & $\sigma(-\gamma^{2k})$ & $q\sigma(-\gamma^{2k})$ & $q\sigma(-\gamma^{2k})$ \\
\hline
$B_{22}(k)$ & $0$ & $q\sigma(-\gamma^{2k+1})$ & $-\sigma(-\gamma^{2k+1})$ & $-q\sigma(-\gamma^{2k+1})$ \\
\hline
$B_3(k)$ & $\frac{1}{2}(q+1)\sigma(\gamma^{2k})$ & $\frac{1}{2}(q-1)\sigma(\gamma^{2k})$ & $\frac{1}{2}(q+1)\sigma(\gamma^{2k})$ & $0$ \\
\hline
$B_{41}(k)$ & $\frac{1}{2}(q+1)\sigma(\gamma^{2k})$ & $-\frac{1}{2}(q+1)\sigma(\gamma^{2k})$ & $-\frac{1}{2}(q-1)\sigma(\gamma^{2k})$ & $0$ \\
\hline
$B_{42}(k)$ & $-\frac{1}{2}(q-1)\sigma(\gamma^{2k})$ & $\frac{1}{2}(q-1)\sigma(\gamma^{2k})$ & $\frac{1}{2}(q+1)\sigma(\gamma^{2k})$ & $0$ \\
\hline
$B_{43}(k)$ & $-\frac{1}{2}(q+1)\sigma(\gamma^{2k+1})$ & $\frac{1}{2}(q+1)\sigma(\gamma^{2k+1})$ & $\frac{1}{2}(q-1)\sigma(\gamma^{2k+1})$ & $0$ \\
\hline
$B_{44}(k)$ & $\frac{1}{2}(q-1)\sigma(\gamma^{2k+1})$ & $-\frac{1}{2}(q-1)\sigma(\gamma^{2k+1})$ & $-\frac{1}{2}(q+1)\sigma(\gamma^{2k+1})$ & $0$ \\
\hline
$B_{51}(k)$ & $\sigma(-\gamma^{2k})$ & $\sigma(-\gamma^{2k})$ & $0$ & $0$ \\
\hline
$B_{52}(k)$ & $0$ & $0$ & $-\sigma(-\gamma^{2k+1})$ & $0$ \\
\hline
$C_1(i,k)$ & $(q+1)\sigma(\gamma^{2k+i})$ & $\sigma(\gamma^{2k+i})$ & $q\sigma(\gamma^{2k+i})$ & $q\sigma(\gamma^{2k+i})$ \\
\hline
$C_{21}(i,k)$ & $2\sigma(-\gamma^{2k+i})$ & $\sigma(-\gamma^{2k+i})$ & $\sigma(-\gamma^{2k+i})$ & $\sigma(-\gamma^{2k+i})$ \\
\hline
$C_{22}(i,k)$ & $0$ & $\sigma(-\gamma^{2k+i+1})$ & $-\sigma(-\gamma^{2k+i+1})$ & $-\sigma(-\gamma^{2k+i+1})$\\
\hline
$C_3(i,k)$ & $\sigma(\gamma^{2k+i})$ & $\sigma(\gamma^{2k+i})$ & $0$ & $0$ \\
\hline
$C_4(i,k)$ & $\sigma(\gamma^{2k})$ & $0$ & $\sigma(\gamma^{2k})$ & $0$ \\
\hline
$C_5(i,k)$ & $(q+1)\sigma(\gamma^{2k})$ & $q\sigma(\gamma^{2k})$ & $\sigma(\gamma^{2k})$ & $q\sigma(\gamma^{2k})$ \\
\hline
$C_6(i,j,k)$ & $2\sigma(\gamma^{2k+i+j})$ & $\sigma(\gamma^{2k+i+j})$ & $\sigma(\gamma^{2k+i+j})$ & $\sigma(\gamma^{2k+i+j})$ \\
\hline
$D_1(i,k)$ & $0$ & $-\sigma(\gamma^{2k+i})$ & $\sigma(\gamma^{2k+i})$ & $-\sigma(\gamma^{2k+i})$ \\
\hline
$D_2(i,k)$ & $0$ & $q\sigma(\gamma^{2k+i})$ & $-\sigma(\gamma^{2k+i})$ & $-q\sigma(\gamma^{2k+i})$ \\
\hline
$D_{31}(i,k)$ & $0$ & $\sigma(-\gamma^{2k+i})$ & $-\sigma(-\gamma^{2k+i})$ & $-\sigma(-\gamma^{2k+i})$ \\
\hline
$D_{32}(i,k)$ & $0$ & $-\sigma(\gamma^{2k+1})$ & $-\sigma(\gamma^{2k+1})$ & $\sigma(\gamma^{2k+1})$ \\
\hline
$D_4(i,j,k)$ & $0$ & $\sigma(\gamma^{2k+i+j})$ & $-\sigma(\gamma^{2k+i+j})$ & $-\sigma(\gamma^{2k+i+j})$ \\
\hline
$D_5(i,k)$ & $0$ & $0$ & $-\sigma(\gamma^{2k+1})$ & $0$\\
\hline
$D_6(i,k)$ & $0$ & $-\sigma(\gamma^{2k})$ & $q\sigma(\gamma^{2k})$ & $-q\sigma(\gamma^{2k})$\\
\hline
$D_7(i,j,k)$ & $0$ & $-\sigma(\gamma^{2k})$ & $-\sigma(\gamma^{2k})$ & $\sigma(\gamma^{2k})$\\
\hline
$D_8(i,k)$ & $0$ & $-\sigma(\gamma^{2k})$ & $0$ & $0$ \\
\hline
$D_9(i,k)$ & $-\sigma(\gamma^{2k})$ & $0$ & $0$ & $\sigma(\gamma^{2k})$\\
\end{longtable}
\end{center}

\begin{center}
\begin{longtable}{|c|l|l|l|l|}
\caption[${\rm Ind}(\omega_\pi\xi_1(n))_a, {\rm Ind}(\omega_\pi\xi_1'(n))_a, {\rm Ind}(\chi\omega_\pi\xi_{21}(n)), {\rm Ind}(\chi\omega_\pi\xi_{41}'(n))$]{\small{${\rm Ind}(\omega_\pi\xi_1(n))_a, {\rm Ind}(\omega_\pi\xi_1'(n))_a, {\rm Ind}(\chi\omega_\pi\xi_{21}(n)), {\rm Ind}(\chi\omega_\pi\xi_{41}'(n))$}} \label{noncuspchartable1}\\
\hline
   \multicolumn{1}{|c|}{\parbox{1cm}{\strut\centering\small\textbf{Class}\strut}} &
   \multicolumn{1}{|c|}{\parbox{4cm}{\strut\centering\small\textbf{${\rm Ind}(\omega_\pi\xi_1(n))_a=\chi_7(\Lambda)$}\strut}} &
   \multicolumn{1}{|c|}{\parbox{4cm}{\strut\centering\small\textbf{${\rm Ind}(\omega_\pi\xi_1'(n))_a=\chi_8(\Lambda)$}\strut}} &
   \multicolumn{1}{|c|}{\parbox{4cm}{\strut\centering\small\textbf{${\rm Ind}(\chi\omega_\pi\xi_{21}(n))=X_3(\Lambda,\xi)$}\strut}} &
   \multicolumn{1}{|c|}{\parbox{4cm}{\strut\centering\small\textbf{${\rm Ind}(\chi\omega_\pi\xi_{41}'(n))=X_3(\Lambda_1,\xi)$}\strut}} \\
\hline
\hline
\endfirsthead
\hline
   \multicolumn{1}{|c|}{\parbox{1cm}{\strut\centering\small\textbf{Class}\strut}} &
   \multicolumn{1}{|c|}{\parbox{4cm}{\strut\centering\small\textbf{${\rm Ind}(\omega_\pi\xi_1(n))_a=\chi_7(\Lambda)$}\strut}} &
   \multicolumn{1}{|c|}{\parbox{4cm}{\strut\centering\small\textbf{${\rm Ind}(\omega_\pi\xi_1'(n))_a=\chi_8(\Lambda)$}\strut}} &
   \multicolumn{1}{|c|}{\parbox{4cm}{\strut\centering\small\textbf{${\rm Ind}(\chi\omega_\pi\xi_{21}(n))=X_3(\Lambda,\xi)$}\strut}} &
   \multicolumn{1}{|c|}{\parbox{4cm}{\strut\centering\small\textbf{${\rm Ind}(\chi\omega_\pi\xi_{41}'(n))=X_3(\Lambda_1,\xi)$}\strut}} \\
\hline
\hline
\endhead
 \hline
\endlastfoot
$A_1(k)$ & $(q-1)(q^2+1)\omega_\pi(\gamma^k)$ & $q(q-1)(q^2+1)\omega_\pi(\gamma^k)$ & $(q^4-1)\chi(\gamma^k)\omega_\pi(\gamma^k)$ & $(q^4-1)\chi(\gamma^k)\omega_\pi(\gamma^k)$ \\
\hline
$A_2(k)$ & $-(q^2-q+1)\omega_\pi(\gamma^k)$ & $-q\omega_\pi(\gamma^k)$ & $-(q^2+1)\chi(\gamma^k)\omega_\pi(\gamma^k)$ & $-(q^2+1)\chi(\gamma^k)\omega_\pi(\gamma^k)$ \\
\hline
$A_{31}(k)$ & $(q-1)\omega_\pi(\gamma^k)$ & $0$ & $(q-1)\chi(\gamma^k)\omega_\pi(\gamma^k)$ & $(q-1)\chi(\gamma^k)\omega_\pi(\gamma^k)$ \\
\hline
$A_{32}(k)$ & $(q-1)\omega_\pi(\gamma^k)$ & $-2q\omega_\pi(\gamma^k)$ & $-(q+1)\chi(\gamma^k)\omega_\pi(\gamma^k)$ & $-(q+1)\chi(\gamma^k)\omega_\pi(\gamma^k)$ \\
\hline
$A_5(k)$ & $-\omega_\pi(\gamma^k)$ & $0$ & $-\chi(\gamma^k)\omega_\pi(\gamma^k)$ & $-\chi(\gamma^k)\omega_\pi(\gamma^k)$ \\
\hline
$B_{11}(k)$ & $(q-1)(1+\omega_\pi(-1))\omega_\pi(\gamma^k)$ & $q(q-1)(1+\omega_\pi(-1))\omega_\pi(\gamma^k)$ & $(-1)^k(q^2-1)((-1)^t$ & $(-1)^k(q^2-1)((-1)^t$ \\
& & & $+\chi(-1)\omega_\pi(-1))\chi(\gamma^k)\omega_\pi(\gamma^k)$ & $+\chi(-1)\omega_\pi(-1))\chi(\gamma^k)\omega_\pi(\gamma^k)$\\
\hline
$B_{22}(k)$ & $(q-1)(\Lambda(\gamma^{1/2})+\Lambda(\gamma^{q/2}))\cdot$ & $-(q-1)(\Lambda(\gamma^{1/2})+\Lambda(\gamma^{q/2}))\cdot$ & $0$ & $0$ \\
& \qquad$\omega_\pi(\gamma^k)$ & \qquad$\omega_\pi(\gamma^k)$ & & \\
\hline
$B_3(k)$ & $((q-1)\omega_\pi(-1)-1)\omega_\pi(\gamma^k)$ & $-q\omega_\pi(-\gamma^k)$ & $(-1)^k((q-1)\chi(-1)\omega_\pi(-1)$ & $(-1)^k((q-1)\chi(-1)\omega_\pi(-1)$ \\
& & & $-(-1)^t(q+1))\chi(\gamma^k)\omega_\pi(\gamma^k)$ & $-(-1)^t(q+1))\chi(\gamma^k)\omega_\pi(\gamma^k)$\\
\hline
$B_{41}(k)$ & $-(\omega_\pi(-1)+1)\omega_\pi(\gamma^k)$ & $0$ & $-(-1)^k((-1)^t+\chi(-1)\cdot$ & $-(-1)^k((-1)^t+\chi(-1)\cdot$ \\
& & & \qquad$\omega_\pi(-1))\chi(\gamma^k)\omega_\pi(\gamma^k)$ & \qquad$\omega_\pi(-1))\chi(\gamma^k)\omega_\pi(\gamma^k)$ \\
\hline
$B_{42}(k)$ & $-(\omega_\pi(-1)+1)\omega_\pi(\gamma^k)$ & $0$ & $-(-1)^k((-1)^t+\chi(-1)\cdot$ & $-(-1)^k((-1)^t+\chi(-1)\cdot$ \\
& & & \qquad$\omega_\pi(-1))\chi(\gamma^k)\omega_\pi(\gamma^k)$ & \qquad$\omega_\pi(-1))\chi(\gamma^k)\omega_\pi(\gamma^k)$ \\
\hline
$B_{52}(k)$ & $-(\Lambda(\gamma^{1/2})+\Lambda(\gamma^{q/2}))\omega_\pi(\gamma^k)$ & $(\Lambda(\gamma^{1/2})+\Lambda(\gamma^{q/2}))\omega_\pi(\gamma^k)$ & $0$ & $0$ \\
\hline
$C_4(i,k)$ & $-\omega_\pi(\gamma^i)$ & $-\omega_\pi(\gamma^i)$ & $-2(-1)^{i+k}\chi(\gamma^k)\omega_\pi(\gamma^k)$ & $-2(-1)^{i+k}\chi(\gamma^k)\omega_\pi(\gamma^k)$ \\
\hline
$C_5(i,k)$ & $(q-1)\omega_\pi(\gamma^i)$ & $(q-1)\omega_\pi(\gamma^i)$ & $2(-1)^{i+k}(q-1)\chi(\gamma^k)\omega_\pi(\gamma^k)$ & $2(-1)^{i+k}(q-1)\chi(\gamma^k)\omega_\pi(\gamma^k)$ \\
\hline
$D_1(i,k)$ & $-(\Lambda(\theta^i)+\Lambda(\theta^{qi}))\omega_\pi(\gamma^k)$ & $-(\Lambda(\theta^i)+\Lambda(\theta^{qi}))\omega_\pi(\gamma^k)$ & $-((-1)^i+(-1)^k)(\Lambda(\theta^i)$ & $-2(1+(-1)^{i+k})\chi(\gamma^k)\omega_\pi(\gamma^k)$ \\
& & & \qquad$+\Lambda(\theta^{qi}))\chi(\gamma^k)\omega_\pi(\gamma^k)$ & \\
\hline
$D_2(i,k)$ & $(q-1)(\Lambda(\theta^i)+\Lambda(\theta^{qi}))\omega_\pi(\gamma^k)$ & $-(q-1)(\Lambda(\theta^i)+\Lambda(\theta^{qi}))\cdot$ & $0$ & $0$ \\
& & \qquad$\omega_\pi(\gamma^k)$ & & \\
\hline
$D_{32}(i,k)$ & $-(\Lambda(\gamma^{1/2}\eta^i)+\Lambda(\gamma^{q/2}\eta^{-i}))\cdot$ & $(\Lambda(\gamma^{1/2}\eta^i)+\Lambda(\gamma^{q/2}\eta^{-i}))\cdot$ & $0$ & $0$ \\
& \qquad$(1+\omega_\pi(-1))\omega_\pi(\gamma^k)$ & \qquad$(1+\omega_\pi(-1))\omega_\pi(\gamma^k)$ & & \\
\hline
$D_5(i,k)$ & $-(\Lambda(\theta^i)+\Lambda(\theta^{qi}))\omega_\pi(\gamma^k)$ & $(\Lambda(\theta^i)+\Lambda(\theta^{qi}))\omega_\pi(\gamma^k)$ & $0$ & $0$ \\
\hline
$D_6(i,k)$ & $(q-1-\Lambda(\eta^i)-\Lambda(\eta^{-i}))\cdot$ & $-(q-1+q(\Lambda(\eta^i)+\Lambda(\eta^{-i})))\cdot$ & $-(-1)^k(q+1)(\Lambda(\eta^i)$ & $-2(-1)^{i+k}\chi(\gamma^k)\omega_\pi(\gamma^k)$ \\
& \qquad$\omega_\pi(\gamma^k)$ & \qquad$\omega_\pi(\gamma^k)$ & \qquad$+\Lambda(\eta^{-i}))\chi(\gamma^k)\omega_\pi(\gamma^k)$ & \\
\hline
\newpage
$D_7(i,j,k)$ & $-(\Lambda(\eta^i)+\Lambda(\eta^{-i})$ & $-(\Lambda(\eta^i)+\Lambda(\eta^{-i})+\omega_\pi(-1)\cdot$& $0$ & $0$ \\
& \qquad$+\Lambda(\eta^j)+\Lambda(\eta^{-j}))\omega_\pi(\gamma^k)$ & $(\Lambda(\eta^j)+\Lambda(\eta^{-j})))\omega_\pi(\gamma^k)$ & & \\
\hline
$D_8(i,k)$ & $-(1+\Lambda(\eta^i)+\Lambda(\eta^{-i}))\omega_\pi(\gamma^k)$ & $\omega_\pi(\gamma^k)$ & $-(-1)^k(\Lambda(\eta^i)+\Lambda(\eta^{-i}))\cdot$ & $-2(-1)^{i+k}\chi(\gamma^k)\omega_\pi(\gamma^k)$ \\
& & & \qquad$\chi(\gamma^k)\omega_\pi(\gamma^k)$ & \\
\end{longtable}
\end{center}
\end{small}
\end{landscape}

\bibliographystyle{amsplain}

\end{document}